\pdfoutput=1
\documentclass[11pt,letterpaper]{article}
\usepackage[margin=1in]{geometry}
\usepackage{amsmath}
\usepackage{amsthm}
\usepackage{newtxtext}
\usepackage{newtxmath}
\usepackage{graphicx}
\usepackage{xcolor}
\usepackage{array}
\usepackage[labelfont=bf,labelsep=quad,font=small]{caption}
\usepackage{natbib}
\bibpunct[, ]{(}{)}{,}{a}{}{,}
\usepackage[colorlinks,allcolors=blue!60!black]{hyperref}
\hypersetup{pdftitle={Stochastic Leadtimes and Order Crossover: A Damped Response to a Weighted Inventory Position},pdfauthor={Tong Wang}}
\newcommand{\E}{\mathbb{E}}
\renewcommand{\P}{\mathbb{P}}
\newcommand{\Var}{\operatorname{Var}}
\newcommand{\cbar}{c}
\newcommand{\LB}{\mathrm{LB}}

\newcommand{\IP}{\mathrm{IP}}

\newcommand{\Lmin}{\underline{L}}
\newcommand{\Lmax}{\overline{L}}
\newcommand{\rng}{R}
\newcommand{\plag}{\hat p}
\newcommand{\Ipre}{\tilde I}

\makeatletter
\newcommand{\nv}[1]{\@ifundefined{nv@#1}{\textcolor{red}{[#1]}}{\@nameuse{nv@#1}}}

\makeatother
\input{sections/gen/numbers}

\newtheorem{proposition}{Proposition}
\newtheorem{lemma}{Lemma}
\newtheorem{corollary}{Corollary}

\makeatletter
\renewcommand{\proof}[1]{\par\addvspace{\medskipamount}\noindent\hspace*{1em}\textit{#1}\enskip\ignorespaces}
\renewcommand{\endproof}{\par\addvspace{\medskipamount}}
\newcommand{\Halmos}{\mbox{\quad$\square$}}
\newcommand{\up}{\rule{0pt}{11pt}}
\newcommand{\down}{\rule[-4pt]{0pt}{0pt}}
\newcommand{\TABLE}[3]{\centering\caption{#1}{\footnotesize #2}\par
  \def\public@note{#3}\ifx\public@note\@empty\else\vspace{3pt}%
  \begin{minipage}{\textwidth}\footnotesize\hspace*{1em}#3\end{minipage}\fi}
\newcommand{\FIGURE}[3]{\centering #1\par\caption{#2}%
  \def\public@note{#3}\ifx\public@note\@empty\else\par\vspace{4pt}\noindent
  \begin{minipage}{\textwidth}\footnotesize\textit{Note.} #3\end{minipage}\fi}
\renewcommand{\subsubsection}{\@startsection{subsubsection}{3}{\parindent}{1.5ex plus .5ex minus .2ex}{-1em}{\normalfont\normalsize\bfseries}}
\let\citepec\citep
\let\citetec\citet
\makeatother

\title{\textbf{Stochastic Leadtimes and Order Crossover:\\ A Damped Response to a Weighted Inventory Position}%
\thanks{The original idea of the policy structure proposed in this paper comes from another project, on using
reinforcement learning to solve complex dynamic programs; see auto-mdp-solver at
\protect\url{https://github.com/tong-wang/auto-mdp-solver}.}}
\author{Tong Wang\\[2pt] {\normalsize National University of Singapore}\\ {\normalsize\href{mailto:tong.wang@nus.edu.sg}{tong.wang@nus.edu.sg}}}
\date{November 2026}

\begin{document}
\maketitle

\begin{abstract}
\noindent Inventory control with stochastic leadtimes is a classical problem of
inventory theory. When orders can cross it remains a challenging one: the
base-stock policy need not be optimal, and the optimal policy is still
unknown. We offer two insights about crossover. Outstanding orders may
arrive late, so they should count for less than their size. An order placed
later may arrive early, so the order placed now should not close the whole
gap to its target. These two discounts give the policy its structure, a
damped response to a weighted inventory position. We propose two ways to
set them, and neither needs a search. The first is a closed form that uses
the leadtime distribution alone: an outstanding order counts with the
probability that it arrives no later than a new one. The second takes one
further step and uses the costs: an outstanding order counts by how often
it is missing when the stock is about to run out. Both policies perform very
well against the optimal cost, lower bounds that need no dynamic program,
and the existing policies. On \nv{bed.n} instances they
are \nv{lr.all.paper.mean}\% and \nv{lr.all.p2.mean}\% above the optimal
cost on average, against \nv{lr.all.bs.mean}\% for the base-stock policy.

\end{abstract}

\section{Introduction}\label{sec:intro}

The leadtime between placing an order and receiving it is at the heart of
inventory management. It determines how far ahead a firm must commit, how much
safety stock it needs, and how the orders that are still on their way should
be counted when the next order is placed. The classical theory answers these
questions for a leadtime that is fixed. In practice leadtimes are rarely
fixed, because ocean transit, port congestion, customs and the use of several
carriers all make the time between an order and its arrival random.

A random leadtime is not a difficulty by itself. As long as orders are
received in the sequence in which they were placed, most of the theory for
fixed leadtimes carries over \citep{kaplan1970,ehrhardt1984}. Every
outstanding order will have arrived by the time an order placed now arrives,
so all of them are worth the same as stock, the inventory position is a
sufficient statistic, and a base-stock policy on it is optimal. The
difficulty begins when order \emph{crossover} can happen, that is, when an
order placed today can arrive before one placed earlier. Orders can cross
whenever the leadtimes of successive orders are independent. An outstanding order may then arrive after the order placed now, in
which case it does not protect the periods that the new order is meant to
protect, so a unit on order is worth less than a unit in stock. The state can
no longer be compressed into the inventory position. It is the whole vector of
outstanding orders, and its dimension grows with the range of the leadtime.
Base-stock policies need not be optimal, and the optimal policy has not been
characterized. \citet[p.~410]{zipkin2000} calls base-stock policies under
independent leadtimes ``heuristic methods for what are, in fact, very complex
problems.''

Research on order crossover has taken two routes. The first keeps the
inventory position and the base-stock policy, and studies how that policy
performs when orders cross
\citep{zalkind1978,robinson2001,bradley2005,gongsun2026}. The second changes
the policy, and two streams of work change it in different ways. One changes
what the order responds to. It keeps the order-up-to rule but applies it to a
weighted inventory position, in which each outstanding order counts with a
weight instead of in full \citep{srinivasan2011,stolyarwang2022}. The other
changes how the order responds. It keeps the inventory position but orders
only a fraction of the gap between the position and its target
\citep{disney2016,wangdisney2017}. We call such a rule a damped response, and
the fraction its damping factor. Both changes lower the cost. In these
papers, however, the weights and the damping factor are found by search, by
tuning or by numerical optimization, and no rule says what they should be.

The two changes are not alternatives. Behind each is one insight about
crossover. An outstanding order may arrive after the order placed now, so it
should count for less than its size, and a weight discounts it. An order
placed later may arrive before the order placed now, so the order placed now
should not close the whole gap to its target, and damping holds it back.
Crossover has both effects at once, so a policy needs both discounts. This
gives the structure of our policy: a damped response to a weighted inventory
position. What remains is to set the weights and the damping factor, which
we call the coefficients of the policy.

We study the standard periodic-review model with backlogging, linear holding
and backorder costs, independent and identically distributed demands and
leadtimes, and the long-run average cost as the objective. The controller knows how old each
outstanding order is but not when it will arrive, which is what a firm's own
records of orders and receipts show. We ask three questions about a policy
with the structure above.

The first question is whether such a policy can be optimal when orders
cross. In general it cannot. We give a counterexample in which no weighted
inventory position with fixed weights reproduces the optimal orders. One set
of weights cannot serve both a deep backlog, where an old outstanding order
is worth a full unit of stock, and the normal range of the stock, where it is
worth less, so the weights that describe the optimal order depend on the
state. The policy is therefore a heuristic.

The second question is how to make the policy a heuristic that is
intuitive, easy to obtain and general. This is a matter of setting its
coefficients, and we set them in two ways. The first is a closed form that uses the leadtime distribution alone. It is intuitive: the
weight of an outstanding order is the probability, given its age, that it
arrives no later than an order placed now, and the damping factor falls from
one toward one half as later orders overtake more often. It is easy to
obtain: each weight and the damping factor have a formula, the level is the
best base-stock level of the classical theory lowered by an explicit amount,
and nothing is searched or simulated and no dynamic program is solved. It is
general: the same formulas apply to every leadtime distribution, and the
same weights and damping factor serve every demand distribution and every
cost ratio. Without crossover the policy is the base-stock policy. With
exponential leadtimes in continuous time every weight is one half, and the
weighted position is that of the generalized base-stock policy of
\citet{stolyarwang2022} at one value of its parameter.

The closed form has a price. Its weights ask only whether an outstanding
order will arrive before the new one. But an order that arrives after the new one,
a late order, costs something only if the stock runs out while it is missing, and when the stock runs low
depends on the cost ratio. With a high cost ratio the safety stock is large.
The stock then runs low only at the end of the periods that an order has to
cover, and by then a late order has usually arrived, so it should count
nearly in full. With a low cost ratio the stock runs low early, while a late
order is still missing, and it should count for less. Our second policy
is cost aware. It counts each outstanding order by how often it is missing
at the moments when the stock is about to run out. It takes its weights and
its damping factor from these moments, and it sets its level for the spread
of the net inventory that they leave. Like the optimal
policy, it replaces less of a late order as the cost ratio rises, and it
treats a late order that will arrive soon differently from one that will
arrive later. The moments depend on the policy itself, so the coefficients
come from a small computation that reproduces them, with no search. So the
closed-form policy gives up some performance for coefficients that can be
written down, and the cost-aware policy gives up the formula for most of the
performance that is left.

The third question is how to guarantee the performance of such a policy
when the optimal policy cannot be computed. We develop lower bounds on the
cost of every policy that need no dynamic program, so the gap between a
policy and a bound is a guarantee that does not require the optimal policy.
As the shortest leadtime grows, this gap vanishes at least as fast as one
over the shortest leadtime for a rule with fixed coefficients at its best
level. On a test bed of \nv{bed.n} instances for which the optimal cost can be
computed, the best bound is on average \nv{lr.all.lbk.mean}\% below the
optimal cost. The closed-form policy is \nv{lr.all.paper.mean}\% above the
optimal cost on average and the cost-aware policy \nv{lr.all.p2.mean}\%,
against \nv{lr.all.bs.mean}\% for the base-stock policy. The two
policies also compare well with the published policies.

The paper is organized as follows. Section~\ref{sec:lit} reviews the
literature. Section~\ref{sec:model} sets up the model and its optimality
equation. Section~\ref{sec:structure} studies the optimal policy and shows
that a policy with our structure is not optimal in general.
Section~\ref{sec:policy} gives the closed-form policy,
Section~\ref{sec:costaware} the cost-aware policy, Section~\ref{sec:bounds}
the lower bounds and Section~\ref{sec:numerical} the numerical study.
Section~\ref{sec:conclusion} concludes. Proofs are in the e-companion.

\section{Literature}\label{sec:lit}

This section first reviews the work on leadtimes without crossover and on
information about leadtimes. It then reviews the work on order crossover:
the base-stock policy, and the two streams that change the policy.

\subsubsection*{Stochastic leadtimes without crossover.} \citet{kaplan1970} shows that
when leadtimes are random but orders cannot cross, the problem reduces to one
with a single state variable, and \citet{ehrhardt1984} extends the reduction to
fixed ordering costs. \citet{zipkin1986} and \citet{songzipkin1996} model supply
systems that generate such leadtimes, and \citet{song1994} studies how
leadtime uncertainty moves the base-stock level.

\subsubsection*{Information about leadtimes.} \citet[p.~410]{zipkin2000}
observes that an optimality theory for independent leadtimes needs a scenario
for when the leadtimes are observed, for instance when they begin or when
they end. We take the second, as \citet{bai2026} do: a leadtime is observed
when the order is received. The policies that have been proposed for
crossing orders, which we review below, also use the number of outstanding
orders or their ages, and none uses arrival dates. Other information has been
studied elsewhere: the state of a supply system that processes orders in
sequence \citep{songzipkin1996}, the leadtime of each order before it is
placed, without crossover \citep{chenyu2005}, the stage that each order has
reached in dual sourcing \citep{songxiaozhang2025}, and leadtimes that vary
but are known in advance \citep{riezebos2009}. We are not aware of a study of
the optimal policy under crossover with either unknown or known arrival
dates.

\subsubsection*{Base-stock policies under crossover.} A long line of work keeps the
order-up-to policy and analyzes it when orders cross
\citep{riezebos2006,hayya2008,bischak2014}. \citet{zalkind1978}
obtains the distribution of the number of outstanding orders and the optimal
order-up-to level, and notes by an example that a policy that looks at the
outstanding orders can do better. \citet{robinson2001} show that the level
must be set from the distribution of the shortfall and not from leadtime
demand, and \citet{bradley2005} give approximations. The best base-stock level of this
literature is one input of our policy (Section~\ref{sec:level}).
\citet{gongsun2026} bound the optimality gap of the best base-stock policy and
study its asymptotics, and their lower bound is the starting point of our
Section~\ref{sec:bounds}, where we show how to sharpen it.

\subsubsection*{Policies that change the position.} The first stream keeps
the order-up-to rule and changes what the order responds to.
\citet{srinivasan2011}
compute the optimal policy of a small lost-sales model with crossover by
dynamic programming, find that it has no simple form, and propose to order up
to a fixed level on a weighted sum of on-hand stock and each outstanding
order. Their weights come from a grid search by simulation with the level
held at the best base-stock level, and they come out close to one except
under very skewed leadtime distributions. Their rule also puts a weight on
the stock. A stock weight below one closes only a part of the gap, as a
damping factor does, but with the level held fixed it also raises the
target, so the two cannot be set separately. \citet{srinivasan2007} explains the weights found
in words, by how soon an order is likely to arrive, and gives no formula for
them. \citet{stolyarwang2022} study a continuous-time model with exponential
leadtimes and propose the generalized base-stock policy, which is a base-stock
policy on the net inventory plus the pipeline divided by a parameter
$\gamma>1$. They tune $\gamma$ and show large savings.

\subsubsection*{Policies that change the response.} The second stream keeps
the inventory position and changes how the order responds to it.
\citet{disney2016} give
exact safety stocks for a proportional order-up-to policy, which closes only a
fraction of the gap between the inventory position and its target, and find
numerically that the best fraction is below one whenever orders cross.
\citet{wangdisney2017} prove, for independent demands, that the order-up-to
policy never minimizes the variance of inventory within this class when
orders cross. In both papers the gap is measured on the unweighted inventory
position, and the best fraction is computed numerically.
\citet{weiyang2026}, in continuous time with exponential leadtimes, propose a
policy with two base-stock levels, a high one when the net inventory is low
and a low one when it is high, which uses neither ages nor weights.

\citet{stolyar2026} replace the rule by a computation. They use the ages of
the units in transit, which is our information setting, in a model in which
each unit has its own leadtime. Their policy solves, in every period, a linear
program for the next orders, in which demand is replaced by its mean and the
inventory projected for each future period counts a unit of a given age with
the probability that it arrives by that period. Our weight for an order of a
given age is the average of these probabilities over the leadtime of a new
order, so the ingredients are related. In their policy, however, the
probabilities enter a program that is solved in every period, and neither a
single position nor a closed form for the order results. They name the case
in which the units of an order share one leadtime, which is ours, as an
extension of their model.

Our paper joins the two streams. Each of them is built on one of the two discounts
that crossover calls for, and our policy makes both: a damped response
to a weighted inventory position. Where these papers search, tune or optimize
numerically, we set the weights and the damping factor without a search,
first in closed form from the leadtime distribution alone and then, one step
further, with the costs.

\section{Model}\label{sec:model}

\subsection{Dynamics and Costs}

Time is discrete, with periods $t=1,2,\dots$ In each
period the controller first places an
order $q_t\ge0$, then the orders due in the period are received, and then the
demand $D_t$ occurs, with unmet demand backlogged. Demands are independent and
identically distributed and nonnegative, with mean $\mu>0$ and variance
$\sigma^2$, and $D(j)$ denotes the sum of $j$ independent demands. A
nonnegative order is placed in every period, and each order has its own
leadtime: the order of period $t$ has leadtime $L_t$ and is received in period
$t+L_t$. Leadtimes are independent of each other and of the demands, and they
have a common distribution with bounded support, as in \citet{zalkind1978},
\citet{disney2016} and \citet{bai2026}. Let $\Lmin\ge0$ be the shortest and
$\Lmax$ the longest leadtime in the support, and write
\[
p_j=\P(L_t=j),\qquad j=\Lmin,\dots,\Lmax .
\]
The probabilities $p_{\Lmin}$ and $p_{\Lmax}$ are positive by definition, and
those in between may be zero. We call $\rng=\Lmax-\Lmin$ the range of the
leadtime. We consider $\rng\ge2$. With $\rng=0$ the leadtime is
fixed, and with $\rng=1$ no two orders cross, and in both cases the problem
is the classical one.
The leadtime of an order is observed when it ends, that is, when the order is
received.

Let $I_t$ be the net inventory at the end of period $t$. The cost of the period
is $\cbar(I_t)=h\,I_t^++b\,I_t^-$ with $h,b>0$, and $\beta=b/(b+h)$ is the
critical ratio. There is no ordering cost. The
objective is to minimize the long-run average cost. For a policy $\pi$ this
cost, $g^\pi$, is the limit superior as $T\to\infty$ of the expected cost of
the first $T$ periods divided by $T$, and $g^*$ is the infimum of $g^\pi$
over all policies. Orders and demands are either all real numbers or all
integers, and we say so where the difference matters.

\subsection{The State}

At the beginning of period $t$ the controller knows the net inventory, which
is $I_{t-1}$, and the orders that are outstanding. For $k\ge1$ let $a_{k,t}$
be the order placed $k$ periods ago if it is still outstanding, and zero
otherwise. We call it the outstanding order of age $k$. An order can be
outstanding up to age $\Lmax$, so the natural state is the net inventory
together with $a_{1,t},\dots,a_{\Lmax,t}$. Given the information of the
controller, the remaining leadtimes of the outstanding orders are
independent.

\begin{figure}
\FIGURE
{\includegraphics[width=0.68\textwidth]{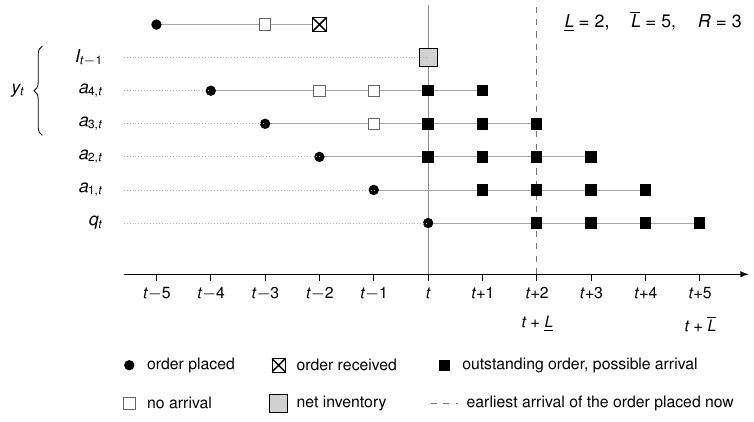}}
{The State at the Beginning of Period $t$ for a Leadtime of Two to Five Periods\label{fig:state}}
{The order of period $t-5$ was received in period $t-2$ and is part of the net
inventory $I_{t-1}$. The orders of the four periods before $t$ are
outstanding. Those of ages four and three arrive no later than period
$t+\Lmin$, the earliest period in which the order placed now can arrive, and
are pooled with the net inventory in the stock
$y_t=I_{t-1}+a_{3,t}+a_{4,t}$. Those of ages two and one can arrive later
than the order placed now and are kept separately, as $a_{2,t}$ and
$a_{1,t}$.}
\end{figure}

The outstanding orders are of two kinds (Figure~\ref{fig:state}). An order
placed in period $t$ arrives in period $t+\Lmin$ at the earliest. An outstanding order
of an age $k\ge\rng$ arrives in period $t+\Lmin$ at the latest, because no
leadtime exceeds $\Lmax=\Lmin+\rng$. So it is received no later than any
order placed from now on, and for the purpose of ordering it is as good as
on-hand stock. An outstanding order of an age $k<\rng$ may arrive before period
$t+\Lmin$ or after it, so an order placed now or later can overtake it.

The orders of the first kind therefore count in full. We say that such an
order is \emph{pooled}: it is added to the net inventory, and its size and
its age are no longer recorded. The net inventory and the pooled orders
together are the stock,
\[
y_t=I_{t-1}+\sum_{k=\rng}^{\Lmax}a_{k,t}.
\]
So an outstanding order is pooled when it reaches age $\rng$, whether or not
it has been received by then. The orders of the second kind have to be kept
one by one. So the state is $(y_t,a_t)$, where $a_t=(a_{1,t},\dots,a_{\rng-1,t})$ is
the vector of the outstanding orders of the ages below $\rng$, and it has
dimension $\rng$.
We drop the index $t$ when the period is clear.

The same information can be written with the inventory position. The
inventory position before ordering is $\IP_t=y_t+\sum_{k<\rng}a_{k,t}$, the
net inventory plus everything that is outstanding, and the inventory
position after ordering is $z_t=\IP_t+q_t$. The pair $(\IP_t,a_t)$ carries
the same information as $(y_t,a_t)$, and choosing an order $q_t\ge0$ is the
same as choosing a position $z_t\ge\IP_t$. The optimality equation below,
the structural results of Section~\ref{sec:structure} and the two policies
are written in $(y,a)$ and the order $q$, because they compare the stock
with the outstanding orders.

\subsection{The Optimality Equation}\label{sec:bellman}

The optimal policy is characterized by an optimality equation. To write it
we need the transition of the state and the cost that a decision controls.
Both
depend on when the outstanding orders will be received, and what the
controller knows about this is contained in one family of probabilities.

An outstanding order of age $k$ with leadtime $L$ is received $L-k$ periods
from now, where zero periods from now is the current period. For $k\ge0$ and
$j\ge0$ let
\begin{equation}\label{eq:B}
B_{k,j}=\mathbf 1\{L\ge k+j+1\}
\end{equation}
indicate the event that the order is received neither in the current period nor in the
next $j$ periods. Age zero is the order placed in the current period. The
controller does not know the leadtime of an outstanding order. All that it knows is that $L\ge k$. So for the controller $B_{k,j}$ is
a Bernoulli random variable, and we write
\begin{equation}\label{eq:m}
m_{k,j}=\P(B_{k,j}=1\mid L\ge k)=\P(L\ge k+j+1\mid L\ge k)
\end{equation}
for the probability that it equals one. Given the information of the
controller, the variables $B_{k,j}$ of different orders are independent of
each other and of the demands, because the leadtimes are. These probabilities
carry the whole effect of crossover. An order placed now with leadtime $j$ is
received $j$ periods from now, so $m_{k,j}$ is the probability that it
overtakes the outstanding order of age $k$.

The two ends of the support fix the boundary values. If $k+j<\Lmin$, the
order is too young to be received within that time, because no leadtime is
shorter than $\Lmin$, so $B_{k,j}=1$ for certain and $m_{k,j}=1$. If
$k+j\ge\Lmax$, it has been received for certain, because no leadtime is
longer than $\Lmax$, so $B_{k,j}=0$ and $m_{k,j}=0$. In between,
$m_{k,j}=\P(L\ge k+j+1)/\P(L\ge k)$ is positive, and it is below one unless
the leadtime takes no value from $k$ to $k+j$. For example,
$m_{1,0}=\P(L\ge2\mid L\ge1)$ is the probability that the order of the
previous period, if it is outstanding, is not received in the current one. The order placed in the current period is not received in it with
probability $m_{0,0}=\P(L\ge1)$.

\subsubsection*{The transition.} Let a prime mark the next period, and let
$D$ be the demand of the current period. The state moves by
\begin{equation}\label{eq:transy}
\begin{aligned}
y'&=y+a_{\rng-1}+\sum_{k=1}^{\rng-2}a_k\big(1-B_{k,0}\big)+q\big(1-B_{0,0}\big)-D,\\
a_1'&=q\,B_{0,0},\qquad a_{k+1}'=a_k\,B_{k,0},\quad k=1,\dots,\rng-2 .
\end{aligned}
\end{equation}
The stock loses the demand. It gains the order of age $\rng-1$, which is
pooled from now on whether or not it has been received, so that
$B_{\rng-1,0}$ does not enter. It also gains each younger outstanding order
that is received in the period, and the order $q$ just placed if that is
received in the period, which is the event $B_{0,0}=0$ and is possible only
with a shortest leadtime of zero. An outstanding order of age $k<\rng-1$
that is not received stays in the vector with age $k+1$, which is the event
$B_{k,0}=1$ and has probability $m_{k,0}$, and the order just placed enters
the vector with age one if it is not received.

\subsubsection*{The cost that a decision controls.} The order of period $t$ arrives in
period $t+\Lmin$ at the earliest, so it cannot affect the cost of an earlier
period. As with a fixed leadtime, we therefore charge to the decision of
period $t$ the cost of period $t+\Lmin$, the first period that it can
affect. By the end of
that period the stock $y$ has been received in full. The outstanding order of
age $k$ has been received unless $B_{k,\Lmin}=1$, which has probability
$m_{k,\Lmin}$.
The order $q$ of period $t$ has been received if its leadtime $L_t$ is the
shortest, and no later order has. The demand of the $\Lmin+1$ periods from
$t$ to $t+\Lmin$ is $D(\Lmin+1)$. Let $\Ipre_{t+\Lmin}$ be the net inventory at
the end of period $t+\Lmin$ before this demand is subtracted, which is what
has been received:
\begin{equation}\label{eq:received}
I_{t+\Lmin}=\Ipre_{t+\Lmin}-D(\Lmin+1),\qquad
\Ipre_{t+\Lmin}=y+\sum_{k=1}^{\rng-1}a_k\big(1-B_{k,\Lmin}\big)+q\,\mathbf 1\{L_t=\Lmin\}.
\end{equation}
The two parts are independent, because the state and the order of period $t$
are set before any of this demand occurs, and the leadtimes do not depend on
the demands. The demand $D(\Lmin+1)$ is what every order has to wait out even
with the shortest leadtime, and no policy can touch it, whereas
$\Ipre_{t+\Lmin}$ collects what the policy and the leadtimes contribute.
Averaging the cost of the period over this demand first, with $\Ipre$ held
fixed, gives a smoothed cost $G$ as a function of $\Ipre$:
\begin{equation}\label{eq:Gamma}
G(\Ipre)=\E\big[\cbar(I)\mid\Ipre\big]=\E\big[\cbar\big(\Ipre-D(\Lmin+1)\big)\mid\Ipre\big].
\end{equation}
The demand of $\Lmin+1$ periods has the same distribution in every period, so
$G$ is the same function in every period, and we have dropped the index of
the period.

\subsubsection*{The equation.} Charging the cost of period $t+\Lmin$ to the
decision of period $t$ shifts the costs by $\Lmin$ periods and does not
change their long-run average. So the problem is a Markov decision problem
with the state $(y,a)$, the order $q\ge0$ as the action, the transition
\eqref{eq:transy} and the stage cost $\E\,G(\Ipre_{t+\Lmin})$.

\begin{proposition}[Optimality equation]\label{prop:acoe}
Let quantities be all real numbers or all integers. There is a function
$V^*\ge0$ of the state such that

(a) at every state,
\begin{equation}\label{eq:acoe}
g^*+V^*(y,a)=\min_{q\ge0}\big\{\E\,G\big(\Ipre_{t+\Lmin}\big)+\E\,V^*(y',a')\big\},
\end{equation}
where $(y',a')$ is the next state \eqref{eq:transy};

(b) $g^*$ is the optimal long-run average cost from every initial state, and
a stationary policy that attains the minimum in (a) at every state is
optimal.
\end{proposition}

We call \eqref{eq:acoe} the optimality equation, $V^*$ the relative value
function and a minimizer of its right side an optimal order. Where several
orders attain the minimum, $q^*(y,a)$ denotes the smallest of them, and
every statement about $q^*$ below also holds for the largest.
$V^*(y,a)$ measures the cost from the state $(y,a)$ relative to $g^*$ per
period.
By
part (b) the optimal cost does not depend on the initial state, and a
stationary policy attains it. The proof, in Section~\ref{sec:ec-acoe} of
the e-companion, constructs $V^*$ as a limit of discounted value functions.
When orders cannot cross, that is, when the range is zero or one, the equation
has the classical solution: the stock $y$ is the inventory position, the
vector $a$ is empty, and a base-stock policy is optimal \citep{kaplan1970}.

\section{The Optimal Policy}\label{sec:structure}

The two discounts of the introduction give a policy of a simple kind. It
orders a fixed fraction of the gap between a target level and a weighted
inventory position, in which each outstanding order counts with a fixed
weight.
Section~\ref{sec:policy} writes it out. Before we set its coefficients, this
section asks whether such a policy can be optimal. In general it cannot. At
range two the optimal order has both discounts in a local form, but they
depend on the state, and beyond range two only a part of this survives.

One property holds at every range. With real-valued orders the relative
value function $V^*$ is convex in $(y,a)$ (Section~\ref{sec:ec-structure} of
the e-companion). Convexity gives a first-order condition for the order, but
it does not say how the optimal order moves with the stock or with an
outstanding order. More can be said at range two, the smallest range at
which orders can cross.

\begin{proposition}[Range two]\label{prop:chain}
Let $\rng=2$. The relative value function $V^*$ is L$^\natural$-convex in
$(a,\,y+a)$, the outstanding order and the inventory position, in the discrete
sense when quantities are integers. At every state, the optimal order
satisfies
\begin{equation}\label{eq:chain}
q^*(y,a)-1\le q^*(y+1,a)\le q^*(y,a+1)\le q^*(y,a),
\end{equation}
and the middle inequality holds with equality if $y+1+q^*(y,a+1)\le0$.
\end{proposition}

The chain \eqref{eq:chain} has a direct reading. Write
$\kappa^*(y,a)=q^*(y,a)-q^*(y+1,a)$ for the local response of the order to a
unit of stock, and $\lambda^*(y,a)=q^*(y,a)-q^*(y,a+1)$ for its response to a
unit on order. Then
\[
0\le\lambda^*\le\kappa^*\le1 .
\]
Where $\kappa^*>0$, the ratio $\theta^*=\lambda^*/\kappa^*$ lies between zero and
one. It is the local worth of a unit on order relative to a unit of stock.
So the optimal order closes at most the whole of a shortfall, and it counts
the outstanding order at no more than its full size. These are the two
discounts in a local form, with $\kappa^*$ in the place of the damping
factor and $\theta^*$ in the place of the weight. In the base-stock policy
both responses are one.

The two responses are not constant, and the equality case of the proposition
shows why. Its condition describes a deep backlog: the stock is so far below
zero that the new order does not bring it above zero. Every unit is then
used the moment it arrives. So a unit ordered now is needed for certain,
whether or not the outstanding order has arrived, and a unit of the
outstanding order is a perfect substitute for a unit of stock. Away from a
deep backlog this is no longer so. The demand then decides whether there
will be stock left over or a backlog, a unit ordered now is more likely to
be needed if the outstanding order has not arrived, and stock can count for
more than the unit on order. So the worth of the outstanding order, relative
to stock, is one in a deep backlog and can be less elsewhere: it depends on
the state.

The proof, in Section~\ref{sec:ec-structure} of the e-companion, relies on
the fact that the state has only one outstanding order when the range is
two. Beyond range two the chain itself fails. Its second inequality says
that a unit of stock lowers the optimal order by at least as much as a unit
of an outstanding order does. In instances of Section~\ref{sec:numerical}
with ranges three and four, this fails for the order of age one at some
states that the optimal policy visits (Section~\ref{sec:ec-structure} of the
e-companion). It survives at every range only for the oldest outstanding
order, with its equality case (Proposition~\ref{prop:oldest} of the
e-companion).

A policy with fixed coefficients is therefore not optimal in general. Its
order is a function of one weighted position, and constant weights cannot
serve both a deep backlog and the normal range of the stock.
Section~\ref{sec:ec-wp} of the e-companion gives a counterexample over a
finite horizon that can be checked by hand, and a numerical check in the
long run. In both, no weighted position with fixed weights reproduces the
optimal orders. So our policy is a heuristic. The next two sections set its
coefficients, and Section~\ref{sec:bounds} gives lower bounds that say how
well it does.

\section{The Closed-Form Policy}\label{sec:policy}

The policy makes the two discounts with constant coefficients. The first
discount is for the outstanding orders, which may arrive after the order
placed now. We give the outstanding order of age $k$ a constant weight
$\theta_k\ge0$, for $k=1,2,\dots,\rng-1$, and call the sum
\begin{equation}\label{eq:wp}
P_t=y_t+\sum_{k=1}^{\rng-1}\theta_k\,a_{k,t}
\end{equation}
a \emph{weighted} inventory position. With every weight equal to one it is
the inventory position. The second discount is for the order placed now.
With crossover it is not the only order that can close the gap between the
target level and the position by the time it arrives, and it may itself
arrive late. So the policy closes a fraction $\kappa$ of the gap now and
leaves the rest to the orders that follow. We call such a rule a
\emph{damped} response and the constant $\kappa$ its damping factor. With
both discounts the policy orders
\begin{equation}\label{eq:policy}
q_t=\kappa\,(S-P_t)^+ ,
\end{equation}
a damped response to a weighted inventory position with a target level $S$.
Under this rule the order falls by $\kappa$ for each unit of stock and by
$\kappa\theta_k$ for each unit of the outstanding order of age $k$. So
$\kappa$ and $\theta_k$ are the constant counterparts of the responses
$\kappa^*$ and $\theta^*=\lambda^*/\kappa^*$ of the optimal order in
Proposition~\ref{prop:chain}, which lie between zero and one at range two
and depend on the state.

This section sets the weights, the damping factor and the target level in
closed form, and the next three subsections give them in turn. The weights
and the damping factor are computed from the leadtime distribution alone,
and the target level also uses the demand distribution and the two cost
rates. Section~\ref{sec:costaware} takes one further step and lets the
weights and the damping factor use the costs as well, at the price of a
small computation.

\subsection{The Weights}

The proposed weight of an outstanding order is the probability that it is received no
later than an order placed now. Let the outstanding order of age $k$ have
leadtime $L$, of which the controller knows only that $L\ge k$, and let $L'$
be the leadtime of an order placed now, which is independent of $L$. The
outstanding order is received $L-k$ periods from now, and the new order $L'$
periods from now. So we set the weight of age $k\ge1$ to
\begin{equation}\label{eq:rk}
\theta_k=\P\big(L-k\le L'\mid L\ge k\big)=1-\sum_jp_j\,m_{k,j},
\end{equation}
where the second form uses that $m_{k,j}$ of \eqref{eq:m} is the probability
that an order placed now with leadtime $j$ overtakes the order of age $k$.

\begin{proposition}[The weights]\label{prop:weights}
(a) With the weights \eqref{eq:rk}, the weighted position $P_t$ is, given
the state, the stock that is expected to have been received by the time an
order placed now is received.

(b) $\theta_k=1$ for $k\ge \rng$. At range two,
$\theta_1=1-p_{\Lmin}p_{\Lmax}/\P(L\ge1)$, which is $1-p_{\Lmin}p_{\Lmax}$ when the
shortest leadtime is at least one period.

(c) If the leadtime has a nondecreasing hazard rate, $\theta_k$ is nondecreasing in
$k$.

(d) Let $\chi$ be the expected number of later orders that arrive strictly
before a given order. Then
\begin{equation}\label{eq:chi}
\chi=\E\,(L-L'-1)^+=\sum_{k\ge1}(1-\theta_k)\,\P(L\ge k) .
\end{equation}
\end{proposition}

Part (a) of Proposition~\ref{prop:weights} is the reasoning behind the
inventory position itself. With a fixed leadtime everything outstanding has
arrived by the time the new order arrives, so it counts in full, and with
crossover it has arrived by then only with the probability $\theta_k$.
Part (b) says that an order old enough
to arrive for certain before a new one counts in full, which is why such orders
were pooled with the stock in Section~\ref{sec:model}. Part (c) says that under
the usual shapes of leadtime distributions an older order counts for more.
\citet{srinivasan2011} state this direction as an expectation for their
searched weights.

Part (d) connects the weights to how often orders cross. Take the order
placed now and the order placed $k$ periods before it. The new order arrives
strictly before the earlier one with probability
$(1-\theta_k)\,\P(L\ge k)$: the earlier order must still be outstanding, and
the new order must then overtake it. Summing over $k$ gives the expected
number of earlier orders that a new order overtakes, which is also the
expected number of later orders that overtake a given order. We call this
number $\chi$ the crossover number. It is zero if and only if orders cannot
cross, and at range two it equals $p_{\Lmin}p_{\Lmax}$. The overtaken orders
are what the weights remove from the position. Their expected number is
$\chi$ and each is of size $\mu$ on average, so the weighted position is
below the inventory position by $\mu\chi$ on average.

Some distributions give simple weights. If the leadtime is uniform on
$\{1,2,3,4\}$, the two weights are $13/16$ and $11/12$. Two other families
lie outside the model of Section~\ref{sec:model}, which has a bounded
support in discrete time, and they only illustrate the formula. If the leadtime is geometric on $\{1,2,\dots\}$ with success
probability $p$, every age has the same weight $(1+p-p^2)/(2-p)$, because the
remaining leadtime does not depend on the age. In continuous time with
exponential leadtimes the weight is one half at every age. The weighted
position is then the net inventory plus half of the pipeline, which is the
generalized base-stock policy of \citet{stolyarwang2022} with their parameter
$\gamma$ equal to two. All of these follow from \eqref{eq:rk} by direct
calculation.

\subsection{The Damping Factor}\label{sec:damping}

When a gap has to be closed, how much of it should be ordered now, and how
much should be kept for the next period? The damping factor is the share
that is ordered now. Take a gap of one unit. Let the order placed now close
the fraction $\kappa$ of it and the next order the rest, with all other
orders unchanged, and let $L$ and $L'$ be the leadtimes of the two orders.
Before either part is received the stock is short by the whole unit, and
after both are received it is not short, whatever $\kappa$ is. In between,
it is short by $1-\kappa$ for $(1+L'-L)^+$ periods if the order placed now
is received first, and by $\kappa$ for $(L-L'-1)^+$ periods if the next
order is received first. These lengths have the means $1+\chi$ and $\chi$,
by \eqref{eq:chi} and $\E(1+L'-L)=1$. If every period is charged the square
of its shortfall, the expected charge is
$(1+\chi)(1-\kappa)^2+\chi\,\kappa^2$ plus a term without $\kappa$. It is
smallest at
\begin{equation}\label{eq:kappa}
\kappa=\frac{1+\chi}{1+2\chi},
\end{equation}
and we use this share as the damping factor.

The order placed now is one period ahead of the next order, so it closes
more than half of the gap. It does not close the whole gap, because the
next order may still be received first, and the part ordered now is then
missing in the meantime. Without crossover $\chi=0$ and the whole gap is
closed now, and as the crossover number grows the share decreases to one
half. Under the uniform distribution it is
$(\rng+2)(\rng+3)/\big(2(\rng^2+2\rng+3)\big)$, which is $10/11$ at range two
and $5/6$ at range three. The shortfall is squared because the expected
cost of a period is convex in the stock, and near its minimum it grows with
the square of a deviation. The calculation counts every period alike and
holds the other orders fixed, and Section~\ref{sec:costaware} relaxes both.

\subsection{The Target Level}\label{sec:level}

Under a base-stock policy with level $S$ on the inventory position, the
position after ordering is always $S$, and the order of each period is the
demand of the period before. The net inventory at the end of a period is
therefore $S$ less the demand of the period and less the orders that are
still outstanding, each of which equals an earlier demand. An order placed $j$
periods earlier is still outstanding with probability $\P(L\ge j+1)$, which is
one for $j<\Lmin$. So the net inventory is $S-Z$ with the shortfall
\[
Z=D(\Lmin+1)+\sum_{j=\Lmin}^{\Lmax-1}D_j'\,\tilde B_j,
\]
where the $D_j'$ are independent demands and the $\tilde B_j$ are independent
Bernoulli variables with means $\P(L\ge j+1)$. With a fixed leadtime $L$ the
shortfall is the leadtime demand $D(L+1)$. With crossover the two differ. The
shortfall has the same mean as the demand over a leadtime and a smaller
variance, because the indicators $\tilde B_j$ are independent, whereas the
corresponding indicators of a single leadtime move together. This is why the
level has to be set from the shortfall and not from the leadtime demand
\citep{robinson2001}. The level with the smallest
long-run average cost is the classical one
\citep{zalkind1978,robinson2001}: $S^{\IP}$ is the $\beta$-quantile of $Z$,
where $\beta=b/(b+h)$ is the critical ratio of Section~\ref{sec:model}. This
level is for the inventory position, and our policy uses the weighted
position, which is lower by $\mu\chi$ on average
(Proposition~\ref{prop:weights}(d)). The first correction lowers the level
by this amount.

The second correction is for the damping. An order-up-to rule brings the
position to its level in every period. A rule that closes only the fraction
$\kappa$ of the gap between its target and the position does not, and it
lets the position sit lower. The order has to replace a demand of $\mu$ per
period on average, so the gap before ordering has to be $\mu/\kappa$ on
average, where it is $\mu$ under an order-up-to rule. To keep the position
where the order-up-to rule has it, we raise the target by the difference,
$(1-\kappa)\mu/\kappa$. The target level of our policy is therefore
\begin{equation}\label{eq:S}
S=S^{\IP}-\mu\chi+\frac{1-\kappa}{\kappa}\,\mu .
\end{equation}

Computing the target level takes the mean demand and one quantile of the
shortfall distribution, which is a convolution. With integer demand we
interpolate the quantile between the integers (Section~\ref{sec:ec-wli} of
the e-companion).

With the weights \eqref{eq:rk}, the damping factor \eqref{eq:kappa} and the
target level \eqref{eq:S}, the rule \eqref{eq:policy} is the closed-form
policy. Without crossover every weight and the damping factor are one, and
it is the base-stock policy with the best level.

\section{The Cost-Aware Policy}\label{sec:costaware}

The cost-aware policy lets the weights and the damping factor depend on the
cost ratio $b/h$:
\begin{equation}\label{eq:policy2}
q_t=\kappa^\beta\Big(S^\beta-y_t-\sum_{k=1}^{\rng-1}\theta^\beta_k\,a_{k,t}\Big)^+ ,
\end{equation}
where the superscript marks a damping factor, weights and a level that
depend on the cost ratio, through the critical ratio $\beta=b/(b+h)$.

Two things differ from Section~\ref{sec:policy}. First, a late arrival alone
costs nothing. A late order costs something only if the stock runs out while
it is missing, and how often that happens depends on the cost ratio. So the
discount of an outstanding order is the probability that it is missing when
the stock just runs out, where it was the probability that it arrives after
a new order. Second, this probability depends on the coefficients
themselves, so they are a fixed point, where they were a closed form.

\subsection{The Cost Charged to an Order}\label{sec:identity}

The cost depends on the net inventory of each period, so we have to say on
which order the net inventory of a period depends. With a fixed leadtime $L$
it is the order of period $t-L$: the net inventory at the end of period $t$
is the inventory position after that order, less the demand since then, and
the classical theory charges the cost of period $t$ to it. We do the same
when leadtimes are random.

\begin{figure}
\FIGURE
{\includegraphics[width=0.5\textwidth]{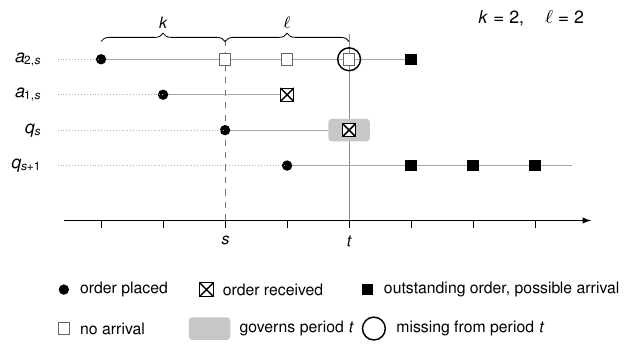}}
{The Order That Governs a Period, and an Order That Is Missing from It\label{fig:governs}}
{The leadtime is two to five periods, and the figure is drawn as of period
$t=s+2$, looking back. The order of period $s$ is the youngest order received
by period $t$, so it governs period $t$, with lag $\ell=2$. The order of age
two at $s$ has not arrived by period $t$, so it is missing from that period.
The order of period $s+1$ cannot arrive by period $t$.}
\end{figure}

We need two terms (Figure~\ref{fig:governs}). Consider the end of a period
$t$, and let $s$ be the period of the order that was placed last among the
orders received by then. We say that the order of period $s$ \emph{governs}
period $t$, and we call $\ell=t-s$ the \emph{lag} of period $t$. With a
fixed leadtime the lag is the leadtime. Here it lies between $\Lmin$ and
$\Lmax$. It is determined by the leadtimes alone, so it does not depend on
the policy. Every order placed after period $s$ is still outstanding at the
end of period $t$. An order placed before period $s$ may be outstanding as
well, and we then say that it is missing from period $t$. We first derive the
distribution of the lag and then the cost of a period given its lag.

Now take the view of period $s$, in which the order $q_s$ is chosen. This
order governs the period $s+\ell$ with probability
\[
\plag_\ell=\P(L\le\ell)\prod_{i=0}^{\ell-1}\P(L>i),\qquad \ell=\Lmin,\dots,\Lmax,
\]
because it must be received within $\ell$ periods, and the order that is
placed $i$ periods before period $s+\ell$ must not be received by then, for
each $i<\ell$. The factor with $i=0$ is one unless the shortest leadtime is
zero. These probabilities sum to one, and $\plag$ is also the distribution
of the lag of a period.

Suppose that the order of period $s$ governs the period $s+\ell$. The
outstanding order of age $k$ at $s$ is missing from that period if it has
still not been received $\ell$ periods later. This is the event
$B_{k,\ell}=1$ of \eqref{eq:B}, and it has probability $m_{k,\ell}$. So the
stock that has been received by the end of the period, before the demand
$D(\ell+1)$ of the periods $s$ to $s+\ell$ is subtracted, is
\begin{equation}\label{eq:received-lag}
\Ipre_\ell=z_s-\sum_{k=1}^{\rng-1}a_{k,s}B_{k,\ell},
\end{equation}
and the net inventory of the period is $\Ipre_\ell-D(\ell+1)$. This is the
counterpart of \eqref{eq:received} at lag $\ell$. Let
$G_\ell(u)=\E\,\cbar\big(u-D(\ell+1)\big)$ be the expected cost of a period
at lag $\ell$ in which the received stock is $u$. The smoothed cost $G$ of
\eqref{eq:Gamma} is the member of this family at the shortest lag,
$G=G_{\Lmin}$. Averaging over the lag and
over the missing orders gives a function of the position $z$ after ordering
and of the outstanding orders $a$:
\begin{equation}\label{eq:Phi}
\Phi(z,a)=\sum_{\ell=\Lmin}^{\Lmax}\plag_\ell\,\E\,G_\ell\Big(z-\sum_{k=1}^{\rng-1}a_kB_{k,\ell}\Big).
\end{equation}

\begin{proposition}[Cost identity]\label{prop:identity}
Under any policy, the expected cost of the periods governed by the order of
period $s$, given the information at $s$ and the order $q_s$, is
$\Phi(z_s,a_s)$.
\end{proposition}

The function $\Phi$ is the cost of a newsvendor whose supply is the position
less the outstanding orders that arrive too late, and whose demand is that
of the governed periods. It is jointly convex in $(z,a)$, because $\cbar$ is
convex and its argument is linear in $(z,a)$ for every outcome of the
leadtimes and the demand. Every period is governed by exactly one order,
apart from the first few periods, in which no order has been received yet.
So the costs $\Phi(z_s,a_s)$ of the successive orders add up to the cost of
the policy.

\begin{corollary}[Long-run average cost]\label{cor:longrun}
The long-run average cost of a policy is the long-run average of
$\E\,\Phi(z_s,a_s)$. So the problem of minimizing it is a Markov decision
problem with the state $(\IP,a)$, the action $z\ge\IP$, the stage cost
$\Phi(z,a)$ and the transition \eqref{eq:transy} written in these
coordinates (Section~\ref{sec:ec-model} of the e-companion).
\end{corollary}

\subsection{The Weights and the Damping Factor}\label{sec:quadmodel}

The two policies start from the same probability and average it in two
ways. The probability is $m_{k,j}$ of \eqref{eq:m}, that the outstanding
order of age $k$ has still not been received $j$ periods from now. The
closed-form weight \eqref{eq:rk} averages it over the leadtime $j$ of the
new order. So it looks at one moment, the arrival of the new order. The
cost-aware weight averages it over the lag $\ell$ of the periods that the
new order governs, with more weight on the periods in which the stock is
close to running out.

\subsubsection*{The discount of an outstanding order.} By
Proposition~\ref{prop:identity}, the order of a period decides the cost
$\Phi(z,a)$ of the periods that it governs. Let $\hat z(a)$ be the position
that minimizes this cost for given outstanding orders $a$, and let the
demand have a density. At this position the holding cost of the outcomes
with stock left over balances the backorder cost of the outcomes with a
backlog. Now let the outstanding order of age $k$ be a little larger, at the
same position. The received stock of a governed period is then smaller by
the same amount if the order is missing from the period. This shifts the
balance only through the outcomes in which the net inventory is zero,
because only there does a small change of the stock turn leftover stock
into a backlog. Differentiating the first-order condition
$\partial\Phi/\partial z=0$ gives the rate at which the position has to rise
to restore the balance:
\begin{equation}\label{eq:foc}
\frac{\partial\hat z}{\partial a_k}
=\frac{\sum_\ell\plag_\ell\,\E\big[B_{k,\ell}\,f_\ell(\Ipre_\ell)\big]}
{\sum_\ell\plag_\ell\,\E\,f_\ell(\Ipre_\ell)} ,
\end{equation}
where $f_\ell$ is the density of the demand $D(\ell+1)$, which enters
because the second derivative of $G_\ell$ is $(b+h)f_\ell$, and
$\Ipre_\ell$ is the received stock \eqref{eq:received-lag} at the position
$\hat z(a)$. The rate is defined where the denominator is positive.
The left side is a rate, and the right side is a conditional probability.
The net inventory of a period at lag $\ell$ is $\Ipre_\ell-D(\ell+1)$, so
its density at zero is $f_\ell(\Ipre_\ell)$. The denominator is therefore
the density of a net inventory of zero in a governed period, and the
numerator is the part of this density that comes from the outcomes in which
the order of age $k$ is missing. So the ratio is the probability that the
order of age $k$ is missing from a governed period, given that the net
inventory of the period is zero. The order then falls by one less this
probability for each unit of the outstanding order, and by one for each
unit of stock. So this is the discount that the cost calls for: an
outstanding order counts for less than stock only to the extent that it is
missing when the stock just runs out.

\subsubsection*{Computing the probability.} The right side of
\eqref{eq:foc} depends on the state, and a rule with constant weights needs
one number for each age. So we evaluate it at typical values of the state.
Let $O$ be the set of the ages $k<\rng$ at which an order is outstanding
when the new order is placed. The expectations in \eqref{eq:foc} are over
the orders that are missing from the period at lag $\ell$. These are the
orders of the ages in a set $B\subseteq O$ with the probability
$\P(B\mid O,\ell)=\prod_{k\in B}m_{k,\ell}\prod_{k\in O\setminus B}(1-m_{k,\ell})$,
and the received stock is then $\Ipre_\ell=z-\sum_{k\in B}a_k$. We replace
every outstanding order by its mean size $\mu$, and the position by its mean
$\bar z_O$ when the orders of the ages in $O$ are outstanding, which
Section~\ref{sec:coef2} gives. So the received stock becomes
$\bar z_O-\mu|B|$. The set $O$ is random too, with the probability
$\P(O)=\prod_{k\in O}\P(L\ge k)\prod_{k\notin O}\P(L<k)$, and we average
over it as well. We call a triple $(O,\ell,B)$ a case. Its term in the sums
of \eqref{eq:foc} is
\begin{equation}\label{eq:caseweight}
w(O,\ell,B)=\P(O)\,\plag_\ell\,\P(B\mid O,\ell)\,f_\ell\big(\bar z_O-\mu|B|\big),
\end{equation}
the probability of the case times the density of a net inventory of zero in
it. With integer demand $f_\ell$ is the probability function of
$D(\ell+1)$, interpolated between the integers, and it takes the place of
the density (Section~\ref{sec:ec-quad} of the e-companion). The numerator of
\eqref{eq:foc} is the sum of these terms over the cases in which the order
of age $k$ is missing. The denominator is the sum over the cases in which it
is outstanding, because the rate matters only when there is an order of age
$k$. This gives $\psi_k$ below, and $\Omega_{jk}$ is the same ratio for the
orders of two ages $j\ne k$ together:
\begin{equation}\label{eq:psi}
\psi_k=\frac{\sum_{k\in B}w(O,\ell,B)}{\sum_{k\in O}w(O,\ell,B)},\qquad
\Omega_{jk}=\frac{\sum_{j,k\in B}w(O,\ell,B)}{\sum_{j,k\in O}w(O,\ell,B)},
\end{equation}
where each sum is over the cases that satisfy the condition under it, and
$\Omega_{kk}=\psi_k$. So $\psi_k$ is the right side of \eqref{eq:foc} at
typical values of the state: the probability that the order of age $k$ is
missing, given that it is outstanding and that the net inventory is zero.
$\Omega_{jk}$ is the probability that two orders are both missing, given
that both are outstanding and that the net inventory is zero, and the
damping factor below needs it.

Without the condition on the net inventory, $\psi_k$ would be
$\sum_\ell\plag_\ell\,m_{k,\ell}$, the probability that the order of age $k$
is missing from the period that the new order governs. This is the
counterpart of $1-\theta_k=\sum_jp_j\,m_{k,j}$ in the closed-form weight
\eqref{eq:rk}, with the distribution of the lag in place of the
distribution of the leadtime. The condition moves the probability toward
the cases in which a net inventory of zero is likely, and this is how the
cost ratio enters. With a high cost ratio the position is high. A net
inventory of zero is then unlikely at a short lag, where little demand has
accumulated, and more likely at a long lag. So when the stock just runs out
the lag is usually long, and at a long lag an outstanding order has usually
arrived. So $\psi_k$ is small. With a low cost ratio the stock also runs
out at short lags, and $\psi_k$ is larger.

\subsubsection*{The coefficients.} The damping factor and the weights of
the cost-aware policy are computed from $\psi$ and $\Omega$ by a recursion
on matrices of size $\rng$. The recursion solves a version of the problem
with a quadratic cost, in which the best policy is a damped response to a
weighted inventory position (Proposition~\ref{prop:quadw} in
Section~\ref{sec:ec-quad} of the e-companion).

\subsection{The Fixed Point}\label{sec:coef2}

The weight \eqref{eq:caseweight} of a case uses the mean position
$\bar z_O$, and this position depends on the coefficients. A rule with the
damping factor $\kappa$ and the weights $\theta_k$ replaces the part
$\kappa(1-\theta_k)$ of the outstanding order of age $k$. So at a given
inventory position its position after ordering is higher by this part of
every outstanding order. An outstanding order has the mean size $\mu$, and
we center the positions at the level $S^{\IP}$ of Section~\ref{sec:level}.
So when the orders of the ages in $O$ are outstanding, we take the mean
position to be
\[
\bar z_O=S^{\IP}+\kappa\mu\sum_{k}(1-\theta_k)\big(1\{k\in O\}-\P(L\ge k)\big),
\]
which averages to $S^{\IP}$ over the sets $O$.

The coefficients therefore enter their own definition: they give the
positions $\bar z_O$, the positions give $\psi$ and $\Omega$ by
\eqref{eq:psi}, and these give new coefficients by the recursion of
Section~\ref{sec:quadmodel}. We define the cost-aware coefficients
$\kappa^\beta$ and $\theta^\beta_k$ as a fixed point of this map. Besides
the critical ratio they depend on the demand distribution, through $f_\ell$,
$\mu$ and $S^{\IP}$.

Two remarks complete the definition. Where no age has $0<\P(L\ge k)<1$, as
at range two with a shortest leadtime of one period, the positions do not
depend on the coefficients, and the map gives the coefficients in one step.
And the recursion needs $\Omega-\psi\psi^\top$ to be positive semidefinite,
which \eqref{eq:psi} guarantees only at range two, so we replace $\Omega$ by
a nearby matrix with this property (Section~\ref{sec:ec-p2alg} of the
e-companion).

A fixed point exists when the demand takes every nonnegative integer value
with positive probability, as Poisson demand does, and the level $S^{\IP}$
is not too low compared with the mean demand. Both hold in every instance of
Section~\ref{sec:numerical} (Proposition~\ref{prop:exist} of the
e-companion). With a shortest leadtime of
zero we found instances in which the map has three fixed points
(Section~\ref{sec:ec-p2alg}). The cost-aware coefficients are the fixed
point that the iteration from the closed-form coefficients reaches, and the
iteration converged in every instance that we computed. The weights are
not restricted to the unit interval, and on the test bed they are slightly
negative in \nv{p2.negw.1} instances (Section~\ref{sec:ec-p2alg}).

\subsection{The Level}\label{sec:level2}

The level of the cost-aware policy is
\begin{equation}\label{eq:S2}
S^\beta=S^{\IP}+\bar z(\kappa^\beta,\theta^\beta)-\bar z(1,1)
-\mu\chi^\beta+\frac{1-\kappa^\beta}{\kappa^\beta}\,\mu,
\end{equation}
where $\chi^\beta=\sum_k(1-\theta^\beta_k)\,\P(L\ge k)$. The last two terms
are as in \eqref{eq:S}: the two corrections, for the weights and for the
damping, with the cost-aware coefficients. The first three terms replace the mean
position $S^{\IP}$ of \eqref{eq:S}, because the coefficients change the
distribution of the net inventory and with it the mean position that is
needed. Here $\bar z(\kappa,\theta)$ is the mean position at which the net
inventory under the rule with the coefficients $(\kappa,\theta)$ is
nonnegative with probability $\beta$, in a model of the net inventory as a
mixture of normal distributions over the cases
(Section~\ref{sec:ec-level2} of the e-companion), and $\bar z(1,1)$ is this
position for the base-stock policy. With the coefficients of the base-stock
policy the level is $S^{\IP}$ exactly.

With the coefficients of Section~\ref{sec:coef2} and this level, the rule
\eqref{eq:policy2} is the cost-aware policy. Nothing in it is searched or
simulated, and on the test bed its computation takes a fraction of a second
for each instance. The box collects the steps.

\smallskip
\begin{center}
\fbox{\begin{tabular}{@{}r@{\ \ }l@{}}
\multicolumn{2}{@{}l@{}}{\textbf{Computing the cost-aware policy} from the leadtime and demand distributions, $h$ and $b$}\\
1:& $(\kappa,\theta)\gets$ the coefficients of the closed-form policy, \eqref{eq:rk} and \eqref{eq:kappa}\\
2:& \textbf{repeat}\\
3:& \quad \textbf{for} every set $O$ of ages with an outstanding order: compute the position $\bar z_O$ (Section~\ref{sec:coef2})\\
4:& \quad \textbf{for} every case $(O,\ell,B)$: compute the weight \eqref{eq:caseweight}\\
5:& \quad compute $\psi$ and $\Omega$ from the weights by \eqref{eq:psi}\\
6:& \quad project $\Omega$ so that $\Omega-\psi\psi^\top$ is positive semidefinite with a margin (Section~\ref{sec:ec-p2alg})\\
7:& \quad $(\kappa',\theta')\gets$ the solution of the quadratic problem with $\psi$ and $\Omega$ (Section~\ref{sec:ec-quad})\\
8:& \quad \textbf{if} $(\kappa',\theta')$ is within $10^{-7}$ of $(\kappa,\theta)$: leave the loop\\
9:& \quad move $(\kappa,\theta)$ toward $(\kappa',\theta')$\\
10:& $(\kappa^\beta,\theta^\beta)\gets(\kappa,\theta)$; compute the level $S^\beta$ by \eqref{eq:S2}; the policy is \eqref{eq:policy2}
\end{tabular}}
\end{center}

\section{Lower Bounds}\label{sec:bounds}

The two policies have to be measured against the optimal cost. Where the
dynamic program can be solved, the optimal cost itself is the benchmark, and
Section~\ref{sec:numerical} uses it at ranges up to four. At longer ranges
the dynamic program cannot be solved, because the state has $\rng$
dimensions. A lower bound on the cost of every policy then takes the place
of the optimal cost. This section gives such bounds, and none of them needs
a dynamic program.

By Corollary~\ref{cor:longrun} the long-run average cost of a policy is the
long-run average of $\Phi(z_s,a_s)$, and $\Phi$ is jointly convex. So by
Jensen's inequality the cost is at least $\Phi$ at the long-run means of the
position and of the outstanding orders. The means of the outstanding orders
are the same under every policy with a finite cost. The average order equals
the mean demand $\mu$, and the size of an order is chosen before its
leadtime is drawn, so the outstanding order of age $k$ has the mean
$\bar a_k=\mu\,\P(L\ge k)$. Only the mean of the position is free, and
minimizing over it gives the bound
\begin{equation}\label{eq:lb0}
\LB_0=\min_z\ \Phi(z,\bar a).
\end{equation}
It is at least the lower bound of \citet[][Lemma~1]{gongsun2026}. It
evaluates $\Phi$ at a constant position and at orders of the mean size, and
no policy has both, because orders and positions are linked by
\begin{equation}\label{eq:link}
q_s=z_s-z_{s-1}+D_{s-1}.
\end{equation}
To use this link, we apply Jensen's inequality given the last $k$ demands,
$D_{s-1},\dots,D_{s-k}$. A policy then enters the bound through its mean
positions given these demands, and \eqref{eq:link} gives its mean orders. We
state the bound at range two with a shortest leadtime of at least one
period. There the outstanding order is the order of the previous period, and
the cost of a policy is the long-run average of $\Phi(z_s,q_{s-1})$. Let
$d=(d_1,\dots,d_k)$ denote values of the last $k$ demands, most recent
first, with probabilities $p(d)$, and let $z(d)$ stand for the mean position
after ordering given $d$. The mean position one period earlier, given the
demands that it had seen, is
$\bar z_1(d_2,\dots,d_k)=\sum_{d'}p(d')\,z(d_2,\dots,d_k,d')$, and $\bar z_2$
is defined from $\bar z_1$ in the same way. By \eqref{eq:link} the last
order has the conditional mean
$\zeta(d)=\bar z_1(d_2,\dots,d_k)-\bar z_2(d_3,\dots,d_k)+d_2$ when $k\ge2$,
and $\zeta=\mu$ when $k\le1$. Define
\begin{equation}\label{eq:lbk2}
\begin{aligned}
\LB_k=\min_{z(\cdot)}\ &\sum_dp(d)\,\Phi\big(z(d),\zeta(d)\big)\\
\text{s.t.}\ \ &\zeta(d)\ge0,\qquad z(d)-\bar z_1(d_2,\dots,d_k)+d_1\ge0 .
\end{aligned}
\end{equation}
The constraints say that orders are nonnegative in conditional mean, and the
second applies for $k\ge1$. With $k=0$ the program is \eqref{eq:lb0}. At a
general range the program has the same form, and it also conditions on which
orders have been received (Section~\ref{sec:ec-lb} of the e-companion).

\begin{proposition}[Lower bounds]\label{prop:lbk}
Let the demand take finitely many values, or let the conditioning be on a
finite partition of its values. For every $k\ge0$, $\LB_k\le g^*$, and
$\LB_k$ is nondecreasing in $k$.
\end{proposition}

With integer demand, $\Phi$ is piecewise linear and \eqref{eq:lbk2} is a
linear program. It grows with the number of demand values to the power $k$,
which limits it to short ranges and moderate demand. The proof of
Proposition~\ref{prop:lbk} is in Section~\ref{sec:ec-lb}.

As the shortest leadtime grows while the range stays fixed, the bound
$\LB_0$ becomes tight. For a rule whose coefficients do not change with the
shortest leadtime, such as the base-stock policy and the closed-form policy,
the cost at the best level exceeds $\LB_0$ by a fraction that falls at least as
fast as one over the shortest leadtime, with a constant that is given by the variance of
the stock that the rule has received (Proposition~\ref{prop:long} in
Section~\ref{sec:long} of the e-companion). So the optimality gap of such a
rule vanishes at least at this rate. That section also gives an upper bound on the
cost in closed form, and it reports
how close both policies come to $\LB_0$ at shortest leadtimes of up to 500
periods.

\section{Numerical Study}\label{sec:numerical}

We evaluate the closed-form policy \eqref{eq:policy}, the cost-aware policy
\eqref{eq:policy2} and the lower bounds of Section~\ref{sec:bounds} on a
test bed of \nv{bed.n} instances, and we compare the two policies with those
that have been published for this problem.

Demand is Poisson, and orders are integers, rounded to the nearest integer
with a half rounded up, except under B3 and B6 of
Section~\ref{sec:published}, which round up as their authors do. Every cost is computed from the transition law and
nothing is simulated, except in the cases that Section~\ref{sec:published}
names, and we call such a cost exact. The computations are numerical, on a
truncated model. So the optimal costs, the lower bounds and the certificates
below are values for that model, not error bounds for the original one
(Section~\ref{sec:ec-pub} of the e-companion).

\subsection{The Test Bed}\label{sec:bed}

The test bed is a full factorial design in four factors
(Table~\ref{tab:bed}): three ranges of the leadtime, nine leadtime
distributions at each range, five cost ratios and three mean demands. The
two-point distributions put all the probability on the shortest and the
longest leadtime, and they vary the skewness. The symmetric distributions put
the same probability on each end and the rest in the middle, and they vary
how much of the probability is concentrated in the middle. The shortest
leadtime is one period. Section~\ref{sec:lmin} varies it from zero to five
periods, and Section~\ref{sec:published} adds published instances with a
shortest leadtime of zero.

\begin{table}
\TABLE
{The Test Bed\label{tab:bed}}
{\begin{tabular}{@{}l@{\quad}l@{\quad}l@{}}
\hline\up
Family & Range & Leadtime distributions\\ \hline\up
Two-point & 2 & $(.9,0,.1)$, $(.7,0,.3)$, $(.5,0,.5)$, $(.3,0,.7)$, $(.1,0,.9)$\\
 & 3 & $(.9,0,0,.1)$, $(.7,0,0,.3)$, $(.5,0,0,.5)$, $(.3,0,0,.7)$, $(.1,0,0,.9)$\\
 & 4 & $(.9,0,0,0,.1)$, $(.7,0,0,0,.3)$, $(.5,0,0,0,.5)$, $(.3,0,0,0,.7)$, $(.1,0,0,0,.9)$\\
Symmetric & 2 & $(.4,.2,.4)$, $(.3,.4,.3)$, $(.2,.6,.2)$, $(.1,.8,.1)$\\
 & 3 & $(.4,.1,.1,.4)$, $(.3,.2,.2,.3)$, $(.2,.3,.3,.2)$, $(.1,.4,.4,.1)$\\
 & 4 & $(.4,0,.2,0,.4)$, $(.3,0,.4,0,.3)$, $(.2,0,.6,0,.2)$, $(.1,0,.8,0,.1)$\down\\ \hline
\end{tabular}}
{A distribution at range $\rng$ lists the probabilities of the leadtimes
$1,\dots,\rng+1$. Each distribution is combined with the cost ratios
$b/h=1,2,4,9,19$, with $h=1$, and with Poisson demand of mean 5, 10 and 20,
which gives 15 instances per distribution and \nv{bed.n} in all.}
\end{table}

\subsection{Optimality Gaps}\label{sec:longrun}

The optimal cost is computed by relative value iteration on the state
$(\IP,a)$ of Corollary~\ref{cor:longrun}. The base-stock policy has its best
level, and nothing in the two policies is searched. As a reference for the
form of the two policies we also search, for each instance, over the damping
factor, the weights and the level of the rule
$q_t=\kappa'(S'-y_t-\sum_k\theta'_ka_{k,t})^+$, and we call the result the
best rule of the form. The search is local, so the true best rule costs no
more.

\begin{table}
\TABLE
{Optimality Gaps and Lower Bounds (Mean and Largest, Percent)\label{tab:stat}}
{\begin{tabular}{@{}l@{\quad}c@{\quad}c@{\ \ }c@{\ \ }c@{\ \ }c@{\qquad}c@{\ \ }c@{}}
\hline\up
 & & \multicolumn{4}{c}{Gap to the optimal cost} & \multicolumn{2}{c}{Bound below}\\
 & & Base & Closed & Cost & Best & \multicolumn{2}{c}{the optimal cost}\\
Instances & Number & stock & form & aware & rule & $\LB_0$ & $\LB_k$\\ \hline\up
All & 405 & \nv{lr.all.bs.mean} (\nv{lr.all.bs.max}) & \nv{lr.all.paper.mean} (\nv{lr.all.paper.max}) & \nv{lr.all.p2.mean} (\nv{lr.all.p2.max}) & \nv{lr.all.lin.mean} (\nv{lr.all.lin.max}) & \nv{lr.all.lbu.mean} (\nv{lr.all.lbu.max}) & \nv{lr.all.lbk.mean} (\nv{lr.all.lbk.max})\down\\%
\hline\up Two-point distributions & 225 & \nv{lr.skew.bs.mean} (\nv{lr.skew.bs.max}) & \nv{lr.skew.paper.mean} (\nv{lr.skew.paper.max}) & \nv{lr.skew.p2.mean} (\nv{lr.skew.p2.max}) & \nv{lr.skew.lin.mean} (\nv{lr.skew.lin.max}) & \nv{lr.skew.lbu.mean} (\nv{lr.skew.lbu.max}) & \nv{lr.skew.lbk.mean} (\nv{lr.skew.lbk.max})\\%
\quad Range 2 & 75 & \nv{lr.r2skew.bs.mean} (\nv{lr.r2skew.bs.max}) & \nv{lr.r2skew.paper.mean} (\nv{lr.r2skew.paper.max}) & \nv{lr.r2skew.p2.mean} (\nv{lr.r2skew.p2.max}) & \nv{lr.r2skew.lin.mean} (\nv{lr.r2skew.lin.max}) & \nv{lr.r2skew.lbu.mean} (\nv{lr.r2skew.lbu.max}) & \nv{lr.r2skew.lbk.mean} (\nv{lr.r2skew.lbk.max})\\%
\quad Range 3 & 75 & \nv{lr.r3skew.bs.mean} (\nv{lr.r3skew.bs.max}) & \nv{lr.r3skew.paper.mean} (\nv{lr.r3skew.paper.max}) & \nv{lr.r3skew.p2.mean} (\nv{lr.r3skew.p2.max}) & \nv{lr.r3skew.lin.mean} (\nv{lr.r3skew.lin.max}) & \nv{lr.r3skew.lbu.mean} (\nv{lr.r3skew.lbu.max}) & \nv{lr.r3skew.lbk.mean} (\nv{lr.r3skew.lbk.max})\\%
\quad Range 4 & 75 & \nv{lr.r4skew.bs.mean} (\nv{lr.r4skew.bs.max}) & \nv{lr.r4skew.paper.mean} (\nv{lr.r4skew.paper.max}) & \nv{lr.r4skew.p2.mean} (\nv{lr.r4skew.p2.max}) & \nv{lr.r4skew.lin.mean} (\nv{lr.r4skew.lin.max}) & \nv{lr.r4skew.lbu.mean} (\nv{lr.r4skew.lbu.max}) & \nv{lr.r4skew.lbk.mean} (\nv{lr.r4skew.lbk.max})\\%
Symmetric distributions & 180 & \nv{lr.mid.bs.mean} (\nv{lr.mid.bs.max}) & \nv{lr.mid.paper.mean} (\nv{lr.mid.paper.max}) & \nv{lr.mid.p2.mean} (\nv{lr.mid.p2.max}) & \nv{lr.mid.lin.mean} (\nv{lr.mid.lin.max}) & \nv{lr.mid.lbu.mean} (\nv{lr.mid.lbu.max}) & \nv{lr.mid.lbk.mean} (\nv{lr.mid.lbk.max})\\%
\quad Range 2 & 60 & \nv{lr.r2mid.bs.mean} (\nv{lr.r2mid.bs.max}) & \nv{lr.r2mid.paper.mean} (\nv{lr.r2mid.paper.max}) & \nv{lr.r2mid.p2.mean} (\nv{lr.r2mid.p2.max}) & \nv{lr.r2mid.lin.mean} (\nv{lr.r2mid.lin.max}) & \nv{lr.r2mid.lbu.mean} (\nv{lr.r2mid.lbu.max}) & \nv{lr.r2mid.lbk.mean} (\nv{lr.r2mid.lbk.max})\\%
\quad Range 3 & 60 & \nv{lr.r3mid.bs.mean} (\nv{lr.r3mid.bs.max}) & \nv{lr.r3mid.paper.mean} (\nv{lr.r3mid.paper.max}) & \nv{lr.r3mid.p2.mean} (\nv{lr.r3mid.p2.max}) & \nv{lr.r3mid.lin.mean} (\nv{lr.r3mid.lin.max}) & \nv{lr.r3mid.lbu.mean} (\nv{lr.r3mid.lbu.max}) & \nv{lr.r3mid.lbk.mean} (\nv{lr.r3mid.lbk.max})\\%
\quad Range 4 & 60 & \nv{lr.r4mid.bs.mean} (\nv{lr.r4mid.bs.max}) & \nv{lr.r4mid.paper.mean} (\nv{lr.r4mid.paper.max}) & \nv{lr.r4mid.p2.mean} (\nv{lr.r4mid.p2.max}) & \nv{lr.r4mid.lin.mean} (\nv{lr.r4mid.lin.max}) & \nv{lr.r4mid.lbu.mean} (\nv{lr.r4mid.lbu.max}) & \nv{lr.r4mid.lbk.mean} (\nv{lr.r4mid.lbk.max})\down\\%
\hline\up $b/h=1$ & 81 & \nv{lr.b1.bs.mean} (\nv{lr.b1.bs.max}) & \nv{lr.b1.paper.mean} (\nv{lr.b1.paper.max}) & \nv{lr.b1.p2.mean} (\nv{lr.b1.p2.max}) & \nv{lr.b1.lin.mean} (\nv{lr.b1.lin.max}) & \nv{lr.b1.lbu.mean} (\nv{lr.b1.lbu.max}) & \nv{lr.b1.lbk.mean} (\nv{lr.b1.lbk.max})\\%
$b/h=2$ & 81 & \nv{lr.b2.bs.mean} (\nv{lr.b2.bs.max}) & \nv{lr.b2.paper.mean} (\nv{lr.b2.paper.max}) & \nv{lr.b2.p2.mean} (\nv{lr.b2.p2.max}) & \nv{lr.b2.lin.mean} (\nv{lr.b2.lin.max}) & \nv{lr.b2.lbu.mean} (\nv{lr.b2.lbu.max}) & \nv{lr.b2.lbk.mean} (\nv{lr.b2.lbk.max})\\%
$b/h=4$ & 81 & \nv{lr.b4.bs.mean} (\nv{lr.b4.bs.max}) & \nv{lr.b4.paper.mean} (\nv{lr.b4.paper.max}) & \nv{lr.b4.p2.mean} (\nv{lr.b4.p2.max}) & \nv{lr.b4.lin.mean} (\nv{lr.b4.lin.max}) & \nv{lr.b4.lbu.mean} (\nv{lr.b4.lbu.max}) & \nv{lr.b4.lbk.mean} (\nv{lr.b4.lbk.max})\\%
$b/h=9$ & 81 & \nv{lr.b9.bs.mean} (\nv{lr.b9.bs.max}) & \nv{lr.b9.paper.mean} (\nv{lr.b9.paper.max}) & \nv{lr.b9.p2.mean} (\nv{lr.b9.p2.max}) & \nv{lr.b9.lin.mean} (\nv{lr.b9.lin.max}) & \nv{lr.b9.lbu.mean} (\nv{lr.b9.lbu.max}) & \nv{lr.b9.lbk.mean} (\nv{lr.b9.lbk.max})\\%
$b/h=19$ & 81 & \nv{lr.b19.bs.mean} (\nv{lr.b19.bs.max}) & \nv{lr.b19.paper.mean} (\nv{lr.b19.paper.max}) & \nv{lr.b19.p2.mean} (\nv{lr.b19.p2.max}) & \nv{lr.b19.lin.mean} (\nv{lr.b19.lin.max}) & \nv{lr.b19.lbu.mean} (\nv{lr.b19.lbu.max}) & \nv{lr.b19.lbk.mean} (\nv{lr.b19.lbk.max})\down\\%
\hline\up Mean demand 5 & 135 & \nv{lr.mu5.bs.mean} (\nv{lr.mu5.bs.max}) & \nv{lr.mu5.paper.mean} (\nv{lr.mu5.paper.max}) & \nv{lr.mu5.p2.mean} (\nv{lr.mu5.p2.max}) & \nv{lr.mu5.lin.mean} (\nv{lr.mu5.lin.max}) & \nv{lr.mu5.lbu.mean} (\nv{lr.mu5.lbu.max}) & \nv{lr.mu5.lbk.mean} (\nv{lr.mu5.lbk.max})\\%
Mean demand 10 & 135 & \nv{lr.mu10.bs.mean} (\nv{lr.mu10.bs.max}) & \nv{lr.mu10.paper.mean} (\nv{lr.mu10.paper.max}) & \nv{lr.mu10.p2.mean} (\nv{lr.mu10.p2.max}) & \nv{lr.mu10.lin.mean} (\nv{lr.mu10.lin.max}) & \nv{lr.mu10.lbu.mean} (\nv{lr.mu10.lbu.max}) & \nv{lr.mu10.lbk.mean} (\nv{lr.mu10.lbk.max})\\%
Mean demand 20 & 135 & \nv{lr.mu20.bs.mean} (\nv{lr.mu20.bs.max}) & \nv{lr.mu20.paper.mean} (\nv{lr.mu20.paper.max}) & \nv{lr.mu20.p2.mean} (\nv{lr.mu20.p2.max}) & \nv{lr.mu20.lin.mean} (\nv{lr.mu20.lin.max}) & \nv{lr.mu20.lbu.mean} (\nv{lr.mu20.lbu.max}) & \nv{lr.mu20.lbk.mean} (\nv{lr.mu20.lbk.max})\down\\%
 
\hline
\end{tabular}}
{Each entry is a mean, with the largest value in parentheses. The gap of a
policy is its cost in excess of the optimal cost, in percent of the optimal
cost. The base-stock policy has its best level. ``Closed form'' is policy
\eqref{eq:policy} with its level from \eqref{eq:S}, and ``Cost aware'' is
policy \eqref{eq:policy2} with its level from \eqref{eq:S2}. ``Best rule''
is the rule of the same form with its damping factor, weights and level
searched for each instance. The last two columns give how far the lower
bounds of Section~\ref{sec:bounds} are below the optimal cost, in percent
of the optimal cost: the bound $\LB_0$, and the bound $\LB_k$ with $k$, the
number of past demands on which it conditions, equal to the range. The rows
by cost ratio and by mean demand are over all leadtime distributions.}
\end{table}

Table~\ref{tab:stat} gives the gaps. The best rule of the form is
\nv{lr.all.lin.mean}\% above the optimal cost on average and
\nv{lr.all.lin.max}\% at most. So a damped response to a weighted inventory
position, with constant coefficients, loses very little to the optimal
policy on this test bed.

The mean gaps of
the base-stock policy, the closed-form policy and the cost-aware policy are
\nv{lr.all.bs.mean}\%, \nv{lr.all.paper.mean}\% and \nv{lr.all.p2.mean}\%,
and the largest gaps \nv{lr.all.bs.max}\%, \nv{lr.all.paper.max}\% and
\nv{lr.all.p2.max}\%. So the closed-form policy removes most of the gap of
the base-stock policy, and the cost-aware policy removes most of what is
left. The medians are smaller, \nv{lr.all.bs.med}\%, \nv{lr.all.paper.med}\%
and \nv{lr.all.p2.med}\%, because many instances are easy for every policy.
The cost-aware policy costs less than the closed-form policy in
\nv{lr.all.p2.belowp1} instances and more in \nv{lr.all.p2.abovep1}, by at
most \nv{lr.all.p2.excessp1} percentage points.

The gaps are largest under the two-point distributions in which an order is
rarely late: under $(.9,0,0,.1)$ the base-stock policy is
\nv{lr.r3q9.bs.mean}\% above the optimal cost on average, the closed-form
policy \nv{lr.r3q9.paper.mean}\% and the cost-aware policy
\nv{lr.r3q9.p2.mean}\%. The
reason is the timing of a replacement. Under a two-point distribution an
order that is not received after the shortest leadtime is known to be late,
and it arrives $\rng$ periods later than the short leadtime would have
brought it. An order placed as soon as this is known arrives $\Lmin+1$
periods after that time if its own leadtime is short. So a replacement can
arrive before the late order when $\rng>\Lmin+1$, which with a shortest
leadtime of one period means a range of three or more. The optimal policy
orders such a replacement, and the base-stock policy cannot, because it
counts the late order in full.

For the lower bounds, conditioning on the past demands helps: over all
instances $\LB_k$ is \nv{lr.all.lbk.mean}\% below the optimal cost on
average, and $\LB_0$ is \nv{lr.all.lbu.mean}\% below it. The bound $\LB_k$
is tight under the symmetric distributions and at range two. It is loose
under the two-point distributions at ranges three and four, where a
replacement can arrive before a late order and the best order depends on
more of the history than the program conditions on. The cost of a policy in
excess of the bound is a certificate: it bounds the gap of the policy
without the optimal cost. It is about the sum of two entries of the table,
the gap of the policy and the distance of the bound below the optimal cost.

The lower blocks of Table~\ref{tab:stat} give the same quantities by cost
ratio and by mean demand. As the cost ratio grows from 1 to 19, the mean gap
of the base-stock policy falls from \nv{lr.b1.bs.mean}\% to
\nv{lr.b19.bs.mean}\%, and that of the closed-form policy rises from
\nv{lr.b1.paper.mean}\% to \nv{lr.b19.paper.mean}\%, because its
coefficients do not depend on the cost ratio. That of the cost-aware policy
is \nv{lr.b1.p2.mean}\% and \nv{lr.b19.p2.mean}\%. The reason is the one
given in Section~\ref{sec:quadmodel}. The optimal policy replaces
a part of a late order at once, and less of it the higher the cost ratio. So
the base-stock policy can cost less than the closed-form policy. It does in
\nv{lr.all.above.t} of the \nv{lr.all.n} instances, among them
\nv{lr.b19.above.t} of the \nv{lr.b19.n} with $b/h=19$. The loss is in the
closed-form coefficients and not in the form, since the best rule of the
form is close to the optimum in these instances. The cost-aware policy is
never above the base-stock policy by more than $0.005$ percentage points.
The mean demand matters less.

\subsection{Other Shortest Leadtimes}\label{sec:lmin}

\begin{figure}
\FIGURE
{\includegraphics{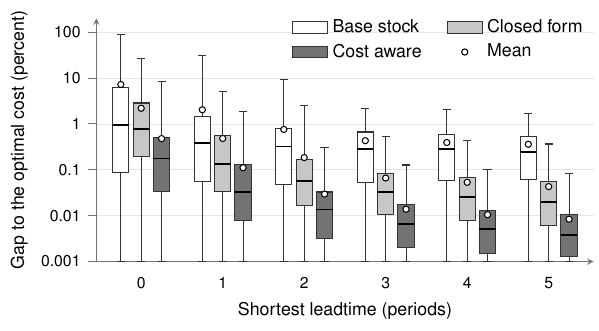}}
{Optimality Gaps by Shortest Leadtime\label{fig:lmin}}
{Each box has the \nv{bed.n} instances of the test bed, with every leadtime
distribution shifted to that shortest leadtime. A box spans the quartiles
of the gap to the optimal cost, with a line at the median. The whiskers
reach the smallest and the largest gap, and the circle is the mean. The
vertical axis is logarithmic, and gaps below $0.001\%$ are drawn at
$0.001\%$.}
\end{figure}

Figure~\ref{fig:lmin} shifts every distribution of the test bed so that the
shortest leadtime runs from zero to five periods, with the coefficients of
the two policies recomputed at each. All three policies come closer
to the optimum as the shortest leadtime grows, and the two policies stay
ahead of the base-stock policy on average. There are two reasons. A longer
shortest leadtime adds demand that no policy can influence. And a
replacement can arrive before a late order only when
$\rng>\Lmin+1$ (Section~\ref{sec:longrun}): at range three the mean gap of
the base-stock policy falls from \nv{lmin.r3.base_stock.1.mean}\% to
\nv{lmin.r3.base_stock.2.mean}\% between a shortest leadtime of one and two
periods, and at range four from \nv{lmin.r4.base_stock.2.mean}\% to
\nv{lmin.r4.base_stock.3.mean}\% between two and three. The gaps keep
falling at longer shortest leadtimes. Section~\ref{sec:long} of the
e-companion proves that the gap of a rule with fixed coefficients, at its
best level, vanishes at least as fast as one over the shortest leadtime,
and it gives the gaps of the
three policies at shortest leadtimes of 20 and 100 periods.

A shortest leadtime of zero is the hardest case. The mean gaps of the
base-stock policy, the closed-form policy and the cost-aware policy are
\nv{lmin.all.base_stock.0.mean}\%, \nv{lmin.all.complete.0.mean}\% and
\nv{lmin.all.p2.0.mean}\% there, and the cost-aware policy is no longer
ahead in every instance. It costs more than the base-stock
policy in \nv{lmin.all.p2.abovebs.0} instances, by at most
\nv{lmin.all.p2.abovebs.0.max} percentage points, and more than the
closed-form policy in \nv{lmin.all.p2.abovep1.0}, by at most
\nv{lmin.all.p2.abovep1.0.max}. Of these \nv{lmin.all.p2.abovep1.0},
\nv{lmin.all.p2.abovep1.0.mid} are under the symmetric distributions and
\nv{lmin.all.p2.abovep1.0.b12} have $b/h\le2$.

\subsection{Published Policies}\label{sec:published}

We now price the policies that have been proposed for random leadtimes with
crossover on the same instances. We compare
six, and we label them B1 to B6. The base-stock policy (B1) orders up to a
level on the inventory position \citep{zalkind1978,robinson2001}. The
proportional order-up-to policy of \citet{disney2016} (B2) orders a fraction of the gap between a level and the
inventory position, which is policy \eqref{eq:policy} with every weight
equal to one. The generalized base-stock policy of \citet{stolyarwang2022},
in the discrete form of \citet{stolyar2026} (B3), orders up to a level on the
outstanding orders plus a multiple $\gamma$ of the net inventory. The
two-level base-stock policy of \citet{weiyang2026} (B4) orders up to a low level
on the inventory position when the net inventory is positive and up to a
high level when it is not. The weighted index policy of
\citet{srinivasan2011} (B5) orders up to the best base-stock level on a weighted
sum of the stock and the outstanding orders, with one weight for each age.
The stock has a weight too, which acts as a damping factor tied to the
level. In our searches it stays at \nv{wi.w0.min} or more, so the policy is
close to an order-up-to rule on a weighted position.
PIPELINE \citep{stolyar2026}, which is B6, places the first order of a linear program
that is solved in every period, in which demand is replaced by its mean and
each outstanding order by its expected arrivals.

Three of these policies were proposed for a setting that differs from ours.
B6 is designed for units that each have their own leadtime. When the units
of an order share one leadtime, the order arrives whole or not at all, and
the first order of its program then varies more than demand does. B3 and B4 were
proposed for continuous time with a leadtime for each unit, and they are asymptotically
optimal in regimes with long leadtimes or high demand rates. On our test bed
an order arrives whole, the leadtimes are short and the demand is moderate.
So the comparison shows how these policies carry over to our model, and it
says nothing about their performance in the settings for which they were
designed.

Each policy is given the best parameters that we could find for each
instance by a local search from several starting points. So each row of
Table~\ref{tab:pub} is an upper bound on the best policy of its class, and a
more thorough search could lower it. B3, B4 and B6 cannot be priced
exactly everywhere and are priced by long simulations there, with a standard
error of at most \nv{net.r4.se} percentage points for B3 and B4. Section~\ref{sec:ec-pub} of the
e-companion gives the searches, the simulations and the adaptations of
these policies in detail.

\begin{table}
\TABLE
{Published Policies: Gap to the Optimal Cost (Percent)\label{tab:pub}}
{\begin{tabular}{@{}l@{\ \ }l@{\ \ }c@{\ \ }c@{\ \ }c@{\ \ }c@{\ \ }c@{\ \ }c@{\ \ }c@{\ \ }c@{\ \ }c@{\ \ }c@{\ \ }c@{}}
\hline\up
 & & & & Parameters & & & & \multicolumn{3}{c}{Mean at range} & \multicolumn{2}{c}{More / less than}\\
 & Policy & Damping & Weights & searched & Median & Mean & Largest & 2 & 3 & 4 & \eqref{eq:policy} & \eqref{eq:policy2}\\ \hline\up
B1 & Base stock & -- & -- & 1 & 0.38 & 2.04 & 30.31 & 0.24 & 1.92 & 3.94 & 278 / 104 & 324 / 0\\%
B2 & Proportional order-up-to & Yes & -- & 2 & 0.23 & 1.75 & 30.30 & 0.15 & 1.64 & 3.47 & 258 / 123 & 307 / 25\\%
B3 & Generalized base stock & -- & One & 2 & 0.38 & 1.51 & 18.77 & 0.24 & 1.52 & 2.75 & 275 / 107 & 324 / 0\\%
B4 & Two-level base stock & -- & -- & 2 & 0.38 & 1.73 & 26.26 & 0.24 & 1.68 & 3.27 & 277 / 105 & 324 / 0\\%
B5 & Weighted index & -- & By age & $R$ & 0.26 & 0.61 & 4.11 & 0.16 & 0.58 & 1.09 & 235 / 150 & 308 / 29\\%
B6 & PIPELINE & -- & -- & 2 & 4.70 & 4.66 & 16.96 & 3.03 & 6.22 & 4.73 & 363 / 41 & 381 / 0\\%
 & Closed-form policy \eqref{eq:policy} & Yes & By age & 0 & 0.13 & 0.48 & 5.07 & 0.09 & 0.54 & 0.81 &  & 326 / 14\\%
 & Cost-aware policy \eqref{eq:policy2} & Yes & By age & 0 & 0.03 & 0.11 & 1.84 & 0.02 & 0.10 & 0.22 & 14 / 326 & \\%
 & Best rule of the form & Yes & By age & $R+1$ & 0.02 & 0.07 & 1.12 & 0.01 & 0.06 & 0.13 & 0 / 351 & 0 / 243\down\\%
 
\hline
\end{tabular}}
{``Damping'' says whether a policy closes only a part of the gap between
its level and its position, and ``Weights'' whether it discounts the
outstanding orders, with one weight for all of them or with one for each
age. The stock weight of B5 acts as a damping factor tied to its level.
The best rule of the form has the form of policies \eqref{eq:policy}
and \eqref{eq:policy2}, with its damping factor, its weights and its level
searched. The last two columns count the instances in which the policy of the row costs more and
less than policy \eqref{eq:policy} and than policy \eqref{eq:policy2}, with
differences of less than $0.005$ percentage points, or of less than two
standard errors for a simulated gap, counted as ties. The counts in the text
follow the same rule.}
\end{table}

Table~\ref{tab:pub} gives the result, and three observations stand out.
First, damping alone helps on average and leaves the large gaps: B2 has a
smaller mean gap than the base-stock policy and nearly the same largest gap.
The large gaps come from late orders that should be replaced, which calls
for telling the outstanding orders apart. Second, weights on the outstanding
orders remove these gaps, as \citet{srinivasan2011} proposed. Their policy
B5, with $\rng$ weights searched for each instance, has the smallest mean
gap and the smallest largest gap among the published policies. The
closed-form policy sets weights of this kind without a search, and it costs
less than B5 in \nv{lr.all.wi.above.t} instances and more in
\nv{lr.all.wi.below.t}. The cost-aware policy costs less than B5 in
\nv{lr.all.wi.abovep2.t} instances and more in \nv{lr.all.wi.belowp2.t}.
Third, B3 and B4 improve on the base-stock policy where a replacement can
arrive before a late order, at ranges three and four. At range two our
searches return the base-stock policy itself, $\gamma=1$ in
\nv{net.r2.gamma1} of the \nv{net.r2.n} instances and two equal levels in
\nv{net.r2.equal}. At ranges three and four they find a multiple $\gamma>1$
under the two-point distributions. An outstanding order then has the weight
$1/\gamma$ relative to a unit of stock, the same at every age, which
replaces a part of a late order. A weight for each age can in addition
spare the orders that are still likely to be on time.

\subsubsection*{Published instances with a shortest leadtime of zero.} The
cases in the literature in which the base-stock policy is far from the
optimum have a shortest leadtime of zero \citep{robinson2001,gongsun2026},
where a replacement can arrive before a late order from range two on.
Table~\ref{tab:zero} gives the gaps on the \nv{zero.n} instances of
\citet[Table~1]{gongsun2026}, with the leadtime distributions $(.7,.1,.2)$
and $(.9,0,.1)$ on $\{0,1,2\}$. We compute the optimal cost and the gaps of
our policies exactly, and our optimal costs differ from those that Gong and
Sun print by at most \nv{zero.opt.diff}.

\begin{table}
\TABLE
{Published Instances with a Shortest Leadtime of Zero: Gap to the Optimal Cost (Mean and Largest, Percent)\label{tab:zero}}
{\begin{tabular}{@{}l@{\quad}c@{\quad}c@{\ \ }c@{\ \ }c@{\ \ }c@{\ \ }c@{\ \ }c@{}}
\hline\up
 & & Base & & & Closed & Cost & Best\\
Instances & Number & stock & B3 & B4 & form & aware & rule\\ \hline\up
All & 36 & 5.60 (43.84) & 2.48 (13.61) & 2.65 (10.91) & 3.94 (27.95) & 0.47 (5.20) & 0.17 (0.50)\\%
\multicolumn{8}{@{}l}{\up By leadtime distribution}\\%
\quad $(.7,.1,.2)$ & 18 & 0.78 (3.90) & 0.64 (2.73) & 0.63 (2.99) & 2.42 (4.74) & 0.32 (1.80) & 0.14 (0.50)\\%
\quad $(.9,0,.1)$ & 18 & 10.42 (43.84) & 4.31 (13.61) & 4.67 (10.91) & 5.45 (27.95) & 0.62 (5.20) & 0.19 (0.45)\\%
\multicolumn{8}{@{}l}{By backorder cost $b$}\\%
\quad $b=5$ & 6 & 18.02 (43.84) & 7.45 (13.61) & 6.51 (10.91) & 7.53 (27.95) & 1.56 (5.20) & 0.23 (0.50)\\%
\quad $b=10$ & 6 & 7.63 (16.70) & 3.34 (8.39) & 3.57 (8.03) & 4.18 (7.44) & 0.40 (0.60) & 0.23 (0.35)\\%
\quad $b=20$ & 6 & 3.84 (9.80) & 1.79 (5.00) & 2.51 (6.30) & 3.13 (4.46) & 0.34 (0.79) & 0.20 (0.45)\\%
\quad $b=40$ & 6 & 2.00 (5.51) & 1.06 (3.37) & 1.55 (4.27) & 2.98 (3.85) & 0.26 (0.48) & 0.14 (0.32)\\%
\quad $b=60$ & 6 & 1.30 (4.02) & 0.72 (2.41) & 1.05 (3.13) & 2.92 (4.42) & 0.15 (0.29) & 0.12 (0.23)\\%
\quad $b=80$ & 6 & 0.84 (3.05) & 0.50 (1.86) & 0.69 (2.43) & 2.87 (4.74) & 0.11 (0.29) & 0.09 (0.22)\down\\%
 
\hline
\end{tabular}}
{Each entry is the mean gap over the instances of the row, with the largest
gap in parentheses. The instances have $h=1$ and mean demand 10, 20 and 50.
B3 and B4 have the gaps that Gong and Sun print, where they are called GBS
and TLBS. ``Best rule'' is the best rule of the form that we found, the
smaller of a search over the damping factor, the weight and the level and of
the two policies.}
\end{table}

The closed-form policy gains on the base-stock policy at a low backorder
cost and loses to it as the backorder cost grows, because its weight does
not follow the cost. The cost-aware policy does, and it costs less than the
base-stock policy, B3 and B4 together in \nv{zero.p2.belowall} of the
\nv{zero.n} instances. The best rule of the form is within
\nv{zero.all.lin.max}\% of the optimum in every instance, so the form is
adequate at a shortest leadtime of zero as well.

\section{Conclusion}\label{sec:conclusion}

With random leadtimes orders can cross, and the classical base-stock policy
on the inventory position need not be optimal. Crossover calls for two
discounts. An outstanding order may not
arrive in time, so it should count for less than its size. And an order
that is placed later may arrive earlier, so the order placed now should not
close the whole gap. The two discounts give the policy its structure:
weights on the outstanding orders and a damped response.

We set the weights and the damping factor in two ways, without a search. In
the closed-form policy they come from the leadtime distribution alone. An
outstanding order counts with the probability that it arrives no later than
a new one, and the response is damped according to the expected number of
orders that overtake a given order. In the cost-aware policy they also use
the cost ratio. An outstanding order counts by how often it is missing when
the stock is about to run out. We also gave lower bounds on the cost of
every policy that need no dynamic program.

On a test bed of \nv{bed.n} instances the closed-form policy is
\nv{lr.all.paper.mean}\% above the optimal cost on average and the
cost-aware policy \nv{lr.all.p2.mean}\%, against \nv{lr.all.bs.mean}\% for
the base-stock policy, and the best rule of the same form is within
\nv{lr.all.lin.max}\% of the optimal cost in every instance. The two
policies also compare well with the published policies on these instances.

Some of our results need neither the long-run criterion nor the linear cost.
The structure of the optimal order at range two holds in every period of a
finite-horizon problem, and the cost identity and the lower bounds hold
under general convex costs. The proofs follow the same arguments and are
omitted.
Our model has no fixed ordering cost. A fixed cost makes orders infrequent,
and two orders that are placed at least as many periods apart as the range
of the leadtime cannot cross. So the problem without a fixed cost is the one
in which orders cross most often. A policy with a reorder point and an
order-up-to level on the weighted inventory position is left to future work.

\bibliographystyle{plainnat}
\bibliography{refs}

\begin{thebibliography}{32}
\providecommand{\natexlab}[1]{#1}
\providecommand{\url}[1]{\texttt{#1}}
\expandafter\ifx\csname urlstyle\endcsname\relax
  \providecommand{\doi}[1]{doi: #1}\else
  \providecommand{\doi}{doi: \begingroup \urlstyle{rm}\Url}\fi

\bibitem[Bai et~al.(2026)Bai, Chen, Li, and Stolyar]{bai2026}
Xingyu Bai, Xin Chen, Menglong Li, and Alexander~L. Stolyar.
\newblock Asymptotic optimality of semi-open-loop policies in inventory models
  with stochastic lead times.
\newblock \emph{Management Science}, 2026.
\newblock \doi{10.1287/mnsc.2023.00539}.
\newblock ePub ahead of print July 1.

\bibitem[Bischak et~al.(2014)Bischak, Robb, Silver, and Blackburn]{bischak2014}
Diane~P. Bischak, David~J. Robb, Edward~A. Silver, and Joseph~D. Blackburn.
\newblock Analysis and management of periodic review, order-up-to level
  inventory systems with order crossover.
\newblock \emph{Production and Operations Management}, 23\penalty0
  (5):\penalty0 762--772, 2014.

\bibitem[Bradley and Robinson(2005)]{bradley2005}
James~R. Bradley and Lawrence~W. Robinson.
\newblock Improved base-stock approximations for independent stochastic lead
  times with order crossover.
\newblock \emph{Manufacturing \& Service Operations Management}, 7\penalty0
  (4):\penalty0 319--329, 2005.

\bibitem[Chen and Yu(2005)]{chenyu2005}
Fangruo Chen and Bin Yu.
\newblock Quantifying the value of leadtime information in a single-location
  inventory system.
\newblock \emph{Manufacturing \& Service Operations Management}, 7\penalty0
  (2):\penalty0 144--151, 2005.

\bibitem[Disney et~al.(2016)Disney, Maltz, Wang, and Warburton]{disney2016}
Stephen~M. Disney, Arnold Maltz, Xun Wang, and Roger D.~H. Warburton.
\newblock Inventory management for stochastic lead times with order crossovers.
\newblock \emph{European Journal of Operational Research}, 248\penalty0
  (2):\penalty0 473--486, 2016.

\bibitem[Ehrhardt(1984)]{ehrhardt1984}
Richard Ehrhardt.
\newblock $(s, {S})$ policies for a dynamic inventory model with stochastic
  lead times.
\newblock \emph{Operations Research}, 32\penalty0 (1):\penalty0 121--132, 1984.

\bibitem[Feinberg and Liang(2022)]{feinbergliang2022}
Eugene~A. Feinberg and Yan Liang.
\newblock On the optimality equation for average cost {M}arkov decision
  processes and its validity for inventory control.
\newblock \emph{Annals of Operations Research}, 317:\penalty0 569--586, 2022.

\bibitem[Feinberg et~al.(2012)Feinberg, Kasyanov, and
  Zadoianchuk]{feinberg2012}
Eugene~A. Feinberg, Pavlo~O. Kasyanov, and Nina~V. Zadoianchuk.
\newblock Average cost {M}arkov decision processes with weakly continuous
  transition probabilities.
\newblock \emph{Mathematics of Operations Research}, 37\penalty0 (4):\penalty0
  591--607, 2012.

\bibitem[Feller(1971)]{feller1971}
William Feller.
\newblock \emph{An Introduction to Probability Theory and Its Applications},
  volume~II.
\newblock Wiley, New York, 2nd edition, 1971.

\bibitem[Gong and Sun(2026)]{gongsun2026}
Xiting Gong and Chenxi Sun.
\newblock Base-stock policies in backlogging inventory systems with stochastic
  lead times: New results and insights.
\newblock Working paper, 2026.

\bibitem[Hayya et~al.(2008)Hayya, Bagchi, Kim, and Sun]{hayya2008}
Jack~C. Hayya, Uttarayan Bagchi, Jeon~G. Kim, and Daewon Sun.
\newblock On static stochastic order crossover.
\newblock \emph{International Journal of Production Economics}, 114\penalty0
  (1):\penalty0 404--413, 2008.

\bibitem[Higham(2002)]{higham2002}
Nicholas~J. Higham.
\newblock Computing the nearest correlation matrix---a problem from finance.
\newblock \emph{IMA Journal of Numerical Analysis}, 22\penalty0 (3):\penalty0
  329--343, 2002.

\bibitem[Huh et~al.(2011)Huh, Janakiraman, and Nagarajan]{huh2011}
Woonghee~Tim Huh, Ganesh Janakiraman, and Mahesh Nagarajan.
\newblock Average cost single-stage inventory models: An analysis using a
  vanishing discount approach.
\newblock \emph{Operations Research}, 59\penalty0 (1):\penalty0 143--155, 2011.

\bibitem[Kaplan(1970)]{kaplan1970}
Robert~S. Kaplan.
\newblock A dynamic inventory model with stochastic lead times.
\newblock \emph{Management Science}, 16\penalty0 (7):\penalty0 491--507, 1970.

\bibitem[Murota(2003)]{murota2003}
Kazuo Murota.
\newblock \emph{Discrete Convex Analysis}.
\newblock SIAM Monographs on Discrete Mathematics and Applications. Society for
  Industrial and Applied Mathematics, Philadelphia, 2003.

\bibitem[Riezebos and Gaalman(2009)]{riezebos2009}
J.~Riezebos and G.~J.~C. Gaalman.
\newblock A single-item inventory model for expected inventory order
  crossovers.
\newblock \emph{International Journal of Production Economics}, 121\penalty0
  (2):\penalty0 601--609, 2009.

\bibitem[Riezebos(2006)]{riezebos2006}
Jan Riezebos.
\newblock Inventory order crossovers.
\newblock \emph{International Journal of Production Economics}, 104\penalty0
  (2):\penalty0 666--675, 2006.

\bibitem[Robinson et~al.(2001)Robinson, Bradley, and Thomas]{robinson2001}
Lawrence~W. Robinson, James~R. Bradley, and L.~Joseph Thomas.
\newblock Consequences of order crossover under order-up-to inventory policies.
\newblock \emph{Manufacturing \& Service Operations Management}, 3\penalty0
  (3):\penalty0 175--188, 2001.

\bibitem[Sch{\"a}l(1993)]{schal1993}
Manfred Sch{\"a}l.
\newblock Average optimality in dynamic programming with general state space.
\newblock \emph{Mathematics of Operations Research}, 18\penalty0 (1):\penalty0
  163--172, 1993.

\bibitem[Song(1994)]{song1994}
Jing-Sheng Song.
\newblock The effect of leadtime uncertainty in a simple stochastic inventory
  model.
\newblock \emph{Management Science}, 40\penalty0 (5):\penalty0 603--613, 1994.

\bibitem[Song and Zipkin(1996)]{songzipkin1996}
Jing-Sheng Song and Paul Zipkin.
\newblock Inventory control with information about supply conditions.
\newblock \emph{Management Science}, 42\penalty0 (10):\penalty0 1409--1419,
  1996.

\bibitem[Song et~al.(2025)Song, Xiao, and Zhang]{songxiaozhang2025}
Jing-Sheng Song, Li~Xiao, and Hanqin Zhang.
\newblock Optimal dual-sourcing inventory policies with order tracking:
  Backlogging and lost sales under uncertain lead times.
\newblock \emph{Operations Research}, 73\penalty0 (6):\penalty0 2867--2885,
  2025.

\bibitem[Srinivasan(2007)]{srinivasan2007}
Mahesh Srinivasan.
\newblock \emph{Optimal and Approximate Policies for Periodic Review Inventory
  Systems: The Case of Order Crossover and Multiple Supply Options}.
\newblock PhD thesis, The Pennsylvania State University, 2007.

\bibitem[Srinivasan et~al.(2011)Srinivasan, Novack, and Thomas]{srinivasan2011}
Mahesh Srinivasan, Robert Novack, and Douglas Thomas.
\newblock Optimal and approximate policies for inventory systems with order
  crossover.
\newblock \emph{Journal of Business Logistics}, 32\penalty0 (2):\penalty0
  180--193, 2011.

\bibitem[Stolyar and Wang(2022)]{stolyarwang2022}
Alexander~L. Stolyar and Qiong Wang.
\newblock Exploiting random lead times for significant inventory cost savings.
\newblock \emph{Operations Research}, 70\penalty0 (4):\penalty0 2496--2516,
  2022.

\bibitem[Stolyar et~al.(2026)Stolyar, Taneja, and Wang]{stolyar2026}
Alexander~L. Stolyar, Aditi Taneja, and Qiong Wang.
\newblock New policies exploiting randomness of lead times in inventory
  systems.
\newblock \emph{Naval Research Logistics}, 73\penalty0 (1):\penalty0 77--98,
  2026.

\bibitem[Wang and Disney(2017)]{wangdisney2017}
Xun Wang and Stephen~M. Disney.
\newblock Mitigating variance amplification under stochastic lead-time: The
  proportional control approach.
\newblock \emph{European Journal of Operational Research}, 256\penalty0
  (1):\penalty0 151--162, 2017.

\bibitem[Wei and Yang(2026)]{weiyang2026}
Zhaoxuan Wei and Yunfang Yang.
\newblock Blessing of order crossover for inventory systems.
\newblock Working paper, 2026.

\bibitem[Zalkind(1978)]{zalkind1978}
David Zalkind.
\newblock Order-level inventory systems with independent stochastic leadtimes.
\newblock \emph{Management Science}, 24\penalty0 (13):\penalty0 1384--1392,
  1978.

\bibitem[Zipkin(1986)]{zipkin1986}
Paul Zipkin.
\newblock Stochastic leadtimes in continuous-time inventory models.
\newblock \emph{Naval Research Logistics Quarterly}, 33\penalty0 (4):\penalty0
  763--774, 1986.

\bibitem[Zipkin(2000)]{zipkin2000}
Paul Zipkin.
\newblock \emph{Foundations of Inventory Management}.
\newblock McGraw-Hill, Boston, 2000.

\bibitem[Zipkin(2008)]{zipkin2008}
Paul Zipkin.
\newblock On the structure of lost-sales inventory models.
\newblock \emph{Operations Research}, 56\penalty0 (4):\penalty0 937--944, 2008.

\end{thebibliography}

\clearpage
\setcounter{section}{0}\renewcommand{\thesection}{EC.\arabic{section}}
\setcounter{equation}{0}\renewcommand{\theequation}{EC.\arabic{equation}}
\setcounter{table}{0}\renewcommand{\thetable}{EC.\arabic{table}}
\setcounter{figure}{0}\renewcommand{\thefigure}{EC.\arabic{figure}}
\setcounter{proposition}{0}\renewcommand{\theproposition}{EC.\arabic{proposition}}
\setcounter{lemma}{0}\renewcommand{\thelemma}{EC.\arabic{lemma}}
\setcounter{corollary}{0}\renewcommand{\thecorollary}{EC.\arabic{corollary}}
\makeatletter
\renewcommand{\theHsection}{EC.\arabic{section}}\renewcommand{\theHequation}{EC.\arabic{equation}}
\renewcommand{\theHtable}{EC.\arabic{table}}\renewcommand{\theHfigure}{EC.\arabic{figure}}
\@ifundefined{theHproposition}{}{\renewcommand{\theHproposition}{EC.\arabic{proposition}}}
\@ifundefined{theHlemma}{}{\renewcommand{\theHlemma}{EC.\arabic{lemma}}}
\@ifundefined{theHcorollary}{}{\renewcommand{\theHcorollary}{EC.\arabic{corollary}}}
\makeatother
\begin{center}
{\Large\bfseries E-Companion: Proofs and Supplementary Material}
\end{center}
\section{Proofs for the Main Text}\label{sec:ec-proofs}

This section proves the propositions of the main text, in their order, and
gives the supporting material that the main text cites with them. Several
proofs use the state in the coordinates
$(\IP,a)$ of Corollary~\ref{cor:longrun}, in which the transition
\eqref{eq:transy} reads
\begin{equation}\label{eq:transition}
\IP'=z-D,\qquad a_1'=(z-\IP)\,B_{0,0},\qquad
a_{k+1}'=a_k\,B_{k,0},\quad k=1,\dots,\rng-2 .
\end{equation}

\subsection{The Optimality Equation: Proof of Proposition~\ref{prop:acoe}}\label{sec:ec-acoe}

The equation is obtained by the vanishing-discount method
\citepec{schal1993,feinberg2012,feinbergliang2022}. The state space and the
costs are unbounded, so the discounted value functions, each less a constant,
must be bounded at each state uniformly in the discount factor, and they must
be equicontinuous. For inventory models with a fixed leadtime this is
verified by \citetec{huh2011} and \citetec{feinbergliang2022}. With random
leadtimes it rests on two estimates: a state far from the origin is costly
(Step~1), and one system can catch up with another at a bounded cost
(Step~3).

\proof{Proof of Proposition~\ref{prop:acoe}.}
We give the proof for real quantities and say at the end what changes for
integers.

\textit{Setup.} A state is $x=(\IP,a)$ with $a\ge0$, an action is an order
$q\ge0$, and the next state $x'$ is given by \eqref{eq:transition}. The stage
cost is that of the optimality equation \eqref{eq:acoe}. By \eqref{eq:received}
with $y=\IP-\sum_ka_k$ it is
\[
C(x,q)=\E\,G\Big(\IP-\sum_{k=1}^{\rng-1}a_kB_{k,\Lmin}+q\,\mathbf 1\{L=\Lmin\}\Big),
\]
where $L$ is the leadtime of the new order. This is the expected cost of the
period that lies $\Lmin$ periods ahead. The stage costs of the first $T$
periods therefore add up to the expected cost of periods $\Lmin+1$ to
$\Lmin+T$, and their long-run average is the long-run average cost. A policy
of the controller may use more than the history of the states and the orders.
But given that history and the current order, the next state has the same
distribution whatever else the controller knows (proof of
Corollary~\ref{cor:longrun}). So every policy induces a randomized policy on
the histories of states and orders with the same joint distribution of states
and orders, and it suffices to consider such policies. Write
$\|x\|=|\IP|+\sum_ka_k$ and $\bar b=\max(h,b)$, so that
$|\cbar(u)-\cbar(v)|\le\bar b\,|u-v|$. For a discount factor $\alpha\in(0,1)$
let $v_\alpha(x)$ be the infimum over all policies of the expected discounted
sum of the stage costs from $x$, and let $m_\alpha$ be the infimum of
$v_\alpha$ over the states.

\textit{Step 0: the discounted problem.} The stage cost is continuous and
nonnegative. Since $p_{\Lmin}>0$ and $G(u)\ge h\,(u-(\Lmin+1)\mu)$, it tends to
infinity with $q$, uniformly on bounded sets of states. The next state is a
continuous function of $(x,q)$ for each outcome of the demand and of the
indicators. So the model satisfies Assumption (W*) of \citetec{feinberg2012}:
the stage cost is bounded below and inf-compact on the product of every
compact set of states with the set of orders, and the transition is weakly
continuous. By their Theorem~2, $v_\alpha$ is lower semicontinuous, it
satisfies $v_\alpha(x)=\min_{q\ge0}\{C(x,q)+\alpha\,\E\,v_\alpha(x')\}$, a
stationary policy that attains the minimum is optimal, and $v_\alpha$ is the
pointwise limit of the value functions of the problems with finitely many
stages,
\begin{equation}\label{eq:Q}
V_0=0,\qquad V_n(x)=\min_{q\ge0}\big\{C(x,q)+\alpha\,\E\,V_{n-1}(x')\big\},\quad n\ge1 .
\end{equation}
Each $V_n$ is convex. For every outcome of the leadtimes and of the demand,
the argument of $G$ in the stage cost and the next state are linear in
$(x,q)$, and $G$ is convex. So if $V_{n-1}$ is convex, the function in braces
is jointly convex in $(x,q)$, and its minimum over the interval $q\ge0$ is
convex in $x$. Hence $v_\alpha$ is convex. It is finite: under the policy
that never orders, the net inventory of the period that is charged $t$ stages
from now is at most $\|x\|$ plus the demand of $t+\Lmin+1$ periods in absolute
value, so the stage cost is at most $\bar b\,(\|x\|+(t+\Lmin+1)\mu)$.

\textit{Step 1: the cost of $\rng+1$ stages.} There are constants $\epsilon>0$
and $A_0$ such that, from every state $x$ and under every policy,
the expected sum of the next $\rng+1$ stage costs is at least
$\epsilon\|x\|-A_0$. To see this, let the current period be period $0$. The
last of the $\rng+1$ stages is charged with period $\Lmax$, by which every
order that is outstanding at the start has been received. The orders placed
since are nonnegative, so the net inventory of that period is at least
$\IP-D(\Lmax+1)$, and the expected cost of the stage is at least
$h\,(\IP-(\Lmax+1)\mu)$. The first stage is charged with period $\Lmin$. With
probability $1-p_{\Lmin}$, which is positive because $\rng\ge2$, the new order
has not been received by then, and the net inventory is then at most
$\IP-\sum_ka_kB_{k,\Lmin}$. So the expected cost of the first stage is at
least $b\,(1-p_{\Lmin})\big(m_0\sum_ka_k-\IP\big)$, where
$m_0=\min_{k<\rng}m_{k,\Lmin}$ is positive because an outstanding order of age
$k<\rng$ can still have the longest leadtime. Stage costs are nonnegative, so
the sum of the $\rng+1$ stage costs is at least every convex combination of
the two bounds. For $\IP<0$ the second bound gives the claim. For $\IP\ge0$ a
convex combination does, with enough weight on the first bound to make the
coefficient of $\IP$ positive.

\textit{Step 2: the infimum is attained in a bounded set.} Let $x^0$ be the
state with $\IP=0$ and no outstanding orders. Under a base-stock policy from
$x^0$ the expected stage cost is at most some constant $A_1$ in every stage,
so $m_\alpha\le v_\alpha(x^0)\le A_1/(1-\alpha)$. Let $\alpha\ge\tfrac12$,
follow an optimal stationary policy from a state $x$, and write $C_t$ for
the cost of stage $t$ and $x_t$ for the state at stage $t$. By the
optimality equation of Step~0,
\[
v_\alpha(x)=\E\sum_{t\le\rng}\alpha^tC_t+\alpha^{\rng+1}\,\E\,v_\alpha(x_{\rng+1})
\ \ge\ \alpha^{\rng}(\epsilon\|x\|-A_0)+\alpha^{\rng+1}m_\alpha ,
\]
by Step~1, because stage costs are nonnegative. Since
$(1-\alpha^{\rng+1})\,m_\alpha\le(\rng+1)(1-\alpha)\,m_\alpha\le(\rng+1)A_1$,
\[
v_\alpha(x)-m_\alpha\ \ge\ 2^{-\rng}(\epsilon\|x\|-A_0)-(\rng+1)A_1
\qquad\text{whenever }\epsilon\|x\|\ge A_0 .
\]
Choose $\rho$ so that the right side is at least one when $\|x\|\ge\rho$, and
let $X_\rho$ be the set of states with $\|x\|\le\rho$, which is compact. Then
$m_\alpha$ is the infimum of $v_\alpha$ over $X_\rho$, and it is attained at
some $x_\alpha\in X_\rho$, because $v_\alpha$ is lower semicontinuous.

\textit{Step 3: catching up.} Let $x^1=(\IP^1,a^1)$ and $x^2=(\IP^2,a^2)$ be
two states with $\IP^1\le\IP^2$. Then
\begin{equation}\label{eq:catchup}
v_\alpha(x^1)\le v_\alpha(x^2)+\rng\,\bar b\,\Big(\IP^2-\IP^1+\sum_k|a^1_k-a^2_k|\Big).
\end{equation}
To see this, run a second system optimally from $x^2$, and let the first
system, which starts at $x^1$, imitate it: its first order is the first order
of the second system plus $\IP^2-\IP^1$, and its later orders are those of the
second system. Give the two systems the same demands, the same leadtimes for
the orders placed from now on, and the same leadtime for the outstanding
orders of each age. The controller of the first system observes all of this,
except the fate of an outstanding order of the second system at an age at
which it holds none itself. It draws the leadtime of such an order for
itself. It can then compute the orders of the second system, so the imitation
is an admissible randomized policy. The positions after ordering of the two
systems coincide from the first stage on. Their net inventories in a period
that is charged to some stage therefore differ only by what is outstanding of
the first order and of the orders outstanding at the start, that is, by at
most $\IP^2-\IP^1+\sum_k|a^1_k-a^2_k|$. And they differ only as long as one of
these orders is outstanding, which is the case in at most $\rng$ of the
charged periods, because the first charged period is period $\Lmin$ and each
of these orders has been received by period $\Lmax$. The cost of a period is
Lipschitz with constant $\bar b$, which gives \eqref{eq:catchup}.

\textit{Step 4: an upper bound.} Let $u_\alpha=v_\alpha-m_\alpha$ for
$\alpha\ge\tfrac12$. There is a constant $A_2$ that does not depend on
$\alpha$ such that
\[
u_\alpha(x)\le U(x):=A_2\,(1+\|x\|^2)
\]
at every state. If $\IP\le-\rho$, this is \eqref{eq:catchup} with $x^1=x$ and
$x^2=x_\alpha$. Otherwise, order nothing until the first stage $N$ at which
the position is at most $-\rho$, and follow an optimal policy from there. The
stage $N$ has a finite mean because $\mu>0$. Before stage $N$ the cumulative
demand is below $\IP+\rho$, so the net inventory of the period charged in
such a stage is at most $2\|x\|+\rho$ plus the demand of the $\Lmin+1$
periods that end with it, in absolute value. The event $\{N>t\}$ depends only
on the demands before stage $t$, so by Wald's identity the expected cost of
the stages before $N$ is at most
$\bar b\,(2\|x\|+\rho+(\Lmin+1)\mu)\,\E N$. Wald's identity for $D$ and for
$D^2$ also gives $\mu\,\E N\le\IP+\rho+\E D_{N-1}$ and
$\E D_{N-1}\le(\E N\,\E D^2)^{1/2}$, where $D_{N-1}$ is the demand of the last
stage before $N$. Hence $\E N$ and $\E D_{N-1}$ are at most a constant times
$1+|\IP|$. At stage $N$ the position lies between $-\rho-D_{N-1}$ and $-\rho$
and the outstanding orders are among those of $x$, so \eqref{eq:catchup} with
$x^2=x_\alpha$ bounds $v_\alpha$ at the state of stage $N$ by
$m_\alpha+\rng\,\bar b\,(2\rho+D_{N-1}+\|x\|)$. Stage costs are nonnegative
and $\alpha\le1$, so these bounds together give the claim.

\textit{Step 5: the optimality equation.} We verify the conditions of
\citetec[Theorem~3.2]{feinbergliang2022}. Their Assumption (W*) is that of
Step~0. Their Assumption (B) asks that some policy have a finite long-run
average cost from some state, which a base-stock policy has, and that
$v_\alpha(x)-m_\alpha$ be bounded in $\alpha\in[0,1)$ at every state. For
$\alpha\ge\tfrac12$ this is Step~4. For $\alpha<\tfrac12$ the difference is
at most $v_\alpha(x)$, which is at most the cost of the policy that never
orders at the discount factor $\tfrac12$. Their Assumption (EC) asks that
the functions $u_{\alpha_n}$ be equicontinuous along some sequence
$\alpha_n\uparrow1$, and that they be dominated by a function that is
integrable under every transition. They are dominated by $U$, and $\E\,U(x')$
is finite for every state and order because the demand has a finite
variance. They are equicontinuous for two reasons. Applied to two states with
the same position, in both orders, \eqref{eq:catchup} gives
$|u_\alpha(\IP,a^1)-u_\alpha(\IP,a^2)|\le\rng\,\bar b\sum_k|a^1_k-a^2_k|$. And
$u_\alpha$ is convex in the position, with values between zero and $U$, so
on every bounded set of states it is Lipschitz in the position with a
constant that does not depend on $\alpha$. The theorem then gives a sequence
$\alpha_n\uparrow1$ along which $u_{\alpha_n}$ converges, uniformly on compact
sets, to a continuous function $V^*\ge0$ that satisfies the equation in (a)
with the constant $\lim(1-\alpha)\,m_\alpha$. It also gives that this
constant is the optimal long-run average cost from every state, so that it
is $g^*$, and that every stationary policy that attains the minimum in (a)
is optimal. This is (a) and (b). In the paper $V^*$ denotes this limit. It
is the pointwise limit of the discounted value functions $v_{\alpha_n}$,
each less the constant $m_{\alpha_n}$.

\textit{Integer quantities.} The proof is the same on the integer states,
without convexity. The state space is then countable and discrete. Compact
sets of states are finite, so Assumption (W*) holds and the infimum in
Step~2 is attained. Every family of functions on a discrete space is
equicontinuous, so the bound of Step~4 is all that Assumption (EC) asks.
Steps~1 to~4 hold as written.
\Halmos
\endproof

\subsection{The Optimal Policy: Proofs for Section~\ref{sec:structure}}\label{sec:ec-structure}

The properties of $V^*$ are proved for the value functions $V_n$ of
\eqref{eq:Q}, by induction over the number of stages, and carried to $V^*$ by
two limits. By Step~0 of the proof of Proposition~\ref{prop:acoe}, the
discounted value function $v_\alpha$ is the pointwise limit of the functions
$V_n$. By Step~5 of that proof, $V^*$ is the pointwise limit of a
sequence of functions $v_\alpha$, each less a constant. So $V^*$ has every
property of the functions $V_n$ that is preserved by pointwise limits and by
subtracting constants.

\proof{Proof that $V^*$ is convex.}
With real quantities each $V_n$ is convex in $(\IP,a)$ (Step~0), and convexity
is preserved by the two limits. Convexity in $(\IP,a)$ and in $(y,a)$ are the
same property, because the two states differ by a linear change of
coordinates. \Halmos
\endproof

The next proposition holds at every range, and Section~\ref{sec:structure}
refers to it. At range two the oldest order is the only outstanding order,
and the proposition is then the middle inequality of \eqref{eq:chain} with
its equality case.

\begin{proposition}[The oldest order]\label{prop:oldest}
Let $e_{\rng-1}$ be the unit vector of the oldest age. At every range and at
every state,
\[
q^*(y+1,a)\ \le\ q^*(y,a+e_{\rng-1}),
\]
with equality if $y+1+\sum_{k<\rng-1}a_k+q^*(y,a+e_{\rng-1})\le0$.
\end{proposition}

\proof{Proof of Proposition~\ref{prop:oldest}.}
Write $Q(y,a,q)$ for the objective of the optimality equation
\eqref{eq:acoe}. It is continuous and tends to infinity with $q$, so it has
minimizers. A unit of stock and a unit of the order of age $\rng-1$ lead to
the same next state, because that order is counted in $y'$ at the next
decision whether or not it has arrived. So the two states have the same
continuation value at every $q$. The oldest order has still not been received
by the end of period $t+\Lmin$, the period that is charged, with probability
$m_{\rng-1,\Lmin}$. So the unit is in the net inventory of that period for
certain at $(y+1,a)$ and with probability $1-m_{\rng-1,\Lmin}$ at
$(y,a+e_{\rng-1})$. Write $\Delta_uG(x)=G(x+u)-G(x)$, and let
$Z=\sum_{k<\rng-1}a_k(1-B_{k,\Lmin})$ be what has been received by then of
the orders of ages below $\rng-1$. Then
$Q(y+1,a,q)-Q(y,a+e_{\rng-1},q)$ equals
\[
m_{\rng-1,\Lmin}\,\E\big[p_{\Lmin}\,\Delta_1G(y+Z+q)+(1-p_{\Lmin})\,\Delta_1G(y+Z)\big].
\]
This is nondecreasing in $q$ because $G$ is convex. Adding a
nondecreasing function of $q$ to an objective cannot raise its smallest or its
largest minimizer, which gives the inequality. Now let $q_1$ be the smallest
minimizer at $(y,a+e_{\rng-1})$ and suppose that
$y+1+\sum_{k<\rng-1}a_k+q_1\le0$. The demand is nonnegative and $\cbar$ is linear
on $(-\infty,0]$, so $G$ is linear there, and the difference above takes
the same value at every $q\le q_1$. So for $q<q_1$,
$Q(y+1,a,q)-Q(y+1,a,q_1)=Q(y,a+e_{\rng-1},q)-Q(y,a+e_{\rng-1},q_1)>0$, and the
smallest minimizer at $(y+1,a)$ is $q_1$. The argument for the largest
minimizer is the same. The proof does not use convexity in $q$, and an
ordering cost that depends on $q$ alone cancels from the difference. \Halmos
\endproof

\proof{Proof of Proposition~\ref{prop:chain}.}
At range two the state is $(\IP,a)$ with $\IP=y+a$. Put $s_1=a$, $s_2=\IP$ and
$s_3=\IP+q$, and for a function $U$ of the state write $U(s_1,s_2)$ for its
value at $a=s_1$ and $\IP=s_2$. By the period that is charged, the
outstanding order has not been received with probability
$m=m_{1,\Lmin}=p_{\Lmax}/\P(L\ge1)$, and the new order has been received with
probability $p_{\Lmin}$, so the stage cost is
\[
p_{\Lmin}m\,G(s_3-s_1)+p_{\Lmin}(1-m)\,G(s_3)
+(1-p_{\Lmin})\,m\,G(s_2-s_1)+(1-p_{\Lmin})(1-m)\,G(s_2).
\]
Each term is a convex function of one coordinate or of the difference of two,
so the sum is L$^\natural$-convex in $(s_1,s_2,s_3)$, and the constraint
$s_3\ge s_2$ defines an L$^\natural$-convex set. The next state is
$(a',\IP')=(s_3-s_2,\ s_3-D)$ if $B_{0,0}=1$, and $(0,\ s_3-D)$ if the order
is received in its own period, which is possible only with a shortest
leadtime of zero. Suppose that $U$ is L$^\natural$-convex. We claim that
$(s_2,s_3)\mapsto U(s_3-s_2,s_3)$ is
L$^\natural$-convex. By definition $\tilde U(u,v,w)=U(u-w,v-w)$ is L-convex in
three coordinates. L-convexity is preserved by permuting the coordinates and
by the reflection $x\mapsto-x$, and $U(s_3-s_2,s_3)=\tilde U(-s_2,0,-s_3)$,
which is the L-convex function $(x_1,x_2,x_3)\mapsto\tilde U(-x_1,-x_3,-x_2)$
evaluated at $(s_2,s_3,0)$. This proves the claim. The second next state gives
$U(0,s_3)$, a convex function of $s_3$ alone, because an L$^\natural$-convex
function is convex along each coordinate, and such a function is
L$^\natural$-convex in $(s_2,s_3)$. A shift by $D$, the expectation and the
mixture over $B_{0,0}$ preserve L$^\natural$-convexity, and a function that does not depend
on $s_1$ is L$^\natural$-convex in $(s_1,s_2,s_3)$ if it is in $(s_2,s_3)$. So
the stage cost plus $\alpha$ times the expectation of $U$ at the next state
is L$^\natural$-convex in $(s_1,s_2,s_3)$ for every $\alpha\in(0,1]$, and
minimizing it over $s_3$ leaves an L$^\natural$-convex function of
$(s_1,s_2)$ \citepec{murota2003,zipkin2008}.

With $U=V_{n-1}$ this is the induction step for the functions $V_n$ of
\eqref{eq:Q}, and $V_0=0$. So every $V_n$ is L$^\natural$-convex in
$(a,\,y+a)$, and so is $V^*$, because L$^\natural$-convexity is preserved by
pointwise limits and by subtracting constants. With $U=V^*$ and $\alpha=1$,
the objective of the optimality equation \eqref{eq:acoe} is
L$^\natural$-convex in $(s_1,s_2,s_3)$.
For an L$^\natural$-convex objective the smallest and the largest minimizers
$s_3^*(s_1,s_2)$ are nondecreasing in each argument and satisfy
$s_3^*(s_1+1,s_2+1)\le s_3^*(s_1,s_2)+1$. With $q^*=s_3^*-s_2$, a unit more of
stock is the move to $(s_1,s_2+1)$, so
$0\le s_3^*(s_1,s_2+1)-s_3^*(s_1,s_2)\le1$, which is
$q^*(y,a)-1\le q^*(y+1,a)\le q^*(y,a)$. A unit more on order is the move to
$(s_1+1,s_2+1)$, and
$s_3^*(s_1,s_2+1)\le s_3^*(s_1+1,s_2+1)\le s_3^*(s_1,s_2)+1$ is
$q^*(y+1,a)\le q^*(y,a+1)\le q^*(y,a)$. The equality case is that of
Proposition~\ref{prop:oldest}. \Halmos
\endproof

\subsubsection*{Beyond range two.} We checked the second
inequality of \eqref{eq:chain}, $q^*(y+1,a)\le q^*(y,a+e_k)$ with $e_k$ the
unit vector of age $k$, on the optimal policies of the test bed at range
three with mean demand 5 and 10 and at range four with mean demand 5, at the
states $(y,a+e_k)$ with a stationary probability above $10^{-9}$. No state
of these instances has more than one optimal order. For the order of age one
the inequality fails in \nv{chain.r3.age1.viol} of the \nv{chain.r3.n}
instances at range three and in \nv{chain.r4.age1.viol} of the
\nv{chain.r4.n} at range four, always by one unit, at states with a
stationary probability of at most \nv{chain.all.age1.mass}\% together in an
instance. At range four it also fails for the order of age two, in
\nv{chain.r4.age2.viol} instances. For the oldest order, which
Proposition~\ref{prop:oldest} covers, it fails at no state, and the first
and the third inequality of \eqref{eq:chain} fail at none of these states.

\subsection{Is a Weighted Position Sufficient? A Counterexample and a Numerical Check}\label{sec:ec-wp}

A policy uses a weighted position if its order is a nonincreasing function
$f_0$ of the position \eqref{eq:wp}. This subsection gives an exact
counterexample at range two, where the position is $y+\theta a$, and a
numerical check of the long-run problem. The example has the leadtime
uniform on $\{1,2,3\}$, the demand uniform on $\{0,1,2\}$, $h=1$ and $b=9$.

\subsubsection*{A counterexample that can be checked by hand.} Consider the
last decision period, in which only the cost of the last period is at stake.
The outstanding order $a$ is received in time with probability $2/3$. One
more unit ordered lowers the expected cost as long as the chance that it is
needed exceeds $h/(h+b)$, so the optimal order is the smallest $q\ge0$ with
\[
\tfrac13\,F(y+q)+\tfrac23\,F(y+a+q)\ \ge\ \tfrac9{10},
\]
where $F$ is the distribution function of the demand of two periods, with the
values $1/9$, $3/9$, $6/9$, $8/9$ and $1$ at $0,\dots,4$. Hence
\[
q^*(y,0)=(4-y)^+,\qquad q^*(y,a)=(3-y)^+\quad\text{for every }a\ge1 .
\]
An outstanding order lowers the optimal order by one unit, whatever its size,
because the new order has to cover the case that the outstanding one comes
too late. Now suppose that $q^*=f_0(y+\theta a)$ with $f_0$ nonincreasing. If
$\theta>0$, the position of the state $(0,a)$ exceeds 4 when $a$ is large, so
the order there is at most the order at the state $(4,0)$, which is zero. But
it is 3. If $\theta=0$, the states $(0,0)$ and $(0,1)$ have the same position,
and their orders are 4 and 3. So no weight reproduces the optimal orders.

\subsubsection*{The long-run problem.} The example is static. To see that its
conclusion is not an effect of the end of the horizon, we examined the
numerical solution of the long-run average problem at range two with a
shortest leadtime of one period, Poisson demand of mean 5, $h=1$ and $b=4$.
A nonincreasing function $f_0$ with $q^*(y,a)=f_0(y+\theta a)$ requires
$y+\theta a>y'+\theta a'$ whenever $q^*(y,a)<q^*(y',a')$, so each such pair
of states bounds $\theta$ from one side. We collected these bounds over all
pairs of states in a box, counting a pair only when the two sets of optimal
orders are disjoint. The bounds contradict each other under each of three
leadtime distributions. The weight $\theta$ would have to be above $0.875$
and below $0.800$ when the leadtime is uniform on $\{1,2,3\}$, above $0.833$
and below $0.600$ when these leadtimes have the probabilities $(.5,0,.5)$,
and above $0.833$ and below $0.750$ when they have $(.6,.2,.2)$. So no
weighted position reproduces the optimal orders of these numerical
solutions either.

\subsection{The Closed-Form Policy: Proof of Proposition~\ref{prop:weights}}\label{sec:ec-wli}

\proof{Proof of Proposition~\ref{prop:weights}.}
(a) Let the order of age $k$ have leadtime $L_k$, and let $L'$ be the leadtime
of an order placed now, in period $t$. The state says of $L_k$ only that
$L_k\ge k$, independently of $L'$. The older order is received in period
$t+L_k-k$ and the new one in period $t+L'$, so the older one is received no
later exactly when $L_k-k\le L'$, which has probability $\theta_k$. So the
expectation of $y+\sum_ka_k\mathbf 1\{L_k-k\le L'\}$ given the state is $P_t$.
(b) If $k\ge \rng$ then $L-k\le\Lmin\le L'$. At range two,
$\theta_1=\P(L'\ge L-1\mid L\ge1)$, and $L'<L-1$ only if $L=\Lmax$ and $L'=\Lmin$,
which has probability $p_{\Lmin}p_{\Lmax}/\P(L\ge1)$ given $L\ge1$.
(c) $\theta_k=\E\,\bar F(R_k)$, where $\bar F(j)=\P(L'\ge j)$ is nonincreasing and
$R_k$ has the distribution of $L-k$ given $L\ge k$. With hazard rates
$\lambda_i=\P(L=i\mid L\ge i)$, $\P(R_k\ge j)=\prod_{i=k}^{k+j-1}(1-\lambda_i)$,
which is nonincreasing in $k$ when the hazard rate is nondecreasing. So $R_k$
is stochastically nonincreasing in $k$ and $\theta_k$ is nondecreasing.
(d) The order placed $k$ periods after a given order, with leadtime $L'$,
arrives strictly before it if $k+L'<L$. So the expected number of such orders
is $\sum_{k\ge1}\P(L-L'\ge k+1)=\E(L-L'-1)^+$. Also
$(1-\theta_k)\P(L\ge k)=\P(L\ge k,\ L'<L-k)=\P(L-L'\ge k+1)$, because $L-L'\ge k+1$
and $L'\ge0$ imply $L\ge k$. Summing over $k$ gives \eqref{eq:chi}. \Halmos
\endproof

\subsubsection*{The level with integer demand.} With integer demand the
$\beta$-quantile of the shortfall moves in steps of one unit, while the
policy has targets that are not integers, so a fraction of a unit in the
level is visible in its cost. For integer demand only, we therefore smooth
$S^{\IP}$: we interpolate the distribution function of the shortfall
linearly between the integers, take its $\beta$-quantile, and add one half.
Let $F_Z$ be the distribution function of the shortfall $Z$ and let $k$ be
its $\beta$-quantile. Then
\[
S^{\IP}=k-\tfrac12+\frac{\beta-F_Z(k-1)}{F_Z(k)-F_Z(k-1)} .
\]
The interpolation is a convention that smooths the level. It is not the
minimizer of a cost.

\subsection{The Cost Charged to an Order: Proof of Proposition~\ref{prop:identity}}\label{sec:ec-model}

\proof{Proof of Proposition~\ref{prop:identity} and Corollary~\ref{cor:longrun}.}
Fix a period $s$ and a lag $\ell$, let $t=s+\ell$, and suppose that the order
of period $s$ governs period $t$. Write $\Lambda_t$ for the period whose order
governs period $t$, so that this is the event $\{\Lambda_t=s\}$.
The orders placed after $s$ have not been received by $t$, and the order of $s$
has. Let $z_s$ be the inventory position after ordering in period $s$. It
counts the net inventory and every order placed up to $s$ that is outstanding.
Hence
\[
I_t=z_s-\sum_{j<s}q_j\,\mathbf 1\{L_j\ge t-j+1\}-D[s,t],
\]
where $D[s,t]$ is the demand of periods $s$ to $t$ and the sum is over the
orders placed before $s$ that are still outstanding after the receipts of
period $t$. An order placed $k=s-j$ periods before $s$ appears in the sum if
$L_j\ge k+\ell+1$. Since $L_j\le\Lmax$ and $\ell\ge\Lmin$, this requires $k\le \rng-1$,
so the sum is $\sum_{k=1}^{\rng-1}a_{k,s}\mathbf 1\{L_{s-k}\ge k+\ell+1\}$.

Condition on the information at $s$ and on $q_s$. Leadtimes are independent
of the demands and of each other, and the policy uses only what has been
observed, so all that the information says about the leadtime of an
outstanding order of age $k$ is that it is at least $k$. The indicators above
are therefore independent Bernoulli variables with means $m_{k,\ell}$. The
event $\{\Lambda_t=s\}$ is $\{L_s\le\ell,\ L_{s+i}>\ell-i\text{ for }1\le
i\le\ell\}$; the condition with $i=\ell$ is on the order of period $t$ itself
and holds for certain unless the shortest leadtime is zero. It has
probability $\plag_\ell$ and is independent of the leadtimes of earlier orders
and of the demands. So
\[
\E\big[\cbar(I_t)\,\mathbf 1\{\Lambda_t=s\}\mid\text{information at }s\big]
=\plag_\ell\,\E\,\cbar\Big(z_s-\sum_ka_{k,s}B_{k,\ell}-D(\ell+1)\Big),
\]
and the sum over $\ell$ is $\Phi(z_s,a_s)$.

This proves the proposition. For the corollary, a period is governed by at
most one order, and a period $t\ge\Lmax+1$ by exactly
one, because the order of period $t-\Lmax$ has been received by then. A period
$t\le\Lmax$ may have no governing order. In that event none of the controller's
orders has been received, and the cost of the period does not depend on the
policy. Let $C_T$ be the expected cost of periods $1$ to $T$, and let $K_0$ be
the expected cost of the periods that have no governing order. The periods
governed by the orders of periods $1$ to $T$ include every governed period up
to $T$ and no period after $T+\Lmax$. Hence
\begin{equation}\label{eq:ec-sandwich}
C_T-K_0\ \le\ \sum_{s\le T}\E\,\Phi(z_s,a_s)\ \le\ C_{T+\Lmax},
\end{equation}
and dividing by $T$ gives the statement on long-run averages. For the
transition, the position before the next
order is $z-D$, the order just placed is outstanding with age one unless it is
received in its own period, which is the event $B_{0,0}=0$, and an
outstanding order of age $k$ is outstanding with age $k+1$ at the next decision
if its leadtime is at least $k+1$, which has conditional probability $m_{k,0}$,
independently across orders and of the demand. The pair $(\IP,a)$ determines
the conditional distribution of the next pair given the action, and the stage
cost depends on $(z,a)$ only. \Halmos
\endproof

The Markov decision problem of Corollary~\ref{cor:longrun} gives a recursion
with stage cost $\Phi$. It has the same long-run average cost as the
optimality equation \eqref{eq:acoe}. With a finite number of stages,
however, it charges some periods after the horizon, so it is not the
finite-horizon problem. This is why the proofs for
Section~\ref{sec:structure} use the stage cost of \eqref{eq:Q}.

\subsection{The Cost-Aware Policy: The Quadratic Model}\label{sec:ec-quad}

\subsubsection*{The local model.} Suppose the position and the outstanding
orders stay near values $\bar z$ and $\bar a$. With a density, to second
order, a change $(\delta z,\delta a)$ then changes the stage cost $\Phi$ of
Corollary~\ref{cor:longrun} by a linear term plus
\begin{equation}\label{eq:quadmodel}
\frac{b+h}{2}\sum_{\ell}\plag_\ell\,\E\Big[f_\ell\Big(\bar z-\sum_k\bar a_kB_{k,\ell}\Big)
\Big(\delta z-\sum_kB_{k,\ell}\,\delta a_k\Big)^2\Big].
\end{equation}
So each case counts with its probability times the density of its net
inventory at zero. With integer demand $G_\ell$ is linear between the
integers, and its second difference plays this role:
$G_\ell(n+1)-2G_\ell(n)+G_\ell(n-1)=(b+h)\,f_\ell(n)$, where $f_\ell(n)$ is the
probability that $D(\ell+1)$ equals $n$. We then use the same expression as
a model of the discrete curvature, with $f_\ell$ interpolated linearly
between the integers.

\subsubsection*{The weighted quadratic cost.} Write a case as $(\ell,B)$, with
$B$ the set of the ages of the missing orders. Let $\omega$ be a probability
distribution over the cases, and let $\omega_0>0$ be the probability that it
gives to the cases in which no order is missing. Replace the stage cost
$\Phi$ of Corollary~\ref{cor:longrun} by
\[
\Psi^\omega(z,a)=\sum_{\ell,B}\omega_{\ell,B}\,\E\Big(z-\sum_{k\in B}a_k-D(\ell+1)\Big)^2 ,
\]
let orders be real numbers of either sign, a negative order being a return
that takes that amount out of the stock when its leadtime has passed, and
drop the constraint $q\ge0$, which a rule that is linear in the state would
violate when the position is high. The second-order part of $\Psi^\omega$ is
$z^2-2z\,\psi^\top a+a^\top\Omega a$, with $\psi_k=\omega(k\in B)$,
$\Omega_{kk}=\psi_k$ and $\Omega_{jk}=\omega(j\in B,\,k\in B)$ for $j\ne k$. A
positive factor in front of $\Psi^\omega$ would multiply the value function
and leave the policy unchanged.

\begin{proposition}[Weighted cases]\label{prop:quadw}
With real-valued orders, without the constraint $q\ge0$ and at every range,
there are constants $\kappa^\omega\in(0,1]$ and
$\theta^\omega_1,\dots,\theta^\omega_{\rng-1}$, which depend on $\omega$ and on
the leadtime distribution only, and a target level $S$ such that the policy
\[
q=\kappa^\omega\Big(S-y-\sum_{k=1}^{\rng-1}\theta^\omega_ka_k\Big)
\]
minimizes the long-run average of $\Psi^\omega$ among all policies under
which the second moments of the state stay bounded.
\end{proposition}

If $\omega$ gives every case its probability,
$\omega_{\ell,B}=\plag_\ell\,\P(B_{\cdot,\ell}=B)$, then $\Psi^\omega$ is the
stage cost of $\cbar(x)=x^2$, $\psi_k=\bar m_k=\sum_\ell\plag_\ell\,m_{k,\ell}$
and $\Omega_{jk}=\sum_\ell\plag_\ell\,m_{j,\ell}m_{k,\ell}$ for $j\ne k$. We then
write $\kappa^{\mathrm Q}$ and $\theta^{\mathrm Q}_k$ for the constants. They
depend on the leadtime distribution only.

\subsubsection*{The recursion.} The state $X=(\IP,a_1,\dots,a_{\rng-1})$ moves
linearly by \eqref{eq:transition}, with the multiplicative noise $B_{k,0}$.
Index the coordinates of $X$ from zero. Let $F$ be the mean transition matrix
of $X$ at $z=0$, with $F_{1,0}=-m_{0,0}$, $F_{k+1,k}=m_{k,0}$ and all other
entries zero; let $g=(1,m_{0,0},0,\dots,0)^\top$ be the mean effect of $z$;
let $\tilde m=(0,\psi)$; let $\tilde\Omega$ be $\Omega$ bordered by a zero row
and column; and for a symmetric matrix $P$ let $N(P)$ be the diagonal matrix
with entries $P_{k+1,k+1}\,m_{k,0}(1-m_{k,0})$ at the positions of $a_k$,
$k\le \rng-2$, and let $n_0(P)=P_{11}\,m_{0,0}(1-m_{0,0})$, the term of the
order just placed. The number $n_0(P)$ is zero unless the shortest leadtime is
zero, because $m_{0,0}=\P(L\ge1)$. With $e_0$ the coordinate vector of $\IP$,
write $d(P)=1+g^\top Pg+n_0(P)$ and $c(P)=\tilde m-F^\top Pg+n_0(P)\,e_0$, and
let $P$ be the limit, from $P=0$, of
\begin{equation}\label{eq:riccati}
P\ \leftarrow\ F^\top PF+\tilde\Omega+N(P)+n_0(P)\,e_0e_0^\top-c(P)\,c(P)^\top/d(P),
\end{equation}
and let $K=c(P)/d(P)$. The coefficients of Proposition~\ref{prop:quadw} are
\[
\kappa^\omega=1-K_0=d(P)^{-1/2},\qquad\theta^\omega_k=1-\frac{K_k}{\kappa^\omega} .
\]
Under the uniform distribution at range three with a shortest leadtime of
at least one period, for example, $\kappa^{\mathrm Q}=0.857$ and
$\theta^{\mathrm Q}=(0.824,0.929)$.

The coefficients $\theta^{\mathrm Q}_k$ need not be nonnegative. With a
shortest leadtime of zero and $\P(L=0)=0.8$, $\P(L=2)=0.15$, $\P(L=6)=0.05$,
so that $\rng=6$, the recursion gives $\kappa^{\mathrm Q}=\nv{quad.ex.kappa}$
and
\[
\theta^{\mathrm Q}=(\nv{quad.ex.th1},\ \nv{quad.ex.th2},\ \nv{quad.ex.th3},\ \nv{quad.ex.th4},\ \nv{quad.ex.th5}).
\]
The coefficient of age three is negative: one more unit of that order raises
the optimal order slightly. The recursion converges to a residual below
$10^{-14}$, far smaller than this coefficient. On the test bed, at each
shortest leadtime that we use, the coefficients $\theta^{\mathrm Q}_k$ lie
between \nv{quad.bed.thmin} and one.

As the shortest leadtime grows, the coefficients of the cost-aware policy
move toward $\kappa^{\mathrm Q}$ and $\theta^{\mathrm Q}$: on the test bed
the largest difference falls from \nv{long.p2conv.5.max} at a shortest
leadtime of five periods to \nv{long.p2conv.500.max} at 500. We have not
proved that they converge, and Proposition~\ref{prop:long} covers rules with
fixed coefficients only.

\begin{lemma}\label{lem:memory}
Consider the system without demand, with the state $X=(\IP,a)$ and the
position after ordering, $z$, as the control, and let $\Phi_0(X,z)$ be its stage
cost, which is $\Psi^\omega$ without demand.

(a) There is a constant $C_0$ such that, on every path and for every sequence
of controls,
\[
\sum_{t\ge0}\|X_t\|^2\le C_0\Big(\|X_0\|^2+\sum_{t\ge0}z_t^2\Big).
\]

(b) $\Phi_0(X,z)\ge\omega_0z^2$. If $\omega$ gives every case its probability,
$\omega_0\ge\plag_{\Lmax-1}\ge p_{\Lmin}p_{\Lmax}^{\rng-1}>0$.
\end{lemma}

\proof{Proof.}
(a) By \eqref{eq:transition} without demand, $\IP_t=z_{t-1}$ and
$a_{1,t}=(z_{t-1}-\IP_{t-1})B_{0,0}$, and an outstanding order moves one age up in every
period, or leaves the state, until it is pooled after $\rng-1$ periods. So for
$t\ge \rng$ the state $X_t$ is a linear function of $z_{t-\rng},\dots,z_{t-1}$ alone,
with coefficients bounded by one in absolute value, and for $t<\rng$ it is such a
function of $X_0$ and $z_0,\dots,z_{t-1}$. (b) Each term of $\Phi_0$ is a
square, and in a case in which no order is missing it is $z^2$. With every
case at its probability, in a period of lag $\Lmax-1$ no outstanding order is
missing, because $m_{k,\Lmax-1}=0$ for $k\ge1$, so $\omega_0\ge\plag_{\Lmax-1}$.
Finally
$\plag_{\Lmax-1}=\P(L\le\Lmax-1)\prod_{i<\Lmax-1}\P(L>i)\ge p_{\Lmin}p_{\Lmax}^{\rng-1}$,
because the factors with $i<\Lmin$ equal one and each of the other $\rng-1$
factors is at least $p_{\Lmax}$. \Halmos
\endproof

\proof{Proof of Proposition~\ref{prop:quadw}.}
\textit{The optimality equation.} The stage cost $\Psi^\omega$ is quadratic
in $(z,a)$, and the transition is linear in $(X,z)$ with coefficients
$B_{k,0}$ that are independent of the state. Try $V(X)=X^\top PX+2v^\top X$ in
the average-cost optimality equation
$\bar J+V(X)=\min_z\{\Psi^\omega(z,a)+\E\,V(X')\}$. Since
$\E[B_{j,0}B_{k,0}]=m_{j,0}m_{k,0}$ for $j\ne k$ and $\E\,B_{k,0}^2=m_{k,0}$, the
function $\E\,V(X')$ is quadratic in $(X,z)$ with second-order part
$(FX+gz)^\top P(FX+gz)+X^\top N(P)X+n_0(P)(z-\IP)^2$, where the last term is the
variance of $a_1'=(z-\IP)B_{0,0}$. Minimizing over $z$ gives $z=K^\top X+k_0$,
and the second-order terms of the equation match if and only if $P$ is a
fixed point of \eqref{eq:riccati}. The demand enters only the terms of first
and zeroth order, through the means $\mu(\ell+1)$ of the cases, and hence only
$v$, $k_0$ and $\bar J$. With $z=K^\top X+k_0$ and $\IP=y+\sum_ka_k$,
\[
q=z-\IP=k_0-(1-K_0)\,y-\sum_k(1-K_0-K_k)\,a_k,
\]
which is the policy of the proposition with $S=k_0/\kappa^\omega$.

\textit{Convergence.} Consider the problem without demand. Its optimal cost
over $n$ periods is $X^\top P_nX$, where $P_n$ is the $n$th iterate of
\eqref{eq:riccati} from zero. It is nondecreasing in $n$, because $\Phi_0$ is
nonnegative, and at most the cost of the control $z\equiv0$, which is at most
$C_1\|X\|^2$ by Lemma~\ref{lem:memory}(a). So $P_n$ converges to a positive
semidefinite matrix $P$, which is a fixed point. The coefficient of $z^2$ in
each minimization is at least one, so the minimizers converge to the linear
function $K^\top X$, and
\begin{equation}\label{eq:quadfix}
X^\top PX=\Phi_0(X,K^\top X)+\E\,(A_\omega X)^\top P\,(A_\omega X),
\end{equation}
where $A_\omega$ is the random matrix that takes $X$ to the next state under
$z=K^\top X$, and $\omega$ in the subscript stands for the indicators
$B_{k,0}$ of the period.

\textit{Stability.} Under $z_t=K^\top X_t$, summing \eqref{eq:quadfix} along the
path and using Lemma~\ref{lem:memory}(b) gives
$\omega_0\sum_t\E\,z_t^2\le\sum_t\E\,\Phi_0(X_t,z_t)\le X_0^\top PX_0$, and
Lemma~\ref{lem:memory}(a) then gives $\sum_t\E\|X_t\|^2\le C_2\|X_0\|^2$ for a
constant $C_2$. Write $X_t=M_tX_0$, where $M_t$ is the product of the matrices
$A_\omega$ of the first $t$ periods. The matrix
$W=\sum_{t\ge0}\E\,M_t^\top M_t$ is then finite, with $I\preceq W\preceq C_2I$,
and it satisfies $W-\E\,A_\omega^\top WA_\omega=I$. Therefore
\[
\E\big[X_{t+1}^\top WX_{t+1}\mid X_t\big]=X_t^\top WX_t-\|X_t\|^2
\le\Big(1-\frac1{C_2}\Big)X_t^\top WX_t ,
\]
and $\E\|X_t\|^2$ decays geometrically. In particular the mean of the state
without demand tends to zero from every initial state, so the matrix
$\E\,A_\omega=F+gK^\top$ has all its eigenvalues inside the unit circle. With
demand, the state follows the same recursion with an input $\xi$ of finite
variance added in every period. The input is independent of the state of
the period, although with a shortest leadtime of zero it shares $B_{0,0}$
with $A_\omega$. By Minkowski's inequality in the norm of $W$,
$(\E\|X_{t+1}\|_W^2)^{1/2}\le(1-1/C_2)^{1/2}(\E\|X_t\|_W^2)^{1/2}+(\E\|\xi\|_W^2)^{1/2}$,
so the second moments stay bounded.

\textit{The terms of first order, and verification.} Matching the terms of
first order in the optimality equation gives a linear equation for $v$ with
the matrix $I-(F+gK^\top)^\top$, which is invertible by the last paragraph, and
the terms of zeroth order then give $\bar J$. So $(\bar J,V)$ solves the
optimality equation and $z=K^\top X+k_0$ attains its minimum. Let a policy keep
the second moments of the state bounded. Summing the inequality
$\bar J+V(X_t)\le\Psi^\omega(z_t,a_t)+\E[V(X_{t+1})\mid X_t,z_t]$ over $t<T$ and
dividing by $T$ shows that its long-run average cost is at least $\bar J$,
because $\E\,V(X_T)$ stays bounded. For the linear policy the inequality is
an equality, and its second moments are bounded by the last paragraph.

\textit{The damping factor.} Write $\kappa=1-K_0$, $d=d(P)$ and $c_0$ for
the entry of $c(P)$ at $\IP$. The column of $\IP$ in $F$ has the single entry
$-m_{0,0}$, in the row of $a_1$, and the stage cost has no term in $\IP$, so
$\tilde m$, $\tilde\Omega$ and $N(P)$ vanish at $\IP$. The entry of
\eqref{eq:riccati} at $(\IP,\IP)$ therefore reads
$P_{00}=m_{0,0}P_{11}-c_0^2/d$, with $c_0=m_{0,0}(P_{01}+P_{11})$ and
$d=1+P_{00}+2m_{0,0}P_{01}+m_{0,0}P_{11}$, so that
$1+P_{00}+m_{0,0}P_{01}=d-c_0$. The entry then gives
$d-c_0-1=c_0-c_0^2/d=c_0(d-c_0)/d$, that is, $(d-c_0)^2=d$. Since
$\kappa=1-c_0/d=(d-c_0)/d$, this is $\kappa^2d=1$, so $|\kappa|\le1$
because $d\ge1$. It remains to see that
$\kappa>0$. Without demand and with $S=0$, the policy orders
$q=-\kappa(y+\sum_k\theta^\omega_ka_k)$, and \eqref{eq:transy} gives, in
the coordinates $(y,a)$,
\begin{align*}
\E[y'\mid X]&=y+a_{\rng-1}+\sum_{k\le \rng-2}(1-m_{k,0})\,a_k+(1-m_{0,0})\,q,\\
\E[a_1'\mid X]&=m_{0,0}\,q,\qquad\E[a_{k+1}'\mid X]=m_{k,0}a_k,\quad k\le \rng-2 .
\end{align*}
Let $A$ be the matrix of this linear map. It is similar to $F+gK^\top$, so its
eigenvalues lie inside the unit circle, and $\det(I-A)>0$. In $I-A$ the
column of $y$ is $\kappa$ times the vector with the entry $1-m_{0,0}$ in the
row of $y$, $m_{0,0}$ in the row of $a_1$ and zeros elsewhere, and the column
of $a_k$ is $\theta^\omega_k$ times the column of $y$ plus a vector that
does not depend on the coefficients. Subtracting $\theta^\omega_k$ times
the column of $y$ from the column of $a_k$, for each $k$, does not change the
determinant, so $\det(I-A)=\kappa\det M$, where $M$ has that vector in the
place of the column of $y$ and the columns of $I-A$ at $\kappa=0$
elsewhere. The rows and columns of
$a$ in $M$ form a lower triangular matrix with unit diagonal, whose inverse
maps the first coordinate vector to $u=(u_1,\dots,u_{\rng-1})$ with
$u_k=\prod_{1\le i<k}m_{i,0}$. The complement of this block gives
$\det M=(1-m_{0,0})+m_{0,0}\big(\sum_{k\le \rng-2}(1-m_{k,0})u_k+u_{\rng-1}\big)$,
and the sum telescopes to $u_1=1$. So $\det M=1$ and
$\kappa=\det(I-A)>0$. With a shortest leadtime of at least one period,
$m_{0,0}=1$, $n_0(P)=0$ and $d(P)=1+g^\top Pg$. \Halmos
\endproof

With the constraint $q\ge0$ the policy can be used as
$q=\big(\kappa^\omega(S-y-\sum_k\theta^\omega_ka_k)\big)^+$. It is then a
heuristic, and neither optimality nor stability is proved for it.

\subsection{The Cost-Aware Policy: The Level}\label{sec:ec-level2}

The two positions $\bar z$ in the level \eqref{eq:S2} come from a model of
the net inventory. Suppose first that the mean position is $S^{\IP}$. In the
case $(O,\ell,B)$ the position is then $\bar z_O$ on average
(Section~\ref{sec:coef2}), the missing orders take $\mu|B|$ from it and the
demand takes $\mu(\ell+1)$. So the net inventory is
$\bar z_O-\mu|B|-\mu(\ell+1)$ on average. We model the net inventory as
normal around this value in each case. The cases have their probabilities
and one common variance, which is chosen so that the mixture has the
long-run variance of the net inventory under the rule. A different mean
position moves the whole mixture by the difference. Then
$\bar z(\kappa,\theta)$ is the mean position at which the mixture is
nonnegative with probability $\beta$, and $\bar z(1,1)$ is this position for
the base-stock policy, with $\kappa=1$ and every weight equal to one. The
next proposition gives the variance.

\begin{proposition}[Variance under a linear rule]\label{prop:linvar}
Under the rule $q=\kappa(S-y-\sum_k\theta_ka_k)$, without the constraint
$q\ge0$, the second moments of the state $(\IP,a)$ solve a linear equation.
If the rule is stable, the long-run variance of the net inventory is
$\sigma^2A+\mu^2A'$, where $A$ and $A'$ depend on the leadtime distribution
and on $(\kappa,\theta)$ only, and neither depends on $S$.
\end{proposition}

The first term is the demand variance that the rule passes on. The second is
the randomness of the arrivals: an outstanding order carries the mean demand,
and whether it has arrived is noise of that size.

The level \eqref{eq:S2} uses this variance: the common variance of the
mixture above is the variance of the proposition less
the variance of the means of the cases
(Section~\ref{sec:ec-p2alg}).

\proof{Proof of Proposition~\ref{prop:linvar}.}
Under the rule the position after ordering is
$z=\kappa S+(1-\kappa)\IP+\kappa\sum_k(1-\theta_k)a_k$, a linear function of
the state $X=(\IP,a)$. By \eqref{eq:transition} the next state is a linear
function of $(X,z)$ and of the demand, with the coefficients $B_{k,0}$, which
are independent of the state and of the demand. Taking expectations of the
next state and of its products of pairs of coordinates gives the first and
second moments of the next state as linear functions of those of the current
state: $\E[B_{k,0}]=m_{k,0}$, $\E[B_{j,0}B_{k,0}]=m_{j,0}m_{k,0}$ for $j\ne k$
and $\E[B_{k,0}^2]=m_{k,0}$. If the rule is stable these moments converge,
and their limit is the solution of the linear equation. By
Proposition~\ref{prop:identity} applied to the functions $x$ and $x^2$ in
place of $\cbar$, the long-run mean and second moment of the net inventory
are the long-run means of $\sum_\ell\plag_\ell\,\E(z-\sum_ka_kB_{k,\ell}-D(\ell+1))$
and of the same with the square, which are linear in the first and second
moments of $(z,a)$.

For the form of the variance, write $x$ for the state less its long-run
mean $\bar X$, and $A_\omega$ for the random matrix that takes the state to
the next one under the rule, at a demand equal to its mean. Then
$x'=A_\omega x+\xi_\omega+e\,(\mu-D)$, where $e$ is the coordinate vector
of $\IP$ and $\xi_\omega$ is the next state from $\bar X$ at the mean
demand less its mean. Three facts give the result. First, the indicators
$B_{k,0}$ multiply the outstanding orders and $B_{0,0}$ multiplies the order
just placed, so $\xi_\omega$ has the entry $(B_{0,0}-m_{0,0})\mu$ in the
row of $a_1$, the entry $(B_{k,0}-m_{k,0})\bar a_k$ in the row of $a_{k+1}$
and zero in the row of $\IP$; with a shortest leadtime of at least one period
$\xi_\omega=(A_\omega-\E A_\omega)\bar X$. These entries do not depend on
$S$: the mean order is $\mu$ by the balance of orders and demand, the mean
$\bar a_1$ of $a_1$ is $m_{0,0}$ times it, and the mean of $a_{k+1}$ is
$m_{k,0}$ times that of $a_k$. So $\xi_\omega$ is $\mu$ times a random
vector whose law depends on the leadtime distribution only, and $S$ enters
neither it nor the rest of the recursion. Second, the indicators and the demand of a period are
independent of the centered state of that period, which has mean zero, and
of each other, so no cross term survives in the second moments:
$\E\,x'x'^\top=\E\,A_\omega(\E\,xx^\top)A_\omega^\top+\mu^2Q+\sigma^2ee^\top$
for a matrix $Q$ that depends on the leadtime distribution only. The
long-run second moments of $x$ are therefore $\sigma^2$ times one matrix
plus $\mu^2$ times another. Third, the net inventory of a period at lag
$\ell$ is $z-\sum_ka_kB_{k,\ell}-D(\ell+1)$. Its variance adds to a quadratic
form in $x$ three terms: the variance of the indicators $B_{k,\ell}$, which
multiply outstanding orders and give a multiple of $\mu^2$ and a quadratic
form in $x$; the variance of the demand, $\sigma^2$ times the mean of
$\ell+1$; and the spread of the mean net inventory across the lags, in which
the mean of $z$ cancels and which is $\mu^2$ times a constant of the
leadtime distribution. \Halmos
\endproof

\subsection{The Cost-Aware Policy: Computation}\label{sec:ec-p2alg}

This subsection states the computation of Section~\ref{sec:costaware} in
full. Let $\beta=b/(b+h)$, let $\mu$ and $\sigma^2$ be the mean and the
variance of demand, let $f_\ell$ be the probability function of
$D(\ell+1)$, interpolated linearly between the integers, and let $S^{\IP}$
be the interpolated quantile of Section~\ref{sec:ec-wli}.

\textit{Cases and shares.} The cases $(O,\ell,B)$, their probabilities and
their weights \eqref{eq:caseweight} are those of
Section~\ref{sec:quadmodel}, with the position $\bar z_O$ of
Section~\ref{sec:coef2} for the coefficients $(\kappa,\theta)$ at hand, and
$\psi$ and $\Omega$ are the shares \eqref{eq:psi}. In a share, the factor
$\P(L\ge k)$ of an age that is outstanding in every case of its two sums
cancels.

\textit{The projection.} The shares $\psi_k$ and $\Omega_{jk}$ are taken
among different sets of cases, so they need not be the probabilities of one
distribution $\omega$, and the form $z^2-2z\,\psi^\top a+a^\top\Omega a$
need not be positive semidefinite. Let $D$ be the diagonal matrix with the
entries $\psi_k(1-\psi_k)$. On the ages with $0<\psi_k<1$ let
$\rho=D^{-1/2}(\Omega-\psi\psi^\top)D^{-1/2}$, which has a unit diagonal,
and let $\rho'$ be the matrix nearest to $\rho$ in the Frobenius norm among
the symmetric matrices with a unit diagonal and every eigenvalue at least
$\varepsilon=10^{-6}$. It is computed by alternating projections
\citepec{higham2002}. The recursion then uses
$\Omega'=\psi\psi^\top+D^{1/2}\rho'D^{1/2}$ in place of $\Omega$, where the
second term has zero rows and columns at the other ages. The diagonal stays $\psi$, and
$\Omega'=\Omega$ when every eigenvalue of $\rho$ is at least $\varepsilon$
already. At range two there is one outstanding order, the form is the
quadratic cost of one distribution, and the projection changes nothing. On
the test bed it changes the form in \nv{p2.proj.n} of the \nv{p2.def.n}
pairs of an instance and a shortest leadtime: \nv{p2.proj.0}, \nv{p2.proj.1}
and \nv{p2.proj.2} at shortest leadtimes of zero, one and two periods, and
none at three, four or five.

\textit{The map.} With $(\psi,\Omega')$, the recursion
\eqref{eq:riccati} is iterated from $P=0$ until no entry changes by more than
$10^{-13}$, and its gains give new coefficients $(\kappa',\theta')$ as in
Section~\ref{sec:ec-quad}.

\textit{The fixed point.} The iteration starts from the coefficients
of the closed-form policy \eqref{eq:policy}. From the current
coefficients $(\kappa,\theta)$ it moves to
$(\kappa,\theta)+\lambda\big((\kappa',\theta')-(\kappa,\theta)\big)$, with
$\lambda=1/2$ at the start. Whenever the largest difference between
$(\kappa',\theta')$ and $(\kappa,\theta)$ exceeds $0.9$ times that of the
previous round, $\lambda$ is halved, down to $0.01$. The iteration stops
when this difference is below $10^{-7}$, that is, when the map returns the
coefficients to within $10^{-7}$, or after 400 rounds.

\textit{The level.} The variance of the net inventory under the rule is
$\sigma^2A+\mu^2A'$ of Proposition~\ref{prop:linvar}, with the second moments
found by iterating their recursion until no entry changes by more than
$10^{-13}$. The net inventory is modeled as a mixture over the cases, with
the mean of Section~\ref{sec:ec-level2} in each case and a common normal
deviation. Its variance $v$ is the variance of the net inventory less the
variance of the means of the cases, and the computation uses
$\max(v,0.05\,\sigma^2)$. The position $\bar z(\kappa,\theta)$ at which the
mixture is nonnegative with probability $\beta$ is found by bisection. The
same computation with $\kappa=1$ and every weight one gives $\bar z(1,1)$,
and \eqref{eq:S2} gives the level.

\textit{Checks.} The policy is defined where these steps succeed: the
iteration reaches a fixed point and $v$ is positive. Under the condition of
Proposition~\ref{prop:exist} a fixed point exists, the recursion converges at
it and the rule is stable in mean square, but the iteration need not find it.
On the \nv{bed.n} instances of Section~\ref{sec:numerical}, at each of the
six shortest leadtimes that it uses, each step succeeds. The condition of
the proposition holds in each of the \nv{p2.exist.n} runs with an age whose
status is uncertain, where $(S^{\IP}+1)/\mu$ exceeds the sum in the
proposition by at least \nv{p2.exist.margin3}. The fixed point is reached in
at most \nv{p2.def.rounds} rounds. At the fixed point the recursion
\eqref{eq:riccati} converges in at most \nv{p2.def.ricits} iterations, and
the coefficients reproduce themselves to $10^{-7}$. The second-moment map of
the rule has a spectral radius of at most \nv{p2.def.radius}. The variance
$v$ is at least \nv{p2.def.raw} times $\sigma^2$, so the floor
$0.05\,\sigma^2$ never binds.

\textit{Uniqueness.} The fixed point is not unique in general. A search
from a grid of starting points found one fixed point in all but
\nv{p2.scan.multi} of \nv{p2.scan.n} runs, the test bed at shortest
leadtimes of zero, one and two periods and the \nv{zero.n} published
instances of Section~\ref{sec:published}. The \nv{p2.scan.multi} have a
shortest leadtime of zero and three fixed points each. The policy uses the
fixed point that the iteration from the closed-form coefficients reaches.

The weights of the cost-aware policy are not restricted to the unit
interval. On the test bed they are negative in \nv{p2.negw.1} of the
\nv{p2.negw.1.n} instances, down to $\nv{p2.negw.1.min}$, in \nv{p2.negw.0}
at a shortest leadtime of zero, down to $\nv{p2.negw.0.min}$, and in none at
two to five periods. All are at range four under the two-point
distributions with $0.9$ or $0.7$ on the shortest leadtime.

\subsection{The Cost-Aware Policy: Existence of the Fixed Point}\label{sec:ec-exist}

The cost-aware coefficients are a fixed point of the map of
Section~\ref{sec:coef2}. The next proposition gives a condition under which
the map has one.

\begin{proposition}[Existence]\label{prop:exist}
Let the demand of a period take every nonnegative integer value with
positive probability, and consider the numbers
$\frac12\sqrt{\rng-1-k}\;\P(L<k)$ over the ages $k<\rng$ with
$0<\P(L\ge k)<1$. If $S^{\IP}+1$ exceeds $\mu$ times the sum of the two
largest of them, or of the one if there is only one, then the map has a
fixed point. Its damping factor lies in $(0,1]$, and without the constraint
$q\ge0$ the rule with these coefficients keeps the second moments of the
state bounded.
\end{proposition}

The condition asks that the level $S^{\IP}$ is not too low compared with the
mean demand. It guarantees that the cases in which no order is missing keep
a positive weight for every set of coefficients that the map can return. It
holds in every instance of Section~\ref{sec:numerical}.

\proof{Proof of Proposition~\ref{prop:exist}.}
\textit{The map.} Call an age $k<\rng$ uncertain if $0<\P(L\ge k)<1$, let
$U$ be the set of the uncertain ages, and write $x_k=\kappa(1-\theta_k)$. An
age at which an order is outstanding for certain adds nothing to $\bar z_O$ of
Section~\ref{sec:coef2}, so
$\bar z_O=S^{\IP}+\mu\sum_{k\in U}x_k\big(1\{k\in O\}-\P(L\ge k)\big)$. The
case weights, the shares $(\psi,\Omega)$ and the projected $\Omega'$ of
Section~\ref{sec:ec-p2alg} therefore depend on the coefficients only through
$x=(x_k)_{k\in U}$. Let $\Gamma(x)$ be the vector of the gains $K_k$,
$k\in U$, that the recursion \eqref{eq:riccati} returns from
$(\psi,\Omega')$. The new coefficients satisfy $\kappa'(1-\theta'_k)=K_k$, so
coefficients reproduce themselves exactly when $x$ is a fixed point of
$\Gamma$. If $U$ is empty there is nothing to reproduce. Let
$\eta_k=\frac12\sqrt{\rng-1-k}$. We show that $\Gamma$ is continuous on the
box of the $x$ with $-\eta_k\le x_k\le1+\eta_k$ for $k\in U$ and maps the box
into itself, so that Brouwer's theorem gives a fixed point.

\textit{Step 1: the recursion for a positive form.} Let $\Lambda$ be the
symmetric matrix with the first row $(1,-\psi^\top)$ and the block
$\Omega'$, so that the stage cost without demand is
$\Phi_0=(z,a^\top)\,\Lambda\,(z,a^\top)^\top$. Suppose that
$\Lambda\succeq\omega_0ee^\top$ for some $\omega_0>0$, with $e$ the first
coordinate vector. Then $\Phi_0\ge\omega_0z^2$, the coefficient of $z^2$ is
one, and $|\Omega'_{jk}|\le(\psi_j\psi_k)^{1/2}\le1$, because
$\Omega'-\psi\psi^\top$ is positive semidefinite. These are the only
properties of $\Psi^\omega$ that the proof of Proposition~\ref{prop:quadw}
uses, with this $\omega_0$ in Lemma~\ref{lem:memory}(b). So the iterates $P_n$
of \eqref{eq:riccati} from zero increase to a fixed point $P$, the gains
$K=c(P)/d(P)$ give a rule that is stable in mean square, and
$1-K_0\in(0,1]$.

\textit{Step 2: continuity in the form.} Let $(\psi,\Omega')$ vary over a set
on which $\Lambda\succeq\omega_0ee^\top$ for one $\omega_0>0$. Each $P_n$ is
continuous in $(\psi,\Omega')$, because a step of \eqref{eq:riccati} divides
by $d(P)\ge1$, so $P=\sup_nP_n$ is lower semicontinuous as a quadratic
form. For a linear rule $K'$ that is stable in mean square, its cost without
demand, $J_{K'}$, is finite and linear in $(\psi,\Omega')$, because the path
under $K'$ does not depend on them. $P_n$ is at most the cost of $K'$ over
$n$ periods, so $P\le J_{K'}$, and summing \eqref{eq:quadfix} along the path
of the rule $K$ gives $P=J_K$. So $P$ is a minimum of continuous functions,
hence upper semicontinuous and therefore continuous, and so is
$K=c(P)/d(P)$.

\textit{Step 3: bounds on the gains.} For $k\le\rng-2$ the gain of $a_k$ is
$K_k=\big(\psi_k-m_{k,0}(Pg)_{k+1}\big)/d(P)$, and
$K_{\rng-1}=\psi_{\rng-1}/d(P)\in[0,1]$. By the Cauchy--Schwarz inequality
in the semi-inner product of $P$,
$|(Pg)_{k+1}|\le(P_{k+1,k+1}\,g^\top Pg)^{1/2}$, and
$d(P)\ge1+g^\top Pg$. Since $\sqrt s/(1+s)\le\frac12$ for $s\ge0$, this gives
$-\frac12P_{k+1,k+1}^{1/2}\le K_k\le1+\frac12P_{k+1,k+1}^{1/2}$. The entry
$P_{k+1,k+1}$ is the optimal cost without demand from one unit outstanding at
age $k+1$ and $\IP=0$. Under $z\equiv0$ that unit costs
$\Omega'_{jj}=\psi_j\le1$ in each period in which it is outstanding at an
age $j\le\rng-1$, and there are at most $\rng-1-k$ such periods. So
$P_{k+1,k+1}\le\rng-1-k$, and $-\eta_k\le K_k\le1+\eta_k$.

\textit{Step 4: the case weights.} Let $x$ lie in the box and $k\in U$. Take
$O$ to contain $k$, every age at which an order is outstanding for certain, and each other
$j\in U$ exactly when $x_j\ge0$. Every term of $\bar z_O-S^{\IP}$ with $j\ne
k$ is then nonnegative, and the term of $k$ is at least
$-\mu\eta_k\P(L<k)$. With two ages $j$ and $k$ in $O$ and the others chosen
in the same way,
$\bar z_O\ge S^{\IP}-\mu\big(\eta_j\P(L<j)+\eta_k\P(L<k)\big)$. By the
condition of the proposition both bounds exceed $-1$. In a period of lag
$\Lmax-1$ no outstanding order is missing, and $\plag_{\Lmax-1}>0$. So the
case $(O,\Lmax-1,\emptyset)$ has the weight
$\P(O)\,\plag_{\Lmax-1}f_{\Lmax-1}(\bar z_O)>0$, because $f_\ell$ is positive
on $(-1,\infty)$ when the demand of a period takes every nonnegative integer
value. Hence every share $\psi_k$ and $\Omega_{jk}$ is a ratio of continuous
functions of $x$ with a denominator that is bounded below on the box. That
case is in the denominator of $\psi_k$ and not in its numerator, so
$\psi_k<1$, and by compactness $\psi_k\le\bar\psi$ on the box for a constant
$\bar\psi<1$.

\textit{Step 5: the projection.} The matrix $\rho'$ has a unit diagonal and
eigenvalues of at least $\varepsilon$. For a position $z$ and orders $a$,
\[
\Phi_0=(z-\psi^\top a)^2+a^\top D^{1/2}\rho'D^{1/2}a
\ge(z-\psi^\top a)^2+\varepsilon\sum_k\psi_k(1-\psi_k)\,a_k^2 .
\]
By the Cauchy--Schwarz inequality the sum is at least $(\psi^\top a)^2/C_3$
with $C_3=\sum_k\psi_k/(\varepsilon(1-\psi_k))\le(\rng-1)/(\varepsilon(1-\bar\psi))$,
and $(z-\psi^\top a)^2+(\psi^\top a)^2/C_3\ge z^2/(1+C_3)$, because the left
side is smallest when $\psi^\top a=zC_3/(1+C_3)$. So
$\Lambda\succeq\omega_0ee^\top$ on the whole box, with
$\omega_0=1/\big(1+(\rng-1)/(\varepsilon(1-\bar\psi))\big)$. The nearest
point of a closed convex set depends continuously on the point, so
$\Omega'$ is continuous in $x$ wherever the shares that are positive stay
bounded away from zero. Let $\psi_k$ tend to zero. The entries of
$D^{1/2}\rho'D^{1/2}$ in the row of $k$ are at most
$\big(\psi_j(1-\psi_j)\psi_k(1-\psi_k)\big)^{1/2}$ and tend to zero. The
numerator of $\Omega_{jk}$ counts cases with order $k$ missing among those
with both orders outstanding, which carry at most $1/\P(L\ge j)$ times the weight of
the numerator of $\psi_k$, and the denominators are bounded below by Step 4.
So $|\Omega_{jk}-\psi_j\psi_k|\le C_4\psi_k$ for a constant $C_4$, and the
entries of $\rho$ between $k$ and the ages whose shares stay bounded away
from zero tend to zero. In the limit $\rho$ is block diagonal. The nearest
point of a block-diagonal matrix is block diagonal, because changing the
signs of the coordinates of one block maps the set onto itself and fixes
the matrix, and its block of the ages whose shares stay positive is the
nearest point of that block alone. Since the nearest point moves by no more
than its argument, $\Omega'$ converges to its value at the limit, and it is
continuous on the box.

\textit{Step 6: conclusion.} By Steps 4 and 5 the pair $(\psi,\Omega')$ is a
continuous function of $x$ on the box, and $\Lambda\succeq\omega_0ee^\top$
there with one $\omega_0>0$. By Steps 1 and 2 the gains are continuous
functions of $x$, so $\Gamma$ is continuous, and by Step 3 it maps the box
into itself. Brouwer's theorem gives a fixed point, and Step 1 gives
$\kappa^\beta\in(0,1]$ and the stability of the rule there. If $U$ is
empty, Steps 4 and 5 hold with the box a single point. \Halmos
\endproof

The proposition is about the coefficients. The level \eqref{eq:S2} and the
constraint $q\ge0$ are applied to them afterwards, as in
Sections~\ref{sec:level2} and~\ref{sec:costaware}.

\subsection{The Lower Bounds: Proof of Proposition~\ref{prop:lbk}}\label{sec:ec-lb}

\begin{proposition}[Jensen bounds]\label{prop:jensen}
Let $\bar M_\ell=\sum_{k\ge1}\P(L\ge k+\ell+1)$ be the expected number of
orders that are missing from a period with lag $\ell$, and let
$\LB=\min_z\sum_\ell\plag_\ell\,G_\ell\big(z-\mu\bar M_\ell\big)$. Then
$\LB\le\LB_0\le g^*$, with $\LB_0$ of \eqref{eq:lb0}.
\end{proposition}

\proof{Proof of Proposition~\ref{prop:jensen}.}
Consider a policy with finite long-run average cost, and let $C_T$ be its
expected cost over periods $1$ to $T$, so that $C_T\le MT$ for some $M$ and all
$T$. Let $\bar z_T$ and $\bar a_{k,T}$ be the averages over $s\le T$ of
$\E\,z_s$ and $\E\,a_{k,s}$. By \eqref{eq:ec-sandwich} and Jensen's inequality
for the convex $\Phi$,
\begin{equation}\label{eq:ec-jensen}
\frac{C_{T+\Lmax}}T\ \ge\ \Phi(\bar z_T,\bar a_T)\ \ge\ \varphi(\bar a_T).
\end{equation}
where $\varphi(a)=\min_z\Phi(z,a)$, which is convex.
We show that $\bar a_T\to\bar a$. In a period of lag $\Lmax-1$ every order placed
before the governing one has been received, so $m_{k,\Lmax-1}=0$ for every $k$
and
\[
\Phi(z,a)\ge\plag_{\Lmax-1}\,G_{\Lmax-1}(z)\ge\plag_{\Lmax-1}\min(h,b)\,\big(|z|-\Lmax\mu\big),
\]
where $\plag_{\Lmax-1}\ge p_{\Lmin}p_{\Lmax}^{\rng-1}>0$. With \eqref{eq:ec-sandwich} this gives
$\sum_{s\le T}\E|z_s|\le M'T$ for some $M'$ and all $T$. From one period to the
next the position after ordering falls by at most the demand,
$z_{s+1}\ge z_s-D_s$. So $\E\,z^+_{s+j}\ge\E\,z^+_s-j\mu$ and
$\E\,z^-_{s-j}\ge\E\,z^-_s-j\mu$. Fix $\varepsilon\in(0,\mu)$ and let $n$ be
the integer part of $\varepsilon T/(2\mu)$. If $\E\,z^+_T\ge\varepsilon T$, then
$\E|z_{T+j}|\ge\varepsilon T/2$ for $0\le j\le n$. If
$\E\,z^-_T\ge\varepsilon T$, then $\E|z_{T-j}|\ge\varepsilon T/2$ for
$0\le j\le n$. In both cases
$\sum_{s\le T+n}\E|z_s|\ge\varepsilon^2T^2/(4\mu)$, which contradicts the linear
bound when $T$ is large. So $\E|z_T|/T\to0$. By flow balance,
$z_T=\IP_1+\sum_{s\le T}q_s-\sum_{s<T}D_s$, so the average of $\E\,q_s$ over
$s\le T$ tends to $\mu$. An order is sized before its leadtime is drawn, so
$\E\,a_{k,s}=\E\,q_{s-k}\,\P(L\ge k)$, and therefore
$\bar a_{k,T}\to\mu\,\P(L\ge k)=\bar a_k$. Finally $\varphi$ is convex and
finite, hence continuous, and \eqref{eq:ec-jensen} gives
$\liminf_TC_T/T\ge\varphi(\bar a)=\LB_0$.
This proves $\LB_0\le g^*$. For $\LB\le\LB_0$, apply
Jensen's inequality to the variables $B_{k,\ell}$ in \eqref{eq:Phi}:
$\Phi(z,\bar a)\ge\sum_\ell\plag_\ell G_\ell\big(z-\sum_k\bar a_km_{k,\ell}\big)$,
and $\sum_k\bar a_km_{k,\ell}=\mu\sum_k\P(L\ge k+\ell+1)=\mu\bar M_\ell$. \Halmos
\endproof

\subsubsection*{The general program.} Let $k\ge0$. Let $O_j$
indicate that the order of age $j$ at $s$ is outstanding, so that the $O_j$
are independent with $\P(O_j=1)=\P(L\ge j)$. Let $j_0=2$ when the shortest
leadtime is at least one period, so that the order of age one is always
outstanding, and $j_0=1$ when it is zero. Write
$O=(O_{j_0},\dots,O_{\rng-1})$ for the pattern of the outstanding orders at
$s$, and $d=(d_1,\dots,d_k)$ for values of $D_{s-1},\dots,D_{s-k}$. For
$j+j_0\le i\le \rng-1$ let
$U^{(j)}_i=\mathbf 1\{L_{s-i}\ge i-j\}$ indicate that the order now of age $i$
was outstanding when the order now of age $j$ was placed. The variables of the
program are $z(O,d)$ and $\tilde q_j(U,d_{j+1},\dots,d_k)\ge0$ for
$j=1,\dots,\rng-1$. Let
\[
\bar z_0(d)=\sum_O\P(O)\,z(O,d),\qquad
\bar z_j(d_{j+1},\dots,d_k)=\sum_{d'}p(d')\,\bar z_{j-1}(d_{j+1},\dots,d_k,d')
\]
for $j=1,\dots,k$,
and $\zeta_j(d)=\bar z_j(d_{j+1},\dots,d_k)-\bar z_{j+1}(d_{j+2},\dots,d_k)+d_{j+1}$ for $j<k$,
$\zeta_j=\mu$ for $j\ge k$. Then
\begin{equation}\label{eq:lbkgen}
\begin{aligned}
&\LB_k=\min\,\sum_{O,d}\P(O)p(d)\Phi\big(z(O,d),\bar a(O,d)\big),\\
&\bar a_j(O,d)=O_j\,\bar q_j(O,d),\qquad
\bar q_j(O,d)=\sum_U\P\big(U^{(j)}=U\mid O\big)\,\tilde q_j(U,\cdot),
\end{aligned}
\end{equation}
where $\tilde q_j(U,\cdot)$ stands for $\tilde q_j(U,d_{j+1},\dots,d_k)$,
subject to $\sum_U\P(U)\,\tilde q_j(U,d_{j+1},\dots,d_k)=\zeta_j(d)$ for every $j$
and, when $k\ge1$, $\bar z_0(d)-\bar z_1(d_2,\dots,d_k)+d_1\ge0$. With $k=0$
the current order has mean $\mu$ and this constraint is absent. Here
$\P(U^{(j)}_i=1\mid O_i=1)=1$ and
$\P(U^{(j)}_i=1\mid O_i=0)=\P(i-j\le L<i)/\P(L<i)$. At range two with a
shortest leadtime of at least one period there is no pattern and
\eqref{eq:lbkgen} is \eqref{eq:lbk2}. For $k\ge1$, $\LB_k$ is the value of
\eqref{eq:lbkgen}, and $\LB_0$ is the bound \eqref{eq:lb0}, which is the
program with $k=0$ and without the pattern. The proof below also shows that
the value does not decrease when the pattern is added to the conditioning or
when the demand values are replaced by a finer partition.

\proof{Proof of Proposition~\ref{prop:lbk}.}
Take a stationary policy. By Corollary~\ref{cor:longrun} its cost is
$\E\,\Phi(z_s,a_s)$, and conditional Jensen with respect to the pattern and the
last $k$ demands gives
$\E\,\Phi(z_s,a_s)\ge\E\,\Phi\big(\E[z_s\mid O,d],\E[a_s\mid O,d]\big)$. It remains
to show that the conditional means are feasible for \eqref{eq:lbkgen}.

Positions. The pattern is a function of leadtimes and is independent of the
demands, so $\sum_O\P(O)\,\E[z_s\mid O,d]=\E[z_s\mid d]=\bar z_0(d)$. The position
$z_{s-j}$ is chosen before $D_{s-j},\dots,D_{s-1}$, so
$\E[z_{s-j}\mid d]=\E[z_{s-j}\mid D_{s-j-1}=d_{j+1},\dots,D_{s-k}=d_k]$, which by
stationarity is $\bar z_0$ with its last $j$ arguments averaged out, that is, $\bar z_j$.

Orders. By \eqref{eq:link}, $\E[q_{s-j}\mid d]=\zeta_j(d)$ for $j<k$, and $\mu$
for $j\ge k$, when none of its three terms sees the conditioning demands. The
order $q_{s-j}$ is independent of the leadtimes $L_{s-j},\dots,L_{s-1}$, hence
of $O_i$ for $i\le j$, and $a_{j,s}=q_{s-j}O_j$, so
$\E[a_{j,s}\mid O,d]=O_j\,\E[q_{s-j}\mid O_{>j},d]$. Of the order now of age
$i>j$ the controller knew at $s-j$ only what it had observed: the arrival
period if $L_{s-i}<i-j$, and otherwise that $L_{s-i}\ge i-j$. Given
$U^{(j)}_i=1$, $q_{s-j}$ is independent of $L_{s-i}$ and hence of $O_i$; given
$U^{(j)}_i=0$, $O_i=0$. So
$\E[q_{s-j}\mid O_{>j},d]=\sum_U\P(U\mid O)\,\tilde q_j(U,\cdot)$ with
$\tilde q_j(U,\cdot)=\E[q_{s-j}\mid U^{(j)}=U,d]\ge0$, which does not depend on
$d_1,\dots,d_j$, and averaging over $U$ gives $\zeta_j$. The last constraint is
$\E[q_s\mid d]\ge0$.

For a policy that is not stationary, apply the same steps to the averages over
$s\le T$ of the conditional means. The pattern and the last $k$ demands have
the same distribution in every period, so Jensen's inequality bounds the
average cost by the objective at these averages. The relations above link the
orders of period $s$ to the positions of earlier periods, so the averages
satisfy the constraints up to terms of order $\E|z_T|/T$. By the proof of
Proposition~\ref{prop:jensen} these terms vanish, and the averages are bounded
because $\sum_{s\le T}\E|z_s|\le M'T$. When the demand takes finitely many
values, or is grouped into finitely many cells as in our computations, a limit
point of the averages is feasible, and since $\Phi$ is continuous its
objective is at most the limit inferior of the average cost.

Monotonicity. Given a feasible point for $k+1$, average it over $d_{k+1}$. The
maps $\bar z_j$ commute with that average, so the result is feasible for $k$, and by
Jensen's inequality its objective is no larger. The same argument covers the
pattern and a finer partition of the demand values; with a partition, $d_{j+1}$
in $\zeta_j$ is the mean of the cell. With $k=0$ and no pattern the program is
$\min_z\Phi(z,\bar a)$. \Halmos
\endproof

\subsection{Computations and Published Policies: Details for Section~\ref{sec:numerical}}\label{sec:ec-pub}

\subsubsection*{Computations.} The Poisson distribution is cut at a tail of
$10^{-12}$ and renormalized, the state space and the orders are bounded, the
iterations stop at a tolerance, and the linear programs of the bounds are
solved by the HiGHS solver at its default tolerances. The optimal cost is
computed by relative value iteration, and the Bellman residuals of the last
step bracket it to within $10^{-9}$ of its value. The levels of the two
policies use the quantile $S^{\IP}$ interpolated between the integers
(Section~\ref{sec:ec-wli}). The linear program of the bound $\LB_k$ groups
the demand values into nine cells at range three and into six at range four.

\subsubsection*{The bounds on the state.} The inventory position and the orders
are kept in a box that grows with the mean demand, and we checked its effect in three
ways. On four instances that we recomputed in a larger box, the optimal cost
and the costs of the base-stock policy and of the closed-form policy changed
by less than $10^{\nv{box.rel.pow}}$ in relative terms. Under the optimal
policy the mean order equals the mean demand to
$10^{\nv{top.lr.balance.pow}}$ in relative terms in every instance of the
test bed, so no stock is lost at the box. And the stationary probability
of each face of the box, that is, of the lowest or the highest inventory
position or of an outstanding order at its cap, is at most
$10^{\nv{edge.lmin1p.pow}}$ under each of the two policies at shortest
leadtimes of one period or more. At a shortest leadtime of zero it reaches
$\nv{edge.lmin0.p2.mant}\times10^{\nv{edge.lmin0.p2.pow}}$ under the
cost-aware policy, whose orders are large there, and its cost is then less
accurate. We have not derived an error bound for the box, so these checks
show that its effect is small where we looked.

\subsubsection*{The searches of the published policies.} We search the level
of the base-stock policy, the fraction and the level of the proportional
policy, the multiple and the level of the generalized base-stock policy, the
two levels of the two-level policy, the weights of the weighted index policy,
and for PIPELINE its centering parameter and the number of periods that its
program plans. The searches are local, from several starting points. At
range four with mean demand 10 and 20 the searches of the proportional
policy, the weighted index policy and the best rule compare their candidates
by simulation with common random numbers, and the parameters that they
return are then priced exactly. We ran both kinds of search on
\nv{simcal.n} smaller instances. For the weighted index policy the rule
found by simulation costs \nv{simcal.wi.mean} percentage points more on
average than the rule found by the exact search, and \nv{simcal.wi.max} at
most. For the proportional policy and the best rule the difference is at
most \nv{simcal.pout.max} and \nv{simcal.lin.max} percentage points.

\subsubsection*{Adaptations.} Three adaptations should be stated. The
weighted index policy was proposed for lost sales, with its level held at the
best base-stock level and its weights searched on a grid, and we keep the
level there as its authors do. They round the order up, and we round it to
the nearest integer. The paper that proposes PIPELINE also proposes a policy
called ADAPTIVE, an adaptive version of the generalized base-stock policy,
which we do not price. PIPELINE was proposed for a model in which
each unit has its own leadtime. Its authors name the case in which the units
of an order travel together as an extension, and we apply the program to it
unchanged. The program can have more than one optimal first order. This
happens, for example, when the critical ratio $b/(b+h)$ is a value of the
leadtime distribution function, as it is in a part of our test bed. The cost
of the policy depends on which one is placed, and its authors do not say. We
price four choices, each with its own parameters: the one that the solver of
the program returns, the smallest, the largest and their midpoint. We report
the one that costs least. Its orders are large in some states and reach the
bound on the orders of our state space. Where they do so with a probability
above $10^{-6}$, in \nv{pipe.all.atbound} instances, its cost is estimated by
a simulation without a bound. At range four with mean demand 20 the state
space is too large to price every candidate on it. There the smallest, the
largest and the midpoint are searched by simulation, each with its own
parameters. The three best candidates are priced by a simulation of
$5\times 10^7$ periods that is paired with the base-stock policy, and we
report the cheapest of the three. Choosing among them on the sample that
prices them is slightly in favor of PIPELINE. On \nv{pipesim.n} smaller instances this procedure gives
a gap that is \nv{pipesim.mean} percentage points above the gap of the exact
procedure on average and \nv{pipesim.max} at most.

\subsubsection*{The policies that use the net inventory.} The generalized and
the two-level base-stock policy use the net inventory by itself, which is not
a function of the state $(\IP,a)$. We price them exactly on a larger state
where it can be computed, and elsewhere by a simulation of $10^8$ periods
that is paired with the base-stock policy. The standard error of a simulated
gap is at most \nv{net.r2.se} percentage points at range two,
\nv{net.r3.se} at range three and \nv{net.r4.se} at range four. The exact
pricing covers \nv{net.exact.n} instances: all of range two and
\nv{net.r3.exact} of range three, at mean demand 5. On these the simulated
gaps differ from the exact ones by at most \nv{net.exact.diff.max}
percentage points.

\section{An Upper Bound and Long Shortest Leadtimes}\label{sec:long}

The bounds of Section~\ref{sec:bounds} are computed instance by instance, and they do not show how
the gap of a policy depends on the leadtime distribution. This section
gives two results that do. The first is an upper bound on the cost of a
policy in closed form. The second is a bound on the gap that falls at least as fast
as one over the shortest leadtime as it grows while the range stays fixed.

Both results use a bound that is weaker than $\LB_0$. Let
$\bar M_\ell=\sum_{k\ge1}\P(L\ge k+\ell+1)$ be the expected number of
orders that are missing from a period with lag $\ell$, and let
\[
\LB=\min_z\ \sum_\ell\plag_\ell\,G_\ell\big(z-\mu\bar M_\ell\big).
\]
This bound replaces the number of missing orders by its mean, and
$\LB\le\LB_0$ (Proposition~\ref{prop:jensen}). Both
results compare the cost of a policy with $\LB$, and both rest on one
observation. By the cost identity, the cost of a policy is the
long-run average of $\sum_\ell\plag_\ell\,\E\,G_\ell(\Ipre_\ell)$, with the
received stock $\Ipre_\ell$ of \eqref{eq:received-lag}. The bound $\LB$ is
the same expression with $\Ipre_\ell$ replaced by its mean, at the best mean
position. So the excess of a policy over $\LB$ is the cost of the variance
of $\Ipre_\ell$, and it is at most half the largest curvature of $G_\ell$
times this variance. The curvature of $G_\ell$ is $(b+h)f_\ell$
(Section~\ref{sec:quadmodel}), and the density $f_\ell$ of the demand of
$\ell+1$ periods falls like one over $\sigma\sqrt{\ell}$. So the excess
shrinks as the shortest leadtime grows, while $\LB$ itself grows like the
square root of the shortest leadtime.

The results are for policies with a level. A policy with a level is a family
of stationary policies $\pi_S$, one for each level $S$ and possibly
randomized, such that $\pi_S$ orders at the state $(y,a)$ what $\pi_0$
orders at $(y-S,a)$. Raising the level by one unit then raises the stock,
the position and every $\Ipre_\ell$ by one unit in the long run, and leaves
the orders as they were. The base-stock policy and the rules
\eqref{eq:policy} and \eqref{eq:policy2}, with their coefficients held
fixed, are policies with a level. With integer quantities the levels differ
by integers. We assume that the state has a stationary distribution with
finite second moments and take the cost and the variances under it. So
$g^\pi$ is the long-run average cost when the system starts from this
distribution, and $g^\pi_*$ is its smallest value over the levels. Such a
cost is never below $g^*$, as shown below under ``Stationary costs.''

We measure the variance of the received stock under $\pi$ by its average
over the lags,
\[
\Var_\pi\Ipre=\sum_\ell\plag_\ell\,\Var_\pi\Ipre_\ell ,
\]
and compare it with
\[
\Var^{0}\Ipre=\sum_\ell\plag_\ell\sum_{k=1}^{\rng-1}\bar a_k^2\,m_{k,\ell}\big(1-m_{k,\ell}\big),
\]
the variance when the position is constant and every outstanding order has
its mean size $\bar a_k$. This is the variance that the bound $\LB_0$ keeps,
and $\Var^{0}\Ipre\le\Var_\pi\Ipre$.
Finally, $\sup f_\ell$ is the largest value of the density of $D(\ell+1)$,
or of its probability function with integer demand.

\begin{proposition}[An upper bound]\label{prop:ub}
Let $\pi$ be a policy with a level. With real quantities, and a demand with
a density,
\[
g^\pi_*\le\LB+\frac{b+h}2\sum_\ell\plag_\ell\,(\sup f_\ell)\,\Var_\pi\Ipre_\ell ,
\]
and with integer quantities the same holds with $\Var_\pi\Ipre_\ell+1$ in
place of $\Var_\pi\Ipre_\ell$. Under the base-stock policy
\[
\Var\Ipre_\ell=\sum_k\P(L\ge k+\ell+1)\,\big(\sigma^2+\P(L\le k+\ell)\,\mu^2\big).
\]
\end{proposition}

With the base-stock policy the proposition gives an interval for the optimal
cost in closed form: $g^*$ lies between $\LB$ and the right side. The
interval is wide when the shortest leadtime is short, because the bound uses
the largest value of the density at every received stock. Its width falls
like one over $\sqrt{\Lmin}$, while $\LB$ grows like $\sqrt{\Lmin}$.

The second result is for a long shortest leadtime. If
$\Lmin\ge\max(1,\rng-1)$, every order in the state is outstanding for
certain, so the state tells the controller nothing about which orders are
late. The transition and the probabilities $m_{k,\ell}$ and $\plag_\ell$
then depend on the leadtimes only through $\ell-\Lmin$. So a rule on the
state is the same rule at every such shortest leadtime, and only the demand
of the governed periods grows with $\Lmin$. Orders still cross, but a policy
can choose only how it sizes its orders.

\begin{proposition}[Long shortest leadtimes]\label{prop:long}
Fix the distribution of $L-\Lmin$, and let $\Lmin\to\infty$.
Let the demand be integer with span one or have a bounded density, and let
$\pi$ be a policy with a level whose rule does not change with $\Lmin$.
The variance $\Var_\pi\Ipre$ does not depend on $\Lmin$ once
$\Lmin\ge\max(1,\rng-1)$, and
\[
\limsup_{\Lmin\to\infty}\Lmin\,\frac{g^\pi_*-\LB_0}{\LB_0}\le\frac{\Var_\pi\Ipre-\Var^{0}\Ipre}{2\sigma^2},
\]
with $\tfrac54$ added to the numerator when quantities are integers. For
these shortest leadtimes, under the base-stock policy,
$\Var\Ipre-\Var^{0}\Ipre=\sigma^2\sum_\ell\plag_\ell\bar M_\ell$.
\end{proposition}

The proposition gives a certificate. Since $\LB_0\le g^*\le g^\pi_*$, the ratio
$(g^\pi_*-\LB_0)/\LB_0$ bounds the optimality gap of the policy, and it
falls at least as fast as one over the shortest leadtime. Its constant is the variance of
the received stock under the policy in excess of the variance that $\LB_0$
keeps, divided by twice the variance of the demand of one period. For the
base-stock policy with real quantities the constant is half the expected
number of orders that are missing from a period, and it does not depend on
the demand distribution. The reason is that every base-stock order equals
the demand of one period, so each missing order adds the variance of one
demand to that of the received stock. A policy that smooths its orders passes
less of the demand variance on to them and lets its position vary instead,
and its constant can be smaller. The proposition holds with integer
quantities. It
applies to the closed-form policy as it stands, because its weights and its
damping factor depend on differences of leadtimes only once
$\Lmin\ge\max(1,\rng-1)$. It is a statement about the best level, and
Table~\ref{tab:long} below shows that the certificates at the level
\eqref{eq:S} fall at the same rate.

The proofs are at the end of this section. The excess of
$\E\,G_\ell(\Ipre_\ell)$ over $G_\ell(\E\,\Ipre_\ell)$ is an integral of
$f_\ell$ against a nonnegative weight whose total is half the variance of
$\Ipre_\ell$. The local limit theorem makes $f_\ell$ nearly constant where
the weight lies, and it gives the order $\sigma\sqrt{\Lmin}$ of the bounds.

Proposition~\ref{prop:long} has three limitations. First, as the shortest
leadtime grows, every policy with a level and a finite variance is optimal
to first order, the base-stock policy included, and the proposition tells
the policies apart only by the constant of its bound. Second, this constant
is mostly due to the bound $\LB_0$ and not to the policy, for the reason
given in Section~\ref{sec:bounds}. Third, the proposition says nothing about
a short shortest leadtime. When $\rng>\Lmin+1$, the state shows which orders
are late, a new order can arrive before a late one, and the certificates
are those of Section~\ref{sec:bounds}.

\subsubsection*{Numerical results.} Table~\ref{tab:long} gives the
certificates of Proposition~\ref{prop:long} at shortest
leadtimes from 5 to 500 periods on the instances of the test bed, and the
gaps to the optimal cost at 20 and 100 periods. The certificates use the
bound $\LB_0$, which is a convex program in one variable. The stationary
distribution of a rule is the same at each of these shortest leadtimes, as
noted above, so we compute it once on the state space of the test bed. The
optimal cost is computed as in Section~\ref{sec:longrun}.

\begin{table}
\TABLE
{Certificates and Optimality Gaps at Long Shortest Leadtimes (Mean and Largest, Percent)\label{tab:long}}
{\begin{tabular}{@{}l@{\quad}c@{\ \ }c@{\ \ }c@{\ \ }c@{\ \ }c@{\quad}c@{\ \ }c@{}}
\hline\up
 & \multicolumn{4}{c}{Certificate, at $\Lmin$ of} & & \multicolumn{2}{c}{Gap, at $\Lmin$ of}\\
 & 5 & 20 & 100 & 500 & Limit & 20 & 100\\ \hline\up
Range 2 (135 instances)\\%
\quad Base stock, best level & \nv{long.r2.bs.5.mean} (\nv{long.r2.bs.5.max}) & \nv{long.r2.bs.20.mean} (\nv{long.r2.bs.20.max}) & \nv{long.r2.bs.100.mean} (\nv{long.r2.bs.100.max}) & \nv{long.r2.bs.500.mean} (\nv{long.r2.bs.500.max}) & \nv{long.r2.bs.lim.mean} (\nv{long.r2.bs.lim.max}) & \nv{long.dp.r2.bs.20.mean} (\nv{long.dp.r2.bs.20.max}) & \nv{long.dp.r2.bs.100.mean} (\nv{long.dp.r2.bs.100.max})\\%
\quad Closed form \eqref{eq:policy} & \nv{long.r2.p1.5.mean} (\nv{long.r2.p1.5.max}) & \nv{long.r2.p1.20.mean} (\nv{long.r2.p1.20.max}) & \nv{long.r2.p1.100.mean} (\nv{long.r2.p1.100.max}) & \nv{long.r2.p1.500.mean} (\nv{long.r2.p1.500.max}) & \nv{long.r2.p1.lim.mean} (\nv{long.r2.p1.lim.max}) & \nv{long.dp.r2.p1.20.mean} (\nv{long.dp.r2.p1.20.max}) & \nv{long.dp.r2.p1.100.mean} (\nv{long.dp.r2.p1.100.max})\\%
\quad Cost aware \eqref{eq:policy2} & \nv{long.r2.p2.5.mean} (\nv{long.r2.p2.5.max}) & \nv{long.r2.p2.20.mean} (\nv{long.r2.p2.20.max}) & \nv{long.r2.p2.100.mean} (\nv{long.r2.p2.100.max}) & \nv{long.r2.p2.500.mean} (\nv{long.r2.p2.500.max}) & \nv{long.r2.p2.lim.mean} (\nv{long.r2.p2.lim.max}) & \nv{long.dp.r2.p2.20.mean} (\nv{long.dp.r2.p2.20.max}) & \nv{long.dp.r2.p2.100.mean} (\nv{long.dp.r2.p2.100.max})\\%
Range 3 (135 instances)\\%
\quad Base stock, best level & \nv{long.r3.bs.5.mean} (\nv{long.r3.bs.5.max}) & \nv{long.r3.bs.20.mean} (\nv{long.r3.bs.20.max}) & \nv{long.r3.bs.100.mean} (\nv{long.r3.bs.100.max}) & \nv{long.r3.bs.500.mean} (\nv{long.r3.bs.500.max}) & \nv{long.r3.bs.lim.mean} (\nv{long.r3.bs.lim.max}) & \nv{long.dp.r3.bs.20.mean} (\nv{long.dp.r3.bs.20.max}) & \nv{long.dp.r3.bs.100.mean} (\nv{long.dp.r3.bs.100.max})\\%
\quad Closed form \eqref{eq:policy} & \nv{long.r3.p1.5.mean} (\nv{long.r3.p1.5.max}) & \nv{long.r3.p1.20.mean} (\nv{long.r3.p1.20.max}) & \nv{long.r3.p1.100.mean} (\nv{long.r3.p1.100.max}) & \nv{long.r3.p1.500.mean} (\nv{long.r3.p1.500.max}) & \nv{long.r3.p1.lim.mean} (\nv{long.r3.p1.lim.max}) & \nv{long.dp.r3.p1.20.mean} (\nv{long.dp.r3.p1.20.max}) & \nv{long.dp.r3.p1.100.mean} (\nv{long.dp.r3.p1.100.max})\\%
\quad Cost aware \eqref{eq:policy2} & \nv{long.r3.p2.5.mean} (\nv{long.r3.p2.5.max}) & \nv{long.r3.p2.20.mean} (\nv{long.r3.p2.20.max}) & \nv{long.r3.p2.100.mean} (\nv{long.r3.p2.100.max}) & \nv{long.r3.p2.500.mean} (\nv{long.r3.p2.500.max}) & \nv{long.r3.p2.lim.mean} (\nv{long.r3.p2.lim.max}) & \nv{long.dp.r3.p2.20.mean} (\nv{long.dp.r3.p2.20.max}) & \nv{long.dp.r3.p2.100.mean} (\nv{long.dp.r3.p2.100.max})\\%
Range 4 (135 instances)\\%
\quad Base stock, best level & \nv{long.r4.bs.5.mean} (\nv{long.r4.bs.5.max}) & \nv{long.r4.bs.20.mean} (\nv{long.r4.bs.20.max}) & \nv{long.r4.bs.100.mean} (\nv{long.r4.bs.100.max}) & \nv{long.r4.bs.500.mean} (\nv{long.r4.bs.500.max}) & \nv{long.r4.bs.lim.mean} (\nv{long.r4.bs.lim.max}) & \nv{long.dp.r4.bs.20.mean} (\nv{long.dp.r4.bs.20.max}) & \nv{long.dp.r4.bs.100.mean} (\nv{long.dp.r4.bs.100.max})\\%
\quad Closed form \eqref{eq:policy} & \nv{long.r4.p1.5.mean} (\nv{long.r4.p1.5.max}) & \nv{long.r4.p1.20.mean} (\nv{long.r4.p1.20.max}) & \nv{long.r4.p1.100.mean} (\nv{long.r4.p1.100.max}) & \nv{long.r4.p1.500.mean} (\nv{long.r4.p1.500.max}) & \nv{long.r4.p1.lim.mean} (\nv{long.r4.p1.lim.max}) & \nv{long.dp.r4.p1.20.mean} (\nv{long.dp.r4.p1.20.max}) & \nv{long.dp.r4.p1.100.mean} (\nv{long.dp.r4.p1.100.max})\\%
\quad Cost aware \eqref{eq:policy2} & \nv{long.r4.p2.5.mean} (\nv{long.r4.p2.5.max}) & \nv{long.r4.p2.20.mean} (\nv{long.r4.p2.20.max}) & \nv{long.r4.p2.100.mean} (\nv{long.r4.p2.100.max}) & \nv{long.r4.p2.500.mean} (\nv{long.r4.p2.500.max}) & \nv{long.r4.p2.lim.mean} (\nv{long.r4.p2.lim.max}) & \nv{long.dp.r4.p2.20.mean} (\nv{long.dp.r4.p2.20.max}) & \nv{long.dp.r4.p2.100.mean} (\nv{long.dp.r4.p2.100.max})\down\\%
 
\hline
\end{tabular}}
{Each entry is the mean over the instances of the range, with the largest
value in parentheses. A certificate is the cost of the policy in excess of
the bound $\LB_0$, in percent of the bound, and it bounds the gap, which is
the excess over the optimal cost in percent of that cost. The policies have the levels of Section~\ref{sec:lmin}, and the
base-stock policy its best level. ``Limit'' is the bound
of Proposition~\ref{prop:long} on the limit of $\Lmin$ times the
certificate at the best level, for real quantities,
$100(\Var_\pi\Ipre-\Var^{0}\Ipre)/(2\sigma^2)$; with integer
quantities the bound adds $\tfrac54$ to the numerator.}
\end{table}

\subsubsection*{Stationary costs.} For a stationary policy $\pi$, possibly
randomized, under which the state has a stationary distribution $\nu$ with
finite second moments, $g^\pi$ above is the stationary
mean of the cost of a period. It is at least $g^*$. By
Proposition~\ref{prop:acoe}(a), at every state and for every order $q\ge0$,
$g^*+V^*(y,a)\le\E\,G(\Ipre_{t+\Lmin})+\E\,V^*(y',a')$. Integrate over
$\nu$ and over the orders of $\pi$. The next state then has the law $\nu$
again, and the mean of $G(\Ipre_{t+\Lmin})$ is the stationary mean of the
cost of a period. The function $V^*$ is a pointwise limit of the functions
$u_\alpha$ of the proof of Proposition~\ref{prop:acoe}, so
$V^*\le A_2(1+\|x\|^2)$ by Step~4 of that proof, and its integral is
finite. It cancels, and $g^*\le g^\pi$. Under $\nu$ the mean order is $\mu$,
since $\IP'=z-D$, and an order is sized before its leadtime is drawn, so
$\E\,a_k=\bar a_k$ and $\E\,\Ipre_\ell=\E\,z-\mu\bar M_\ell$, as in the
proof of Proposition~\ref{prop:jensen}. Given the state, the indicators
$B_{k,\ell}$ add $\sum_ka_k^2m_{k,\ell}(1-m_{k,\ell})$ to the variance of
$\Ipre_\ell$, and $\E\,a_k^2\ge\bar a_k^2$, so
$\Var^{0}\Ipre\le\Var_\pi\Ipre$.

\subsubsection*{The excess over the mean.} Let $G$ be convex with
$G''=(b+h)f$ for a density $f$, and let $Y$ have a finite variance.
Taylor's formula with integral remainder at $\E\,Y$ gives
\begin{equation}\label{eq:ec-excess}
\E\,G(Y)-G(\E\,Y)=(b+h)\int f(u)\,\varpi(u)\,du,\qquad
\varpi(u)=\E(Y-u)^+-(\E\,Y-u)^+ .
\end{equation}
The weight $\varpi$ is nonnegative by Jensen's inequality, its integral is
$\tfrac12\Var Y$, and it vanishes outside the smallest interval that contains
the values of $Y$. For $\rho\ge1$ its integral over $|u-\E\,Y|>\rho$ is at
most $\E[(Y-\E\,Y)^2;|Y-\E\,Y|>\rho]$. With integer demand $G$ is linear
between the integers, with a kink of size $(b+h)f(n)$ at each integer $n$,
where $f$ is the probability function, and \eqref{eq:ec-excess} holds with
the integral replaced by the sum over the integers. A direct computation
gives
\begin{equation}\label{eq:ec-lattice}
\sum_n\varpi(n)=\tfrac12\Big(\Var Y+\E\,\langle Y\rangle\big(1-\langle Y\rangle\big)
-\langle\E\,Y\rangle\big(1-\langle\E\,Y\rangle\big)\Big),
\end{equation}
where $\langle\cdot\rangle$ is the fractional part. So the sum is at most
$\tfrac12\Var Y$ when $Y$ is integer, and at least $\tfrac12(\Var Y-\tfrac14)$
in every case. The bound on the tail holds for the sum as well.

\proof{Proof of Proposition~\ref{prop:ub}.}
Take expectations under the stationary distribution. By
Corollary~\ref{cor:longrun},
$g^\pi=\E\,\Phi(z,a)=\sum_\ell\plag_\ell\,\E\,G_\ell(\Ipre_\ell)$, where the
indicators $B_{k,\ell}$ in $\Ipre_\ell$ are independent of the state. Let
$\bar\Phi(z)=\sum_\ell\plag_\ell G_\ell(z-\mu\bar M_\ell)$, so that
$\LB=\min_z\bar\Phi(z)=\bar\Phi(z^{\LB})$. By \eqref{eq:ec-excess}, with
$f_\ell\le\sup f_\ell$ and with \eqref{eq:ec-lattice} for the integer
$\Ipre_\ell$,
\[
g^\pi\le\bar\Phi(\E\,z)+\frac{b+h}2\sum_\ell\plag_\ell\,(\sup f_\ell)\,\Var_\pi\Ipre_\ell .
\]
A change of the level changes $\E\,z$ by the same amount and leaves the
variances as they are. With real quantities choose the level so that
$\E\,z=z^{\LB}$. With integer quantities choose it so that
$|\E\,z-z^{\LB}|\le\tfrac12$. The function $\bar\Phi$ is convex and
piecewise linear, its left slope at $z^{\LB}$ is at most zero and its right
slope at least zero, and a closed interval of length $\tfrac12$ contains at
most one kink of each $G_\ell(\cdot-\mu\bar M_\ell)$, of size at most
$(b+h)\sup f_\ell$. So between $z^{\LB}$ and $\E\,z$ the slope of
$\bar\Phi$ is at most $(b+h)\sum_\ell\plag_\ell\sup f_\ell$ in absolute
value, and $\bar\Phi(\E\,z)-\LB\le\tfrac12(b+h)\sum_\ell\plag_\ell\sup f_\ell$.
This is the added one.

Under the base-stock policy, from the second period on, the position is the
level and each order equals the demand of the period before. So
$a_{k,s}B_{k,\ell}$ is a demand times the indicator that the leadtime of its
order is at least $k+\ell+1$, and the terms for different $k$ are
independent. The variance of such a term is
$\P(L\ge k+\ell+1)\,\sigma^2+\P(L\ge k+\ell+1)\P(L\le k+\ell)\,\mu^2$.
\Halmos
\endproof

\proof{Proof of Proposition~\ref{prop:long}.}
Let $\Lmin\ge\max(1,\rng-1)$ throughout.

\textit{No news.} For an age $k\le\rng-1\le\Lmin$, $\P(L\ge k)=1$, so
$\bar a_k=\mu$, the indicators $B_{k,0}$ of \eqref{eq:transy} equal one for
$k\le\rng-2$, and $B_{0,0}=1$ because $\Lmin\ge1$. The transition is
therefore $y'=y+a_{\rng-1}-D$, $a_1'=q$, $a_{k+1}'=a_k$ for every such
$\Lmin$. With $\ell=\Lmin+j$, the probabilities
$m_{k,\ell}=\P(L-\Lmin\ge k+j+1)$ and $\plag_\ell$ depend on $j$ only. Under
a rule that does not change with $\Lmin$, the stationary law of $y-S$, $a$
and the order is the same for every $\Lmin$, and so is the law of
$\Ipre_{\Lmin+j}-S$ for each $j$. This is the first statement, and
for the same reason $\Var^{0}\Ipre$ does not depend on $\Lmin$. Since $\rng\ge2$, each
$m_{k,\Lmin}=\P(L-\Lmin\ge k+1)$ with $1\le k\le\rng-1$ lies strictly between
zero and one, so $\Var^{0}\Ipre>0$. Under the base-stock policy,
Proposition~\ref{prop:ub} and $\bar a_k=\mu$ give
$\Var\Ipre_\ell-\mu^2\sum_km_{k,\ell}(1-m_{k,\ell})=\sigma^2\sum_km_{k,\ell}=\sigma^2\bar M_\ell$.

\textit{Scale.} Let $f_{(n)}$ be the density or the probability function of
$D(n)$, and let $\phi$ be the standard normal density. Under the condition
on the demand the local limit theorem holds,
$\sup_u|\sigma\sqrt n\,f_{(n)}(u)-\phi((u-n\mu)/(\sigma\sqrt n))|\to0$
\citepec[Section~XV.5]{feller1971}. Let $\Lmin\to\infty$. By the central limit
theorem, since $D(\ell+1)$ has variance $(\ell+1)\sigma^2$ and $\cbar$ is
Lipschitz and positively homogeneous, $G_\ell$ at the position
$\mu(\Lmin+1)+\sigma\sqrt{\Lmin}\,\upsilon$ less a bounded amount, divided by
$\sigma\sqrt{\Lmin}$, tends to $\E\,\cbar(\upsilon-N)$ for a standard normal
$N$, uniformly for $\upsilon$ in compact sets and $\ell-\Lmin\le\rng$. The
limit is convex in $\upsilon$, with the minimum $(b+h)\phi_\beta$ at the
$\beta$-quantile of $N$, where $\phi_\beta$ is the value of $\phi$ at that
quantile. Hence $\LB$, $\LB_0$ and $g^*$ are
$(b+h)\phi_\beta\sigma\sqrt{\Lmin}\,(1+o(1))$, the last because $g^*$ lies
between $\LB$ and the cost of the base-stock policy at its best level,
which exceeds $\LB$ by a term of the order $\Lmin^{-1/2}$ by
\eqref{eq:ec-long-up} below. The minimizers $z^{\LB}$ of $\bar\Phi$ and
$z^{0}$ of $\Phi(\cdot,\bar a)$, standardized in the same way, tend
to the quantile. By the local limit theorem,
$\sigma\sqrt{\Lmin}\,f_\ell(u)\to\phi_\beta$ uniformly in $\ell$ and in $u$
within any fixed distance of $z^{\LB}-\mu\bar M_\ell$ or of
$z^{0}-\mu\bar M_\ell$, while
$\sigma\sqrt{\Lmin}\,(\sup f_\ell)\to\phi(0)$.

\textit{The bound.} Set the level as in the proof of
Proposition~\ref{prop:ub}. For $\rho\ge1$,
\eqref{eq:ec-excess}, \eqref{eq:ec-lattice} and the bound on the tail give
\begin{align*}
\E\,G_\ell(\Ipre_\ell)-G_\ell(\E\,\Ipre_\ell)\le{}&\tfrac12(b+h)\max_{|u-\E\Ipre_\ell|\le\rho}f_\ell(u)\,\Var_\pi\Ipre_\ell\\
&+(b+h)\,(\sup f_\ell)\,\E\big[(\Ipre_\ell-\E\,\Ipre_\ell)^2;\,|\Ipre_\ell-\E\,\Ipre_\ell|>\rho\big].
\end{align*}
The law of $\Ipre_\ell-\E\,\Ipre_\ell$ does not depend on $\Lmin$, so the
last expectation tends to zero as $\rho\to\infty$, uniformly in $\Lmin$. The
kinks that enter $\bar\Phi(\E\,z)-\LB$ lie within $\tfrac12$ of
$z^{\LB}-\mu\bar M_\ell$. With the scale step this gives
\begin{equation}\label{eq:ec-long-up}
g^\pi_*-\LB\le\frac{(b+h)\phi_\beta}{2\sigma\sqrt{\Lmin}}\big(\Var_\pi\Ipre+\iota\big)+o\big(\Lmin^{-1/2}\big),
\end{equation}
with $\iota=1$ for integer quantities and $\iota=0$ for real ones. For
$\LB_0$, let $\Ipre^{\,0}_\ell=z^{0}-\sum_k\bar a_kB_{k,\ell}$,
so that $\LB_0=\sum_\ell\plag_\ell\E\,G_\ell(\Ipre^{\,0}_\ell)$,
$\E\,\Ipre^{\,0}_\ell=z^{0}-\mu\bar M_\ell$, and
$\bar\Phi(z^{0})\ge\LB$. The values of $\Ipre^{\,0}_\ell$
lie within $(\rng-1)\mu$ of its mean, and the weight $\varpi$ vanishes
beyond. So \eqref{eq:ec-excess}, \eqref{eq:ec-lattice} and the scale step
give
\[
\LB_0-\LB\ge\frac{(b+h)\phi_\beta}{2\sigma\sqrt{\Lmin}}\big(\Var^{0}\Ipre-\iota'\big)+o\big(\Lmin^{-1/2}\big),
\]
with $\iota'=\tfrac14$ for integer quantities and $\iota'=0$ for real ones.
Subtracting and dividing by $\LB_0$ proves the bound, with
$\iota+\iota'=\tfrac54$. \Halmos
\endproof

\end{document}